\documentclass[11pt]{article}
\usepackage{amssymb}
\usepackage{mathrsfs}
\usepackage{latexsym,amsmath,amssymb,amsfonts,epsfig,graphicx,cite,psfrag}
\usepackage{comment}

\newtheorem{theo}{Theorem}[section]

\newtheorem{lem}[theo]{Lemma}

\def\qed{\hfill \rule{4pt}{7pt}}

\def\pf{\noindent {\it Proof. }}

\begin{document}

\title{The Kelmans-Seymour conjecture III: 3-vertices in $K_4^-$}

\author{Dawei He\footnote{dhe9@math.gatech.edu; Partially supported by NSF grant through X. Yu},
Yan Wang\footnote{ywang416@math.gatech.edu; Partially supported by NSF grant through X. Yu},
Xingxing Yu\footnote{yu@math.gatech.edu; Partially supported by NSF grants DMS-1265564 and DMS-1600738}\bigskip\\
School  of Mathematics\\
Georgia Institute of Technology\\
Atlanta, GA 30332-0160, USA}

\date{}

\maketitle

\begin{abstract}
Let $G$ be a 5-connected nonplanar graph and let $x_1,x_2,y_1,y_2\in V(G)$ be distinct, such that $G[\{x_1,x_2,y_1,y_2\}]\cong K_4^-$ and $y_1y_2\notin E(G)$. 
We show that one of the following holds:  $G-x_1$ contains $K_4^-$, or 
$G$ contains a $K_4^-$ in which $x_1$ is of degree 2, or $G$ contains a $TK_5$ in which $x_1$ is not 
a branch vertex, or $\{x_2,y_1,y_2\}$ may be chosen so that for any 
distinct $z_0, z_1\in N(x_1)-\{x_2,y_1,y_2\}$, $G-\{x_1v:v\notin \{z_0, z_1,x_2,y_1,y_2\}\}$ contains $TK_5$.


\end{abstract}

\newpage 

\section{Introduction}

We use notation and terminology from \cite{HWY15I, HWY15II}. For a graph $K$, we use $TK$ to denote  a \textit{subdivision} of $K$. 
The vertices of $TK$ corresponding to the vertices of $K$ are its \textit{branch vertices}. 
Kelmans~\cite{Ke79} and, independently, Seymour~\cite{Se77} conjectured that every 5-connected nonplanar graph contains $TK_5$. 
In \cite{MY10,MY13}, this conjecture is shown to be true for  graphs containing $K_4^-$.

In \cite{HWY15I} we outline a strategy to prove the Kelmans-Seymour conjecture for graphs not containing  $K_4^-$. 
Let $G$ be a 5-connected nonplanar graph not containing $K_4^-$. Then by a result of Kawarabayashi~\cite{Ka02}, $G$ contains an edge $e$ 
such that $G/e$ is 5-connected. If $G/e$ is planar, 
we can apply a discharging argument. So assume $G/e$ is not planar. Let $M$ be a connected subgraph of $G$ such that $|V(M)|\ge 2$ and 
$G/M$ is 5-connected and nonplanar. 
Let $z$ denote the vertex representing the contraction of $M$, and let $H=G/M$. Then one of the following holds:
\begin{itemize}
\item [(a)] $H$ contains a $K_4^-$ in which  $z$ is of degree 2.
\item [(b)] $H$ contains a $K_4^-$ in which  $z$ is of degree 3.
\item [(c)] $H$ does not contain $K_4^-$, and there exists $T\subseteq H$ such that $z\in V(T)$, $T\cong K_2$ or $T\cong K_3$, and $H/T$ is 5-connected and planar.
\item [(d)] $H$ does not contain  $K_4^-$, and for any $T\subseteq H$ with $z\in V(T)$ and $T\cong K_2$ or $T\cong K_3$, $H/T$ is not 5-connected.
\end{itemize}
Note that local structure around $z$ (in particular, $K_4^-$ containing $z$) will help us find $TK_5$ in $G$ from certain  $TK_5$ in $H$.  

In \cite{HWY15I} we deal with certain special separations and the results can be used to take care of (c). In \cite{HWY15II} we prove results that can be used to 
take care of (a). In this paper, we prove the following, which can be used to take care of (b).

\begin{theo}\label{x1a}
Let $G$ be \textsl{a} $5$-connected nonplanar graph and $x_1,x_2,y_1,y_2 \in V(G)$ be distinct such that 
$G[\{x_1 , x_2 , y_1 , y_2 \}]\cong K_4^-$ and $y_1y_2\notin E(G)$. Then one of the following holds:
\begin{itemize}
\item [$(i)$] $G$ contains \textsl{a} $TK_5$ in which $x_1$ is not \textsl{a}  branch vertex.
\item [$(ii)$] $G-x_1$ contains $K_4^-$, or $G$ contains $K_4^-$ in which $x_1$ is of degree 2. 
\item [$(iii)$] $x_2,y_1,y_2$ may be chosen so that for any distinct $z_0, z_1\in N(x_1)-\{x_2,y_1,y_2\}$, 
$G-\{x_1v:v\notin \{z_0, z_1, x_2,y_1,y_2\}\}$ contains $TK_5$. 
\end{itemize}
\end{theo}

This paper is organized as follows. In Section 2, we list a number of known results that will be used in the proof of Theorem~\ref{x1a}. 
The steps we take to prove  Theorem~\ref{x1a} is quite similar to the arguments in \cite{HWY15II}. First, we find a path in $G$ from $x_1$ to $x_2$ 
such that the graph obtained from $G$ by removing that path satisfies certain connectivity requirements. What is different here is that we need the 
path to include $x_1z_0$ or $x_1z_1$. We find this path in
Section 3, see Figures~\ref{path} and \ref{blocks}. 
In Section 4, we derive further structural information of  the graph $G$. In
Section 5, we find a substructure of $G$ consisting of five additional
paths, see Figure~\ref{structure}. 
In Section 6, we use this substructure to find a $TK_5$ in $G-\{x_1v:v\notin \{z_0, z_1, x_2,y_1,y_2\}\}$. 

\section{Lemmas}

For each positive integer $m$, let $[m]=\{1, \ldots, m\}$.  For convenience, we recall a technical notion from \cite{HWY15I}
(originated from \cite{Se80}). 
A {\it 3-planar graph} $(G,{\cal A})$ consists
of a graph $G$ and  a set ${\cal A}=\{A_1,\ldots,
A_k\}$ of pairwise 
disjoint subsets of $V(G)$ $($possibly ${\cal A}=\emptyset)$ such that
\begin{itemize} 
\addtolength{\baselineskip}{-1ex}

\item for distinct $i, j\in [k]$, $N(A_i)\cap A_j=\emptyset$, 

\item for $i\in [k]$,  $|N(A_i)|\le 3$, and 

\item  if $p(G,{\cal A})$ denotes the 
graph obtained from $G$ by (for each $i\in [k]$) deleting $A_i$ and
adding new edges joining every pair of distinct vertices in  $N(A_i)$,
then $p(G,{\cal A})$ can be drawn in a closed disc in the plane  with no edge
crossing.
\end{itemize}
\noindent If, in addition,  $b_1, \ldots, b_n$ are vertices in $G$ such that $b_i\notin A_j$
for all $i\in [n]$ and $j\in [k]$, $p(G,{\cal A})$ can be drawn in a closed disc in the plane with
no edge crossing, and $b_1,\ldots,b_n$ occur on the boundary of the disc in this cyclic order,
then we say that $(G,{\cal A}, b_1,\ldots,b_n)$ is
{\it 3-planar}. If there is no need to specify ${\cal A}$, we will simply say 
that $(G,b_1,\ldots,b_n)$ is 3-planar.

We can now state the following result of Seymour 
\cite{Se80}; equivalent versions can be found in \cite{CR79, Th80, Sh80}.

\begin{lem}\label{2path} 
Let $G$ be \textsl{a} graph and $s_1,s_2,t_1,t_2$ be distinct vertices of $G$. Then
exactly one of the following holds: 
\begin{itemize}
\item [$(i)$] $G$ contains disjoint paths from $s_1$ to $t_1$ and from
  $s_2$ to $t_2$, respectively. 
\item [$(ii)$] $(G,s_1,s_2,t_1,t_2)$ is 3-planar.
\end{itemize}
\end{lem}

We also state a generalization of Lemma~\ref{2path},  which is a consequence of Theorems 2.3 and 2.4 in \cite{RS90}.

\begin{lem}\label{society}
Let $G$ be \textsl{a} graph,  $v_1,\ldots, v_n\in V(G)$ be distinct, and $n\ge 4$. Then 
exactly one of the following holds:
\begin{itemize}
\item [$(i)$] There exist $1\le i<j<k<l\le n$ such that $G$ contains disjoint paths from $v_i,v_j$ to $v_k, v_l$, respectively. 
\item [$(ii)$] $(G,v_1,v_2,\ldots,v_n)$ is 3-planar.
\end{itemize}
\end{lem}

We will make use of the following result of Perfect \cite{Pe68}. A collection of paths in a graph are said to be {\it independent}
if no internal vertex of any path in this collection belongs to another path in the collection. 

\begin{lem}
\label{perfect}
Let $G$ be \textsl{a} graph, $u\in V(G)$, and $A\subseteq V(G-u)$. Suppose there exist $k$ independent paths from $u$ to distinct $a_1,\ldots, a_k\in A$, respectively, 
and otherwise disjoint from $A$. Then for any $n\ge k$, if there exist $n$
independent paths  $P_1,\ldots, P_n$ in $G$ from $u$ to $n$ distinct vertices in $A$ and otherwise disjoint from $A$ then  $P_1,\ldots, P_n$ may be chosen so that $a_i\in V(P_i)$ for $i\in [k]$. 
\end{lem}

We will also use a result of Watkins and Mesner \cite{WM67} on cycles through three vertices. 
\begin{lem}
\label{Watkins}
Let $G$ be \textsl{a} $2$-connected graph and let $y_1, y_2, y_3$ be three
distinct vertices of $G$. There is no cycle in $G$ containing $\{y_1, y_2, y_3\}$ 
if, and only if, one of the following holds: 
\begin{itemize}
\item [$(i)$] There exists \textsl{a} 2-cut $S$ in $G$ and there exist pairwise disjoint subgraphs $D_{y_i}$ of $G - S$, $i\in [3]$, such that 
$y_i\in V(D_{y_i})$ and each $D_{y_i}$ is  \textsl{a} union of components of $G - S$. 
\item [$(ii)$] There exist 2-cuts $S_{y_i}$ of $G$, $i\in [3]$, and pairwise disjoint subgraphs $D_{y_i}$ of $G$, such that 
$y_i \in V(D_{y_i})$, each $D_{y_i}$ is \textsl{a} union of components of $G-S_{y_i}$, there exists $z\in S_{y_1} \cap S_{y_2} \cap S_{y_3}$, 
and $S_{y_1} - \{z\}, S_{y_2} - \{z\}, S_{y_3} - \{z\}$ are pairwise disjoint. 
\item [$(iii)$] There exist pairwise disjoint $2$-cuts $S_{y_i}$ in $G$, $i\in [3]$, and pairwise disjoint subgraphs $D_{y_i}$ of 
$G - S_{y_i}$ such that $y_i \in V(D_{y_i})$, 
$D_{y_i}$ is \textsl{a} union of components of $G - S_{y_i}$, and $G - V(D_{y_1} \cup D_{y_2} \cup D_{y_3})$ has precisely two components, 
each containing exactly one vertex from $S_{y_i}$ for $i\in [3]$.
\end{itemize}
\end{lem}

The next result is Theorem 3.2  from \cite{MY10}.
\begin{lem}\label{2-connected}
Let $G$ be \textsl{a} $5$-connected nonplanar graph and let $x_1 , x_2
, y_1 , y_2\in V(G)$ be 
distinct such that $G[\{x_1 , x_2 , y_1 , y_2\}]\cong K_4^-$ and  $y_1y_2\notin E(G)$. 
Suppose $G-x_1x_2$ contains \textsl{a} path $X$ between $x_1$ and $x_2$ 
such that $G-X$ is 2-connected, $X-x_2$ is induced in $G$, and $y_1,y_2\notin V(X)$. 
Let $v\in V(X)$ such that $x_2v\in E(X)$. Then $G$ contains \textsl{a} $TK_5$ in which $x_2v$ is an edge and 
$x_1,x_2,y_1,y_2$ are branch vertices. 
\end{lem}

It is easy to see that under the conditions of
Lemma~\ref{2-connected}, $G-\{x_2u:u\notin \{v,x_1,y_1,y_2\}\}$
contains $TK_5$. The next result is Corollary 2.11 in \cite{KMY15}. For a graph $G$ and $A\subseteq V(G)$, we say that $(G,A)$ is {\it plane} if $G$ is drawn 
in the plane with no edge crossings, and the vertices in $A$ are incident with the outer face of $G$; and we say that  $(G,A)$ is {\it planar}
if $G$ admits such a planar drawing.

\begin{lem}\label{apex1}
Let $G$ be \textsl{a} connected graph with $|V(G)|\ge 7$, $A\subseteq V(G)$ with $|A|=5$, and $a\in A$, such that  
$G$ is $(5,A)$-connected, $(G-a,A-\{a\})$ is plane, and
$G$ has no 5-separation $(G_1,G_2)$ with $A\subseteq G_1$
and $|V(G_2)|\ge 7$. Suppose there exists $w\in N(a)$ such that $w$ is not incident with  the outer face of $G-a$.
Then
\begin{itemize}
\item [$(i)$] the vertices of $G-a$ cofacial with $w$ induce \textsl{a} cycle $C_w$ in $G-a$, and
\item [$(ii)$] $G-a$ contains paths $P_1,P_2, P_3$ from $w$ to $A-\{a\}$  such that $V(P_i\cap P_j)=\{w\}$
for $1\le i<j\le 3$, and $|V(P_i\cap C_w)|=|V(P_i)\cap A|=1$ for $i\in [3]$.
\end{itemize}
\end{lem}

  The next four results are Theorem 1.1, Theorem 1.2, Proposition 2.3 and Proposition 4.2, respectively,  in \cite{HWY15I}.
Note that condition $(iii)$ in three of these four results (Theorem 1.1, Theorem 1.2 and  Proposition 4.2 in \cite{HWY15I}) states that 
 $G$ has a $5$-separation $(G_1',G_2')$ such that $V(G_1'\cap G_2')=\{a, a_1,a_2,a_3,a_4\}$ and $G_2'$ is 
the graph obtained from the edge-disjoint union of the $8$-cycle $a_1b_1a_2b_2a_3b_3a_4b_4\allowbreak a_1$ 
and the  $4$-cycle $b_1b_2b_3b_4b_1$ by adding $a$ and the edges $ab_i$ for $i\in [4]$.
This condition implies that $G$ contains a $K_4^-$ in which $a$ is of
degree 2. We only need  the weaker versions of these results.

\begin{lem}\label{apex}
Let $G$ be \textsl{a} $5$-connected nonplanar graph and let $(G_1, G_2)$ be \textsl{a} $5$-separation in $G$. 
Suppose $|V(G_i)|\geq 7$ for $i\in [2]$, $a\in V(G_1\cap G_2)$, and $(G_2-a,V(G_1\cap G_2)-\{a\})$ is planar. Then one of the following holds:
\begin{itemize}
\item [$(i)$]  $G$ contains \textsl{a} $TK_5$ in which $a$ is not \textsl{a} branch vertex.
\item [$(ii)$] $G-a$ contains $K_4^-$, or $G$ contains \textsl{a} $K_4^-$ in which $a$ is of degree 2.

\end{itemize}
\end{lem}

\begin{lem}
\label{5cut_triangle}
Let $G$ be \textsl{a} $5$-connected graph and $(G_1,G_2)$ be \textsl{a} $5$-separation in $G$. Suppose that $|V(G_i)|\ge 7$ for $i\in [2]$ and 
$G[V(G_1\cap G_2)]$ contains \textsl{a} triangle $aa_1a_2a$. Then one of the following holds:
\begin{itemize}
\item [$(i)$]  $G$ contains \textsl{a} $TK_5$ in which $a$ is not \textsl{a} branch vertex.
\item [$(ii)$] $G-a$ contains $K_4^-$, or $G$ contains \textsl{a} $K_4^-$ in which $a$ is of degree 2.
\item [$(iii)$] For any distinct $u_1,u_2,u_3\in N(a)-\{a_1,a_2\}$, $G-\{av: v \not\in \{a_1,a_2,u_1,u_2,u_3\}\}$ contains $TK_5$. 
\end{itemize}
\end{lem}

\begin{lem}\label{6-cut2}
Let $G$ be \textsl{a} graph, $A\subseteq V(G)$, and $a\in A$  such that $|A|=6$, $|V(G)|\geq 8$, 
$(G-a, A-\{a\})$ is planar, and $G$ is $(5, A)$-connected.
Then one of the following holds: 
\begin{itemize}
 \item[$(i)$] $G-a$ contains $K_4^-$, or $G$ contains \textsl{a} $K_4^-$ in which the degree of $a$ is $2$.
 \item[$(ii)$] $G$ has \textsl{a} $5$-separation $(G_1,G_2)$ such that $a\in V(G_1\cap G_2)$, $A\subseteq V(G_1)$, $|V(G_2)|\ge 7$, and 
$(G_2-a, V(G_1\cap G_2)-\{a\})$ is planar. 
\end{itemize}
\end{lem}

\begin{lem} \label{apexvertex}
Let $G$ be \textsl{a} $5$-connected nonplanar graph and $a\in V(G)$ such that $G-a$ is planar. Then one of the following holds:
\begin{itemize}
 \item[$(i)$] $G$ contains \textsl{a} $TK_5$ in which $a$ is not \textsl{a} branch vertex.
 \item[$(ii)$]  $G-a$ contains $K_4^-$, or $G$ contains \textsl{a} $K_4^-$ in which $a$ is of degree 2.
\end{itemize}
\end{lem}

We also need Lemma 3.1 in \cite{HWY15II}.  Let $G$ be a graph and $\{u, v\}\subseteq V(G)$. 
We say that a sequence of blocks $B_1, \ldots, B_k$ in $G$ 
is {\it a chain of blocks} from $u$ to $v$ if $|V(B_i)\cap
V(B_{i+1})|=1$ for $i\in [k-1]$,  $V(B_i)\cap V(B_j)=\emptyset$ for
any $1 \leq i < i+1 < j\leq k$, $u,v\in V(B_1)$ are distinct when
$k=1$, and $u\in V(B_1)-V(B_2)$ and $v\in V(B_k)-V(B_{k-1})$ when
$k\ge 2$. A block is {\it nontrivial} if it is 2-connected.

\begin{lem}\label{4A} 
Let $G$ be \textsl{a} graph and $A=\{x_1, x_2, y_1, y_2\}\subseteq V(G)$ such that $G$ is $(4, A)$-connected. Suppose 
there exists \textsl{a} path $X$ in $G-x_1x_2$ from $x_1$ to $x_2$ such that $G-X$ contains \textsl{a} chain of blocks $B$ from $y_1$ to $y_2$. 
Then one of the following holds:
\begin{itemize}
 \item[$(i)$] There is \textsl{a} $4$-separation $(G_1, G_2)$ in $G$ such that $B+\{x_1,x_2\}\subseteq G_1$, $|V(G_2)|\geq 6$, and $(G_2, V(G_1\cap G_2))$ is planar.
 \item[$(ii)$] There exists an induced path $X'$ in $G-x_1x_2$ from $x_1$ to $x_2$ such that $G-X'$ is \textsl{a} 
chain of blocks from $y_1$ to $y_2$ and contains $B$.
\end{itemize}
\end{lem}

\section{Nonseparating paths}

Let $G$ be a $5$-connected nonplanar graph and $x_1,x_2,y_1,y_2 \in V(G)$ be distinct such that 
$G[\{x_1 , x_2 , y_1 , y_2 \}]\cong K_4^-$ and $y_1y_2\notin E(G)$. To take care of case (b) described  in Section 1, 
we need to find a path in $G$  satisfying certain  properties (see $(iv)$ of Lemma~\ref{classify2}).  As a first step, we prove the following.

\begin{lem}\label{classify1} 
Let $G$ be \textsl{a} $5$-connected nonplanar graph and $x_1,x_2,y_1,y_2 \in V(G)$ be distinct such that 
$G[\{x_1 , x_2 , y_1 , y_2 \}]\cong K_4^-$ and $y_1y_2\notin E(G)$. Let  $z_0, z_1\in N(x_1)-\{x_2,y_1,y_2\}$ be distinct. 
Then one of the following  holds:
\begin{itemize}
 \item [$(i)$] $G$ contains \textsl{a} $TK_5$ in which $x_1$ is not a branch vertex. 
\item [$(ii)$] $G-x_1$ contains $K_4^-$, or $G$ contains \textsl{a} $K_4^-$ in which $x_1$ is of degree 2.
\item [$(iii)$] There exist $i\in \{0, 1\}$ and an induced path $X$ in $G-x_1$ from $z_i$ to $x_2$ such that $(G-x_1)-X$ is a chain of blocks 
from $y_1$ to $y_2$, $z_{1-i}\notin V(X)$, and one of $y_1, y_2$ is contained in a nontrivial block of $(G-x_1)-X$.
\end{itemize}
\end{lem} 

\pf We may assume  $G-x_1$ contains disjoint paths $X, Y$ from $z_1, y_1$
to $x_2 , y_2$, respectively. For, otherwise, since $G$ is 5-connected, it follows from Lemma~\ref{2path} that  $(G-x_1,z_1,y_1,x_2,y_2)$ 
is planar; so $(i)$ or $(ii)$ holds by Lemma~\ref{apexvertex}. 

Hence $(G-x_1)-X$ contains a chain of blocks from $y_1$ to $y_2$, say $B$. We may assume that $(G-x_1)-X$ is a chain of blocks from $y_1$ to $y_2$. 
For otherwise, we may apply Lemma~\ref{4A} to conclude that 
$G$ has a $5$-separation $(G_1, G_2)$ such that $x_1\in V(G_1\cap
G_2)$, $B+\{x_1,x_2,z_1\}\subseteq G_1$, $|V(G_2)|\geq 7$, and
$(G_2-x_1, V(G_1\cap G_2)-\{x_1\})$ is planar. If $|V(G_1)|\ge 7$ then 
 $(i)$ or $(ii)$  follows from Lemma~\ref{apex}. So assume
$|V(G_1)|\le 6$. Since $y_1y_2\notin E(G)$, $|V(G_1)|=6$ and
$|V(B)|=3$.  Let $V(B)=\{y_1,y_2,v\}$. Since $G$ is 5-connected and
$y_1y_2\notin E(G)$, $\{x_1,x_2,y_1,y_2,z_1\}= V(G_1\cap G_2)=N(v)$. Hence,
$G[\{v,x_1,x_2,y_1\}]-x_1x_2$ is a $K_4^-$ in which $x_1$ is of degree
2, and $(ii)$ holds. 

We may further assume that $z_0\notin V(X)$. For, suppose $z_0\in V(X)$. Since $G$ is 5-connected and $X$ is induced in 
$G-x_1$, every vertex of $X$ has at least two neighbors in $(G-x_1)-X$. Hence,  
 $(G-x_1)-z_0Xx_2$ is also a chain of  blocks from $y_1$ to $y_2$. So   we may use  $z_0Xx_2$ as $X$.

Let $B_1 , B_2$ be the blocks in $(G-x_1)-X$ containing $y_1, y_2$, respectively. If one of $B_1, B_2$ is nontrivial, then
$(iii)$ holds. So we may assume that $|V(B_1)|=|V(B_2)|=2$. Since $X$
is induced and $G$ is 5-connected, there exists $z\in N(x_2)-(\{x_1,y_1,y_2\}\cup V(X))$, and 
$y_1$ and $y_2$ each have at least two neighbors on $X-x_2$. 
Let $Z$ be a path in $(G-x_1)-X-\{y_1, y_2\}$ from $z_0$ to $z$. 
Then  $(G-x_1)-Z$ contains a chain of blocks, say $B$,  from $y_1$ to $y_2$, and the blocks in $(G-x_1)-Z$ containing $y_1$ or $y_2$ are nontrivial. 
Thus, we may apply  Lemma~\ref{4A} to $G$, $Z$ and $B$. If $(ii)$ of Lemma~\ref{4A} holds,
we have $(iii)$. So assume   $(i)$ of Lemma~\ref{4A} holds. Then, as
in the second paragraph of this proof, 
$(i)$ or $(ii)$ follows from Lemma~\ref{apex}. \qed

\medskip

\begin{figure}
\begin{center}
\includegraphics[scale=0.25]{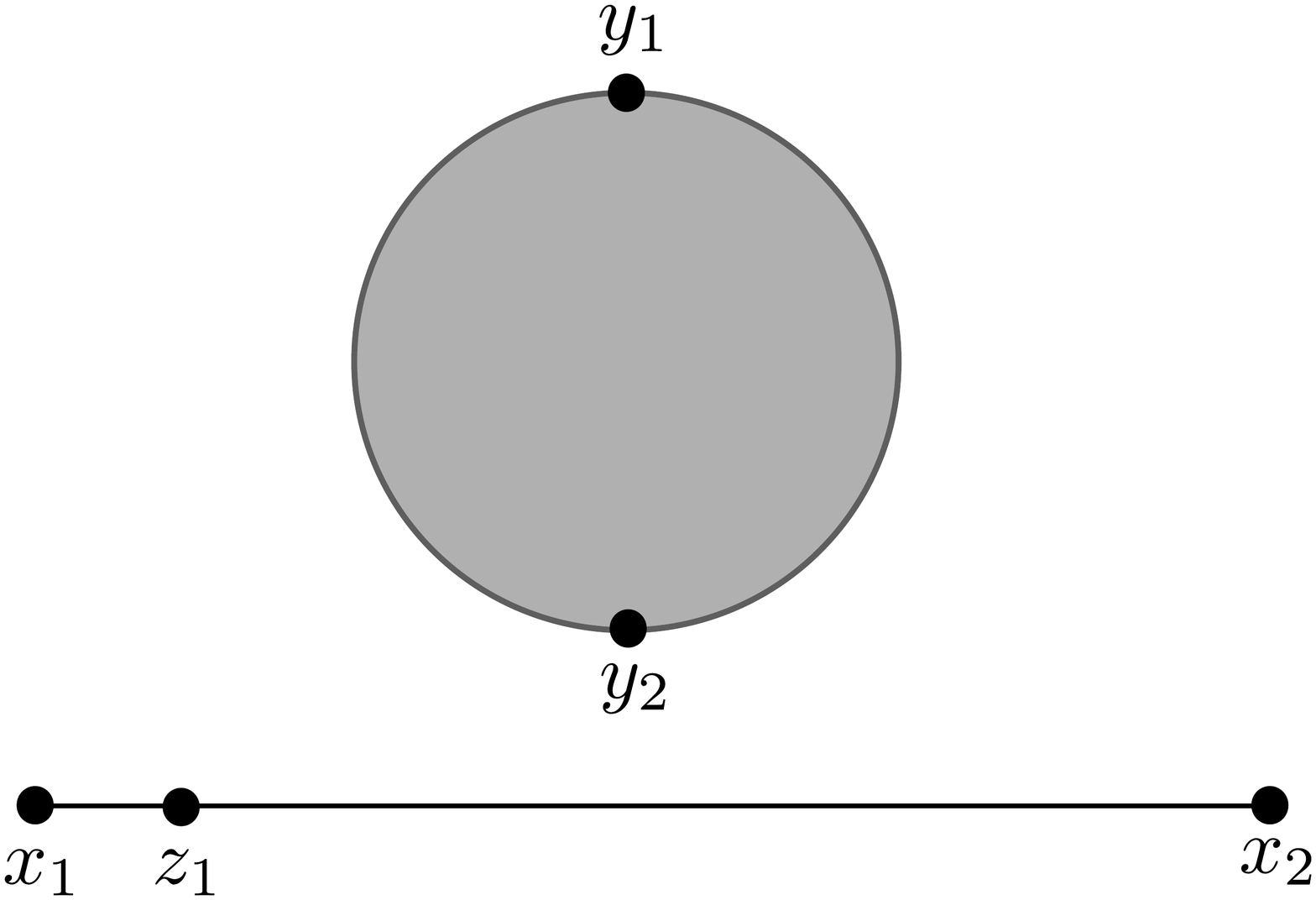}
\includegraphics[scale=0.25]{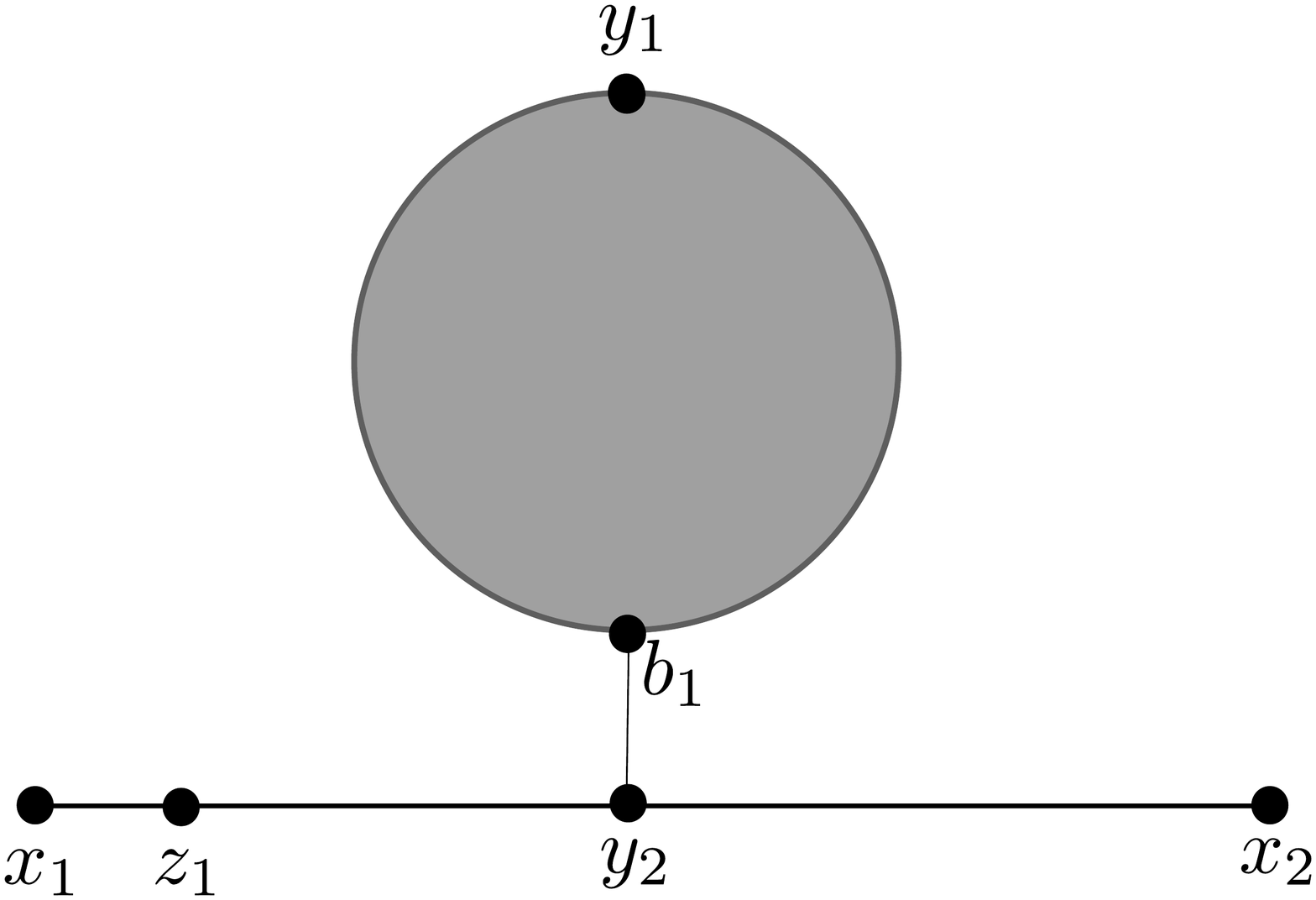}
\caption{\label{path} Structure of $G$ in $(iii)$ of Lemma~\ref{classify2}. }
\end{center}
\end{figure}

\begin{figure}
\begin{center}
\includegraphics[scale=0.25]{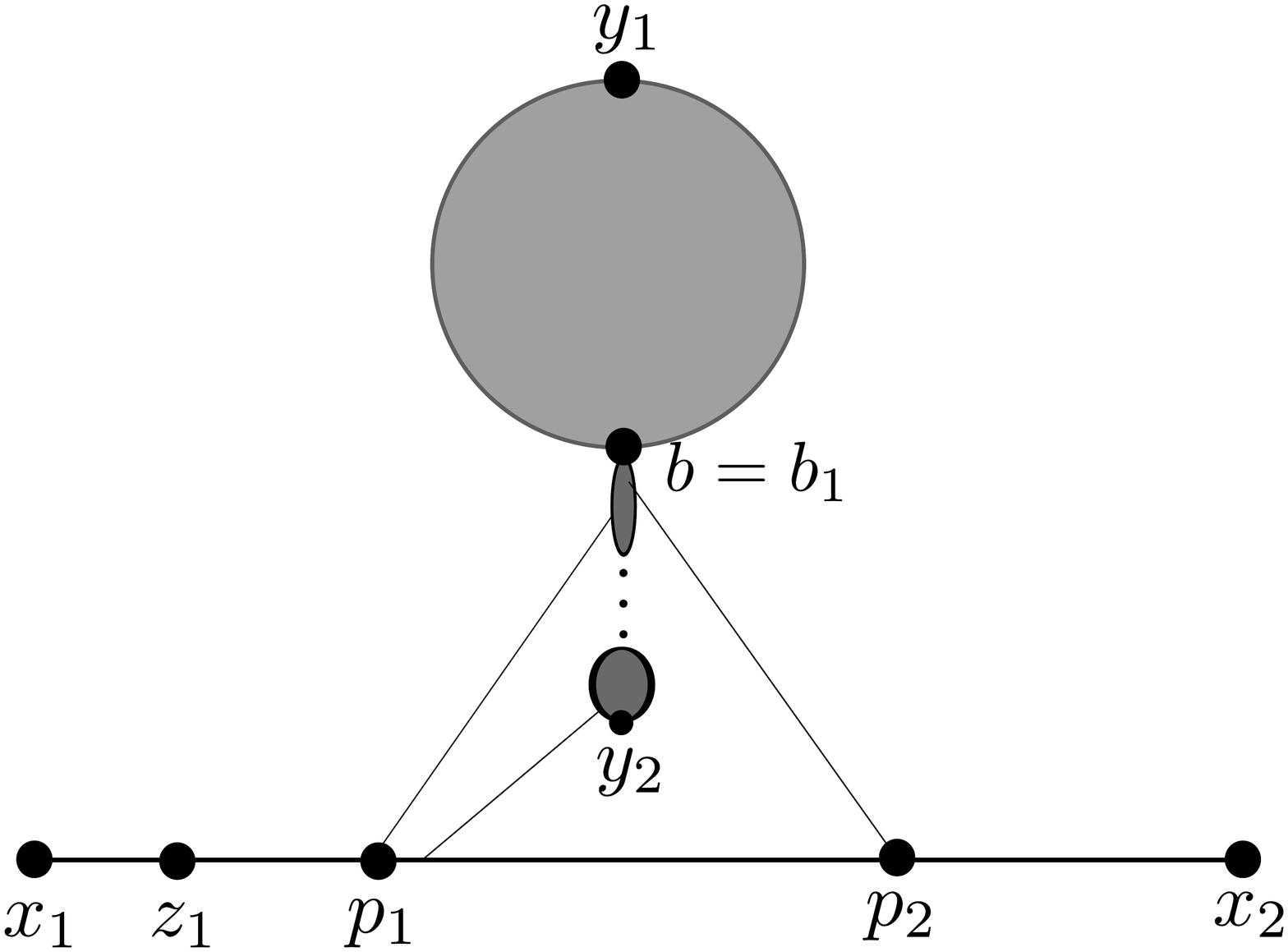}
\includegraphics[scale=0.25]{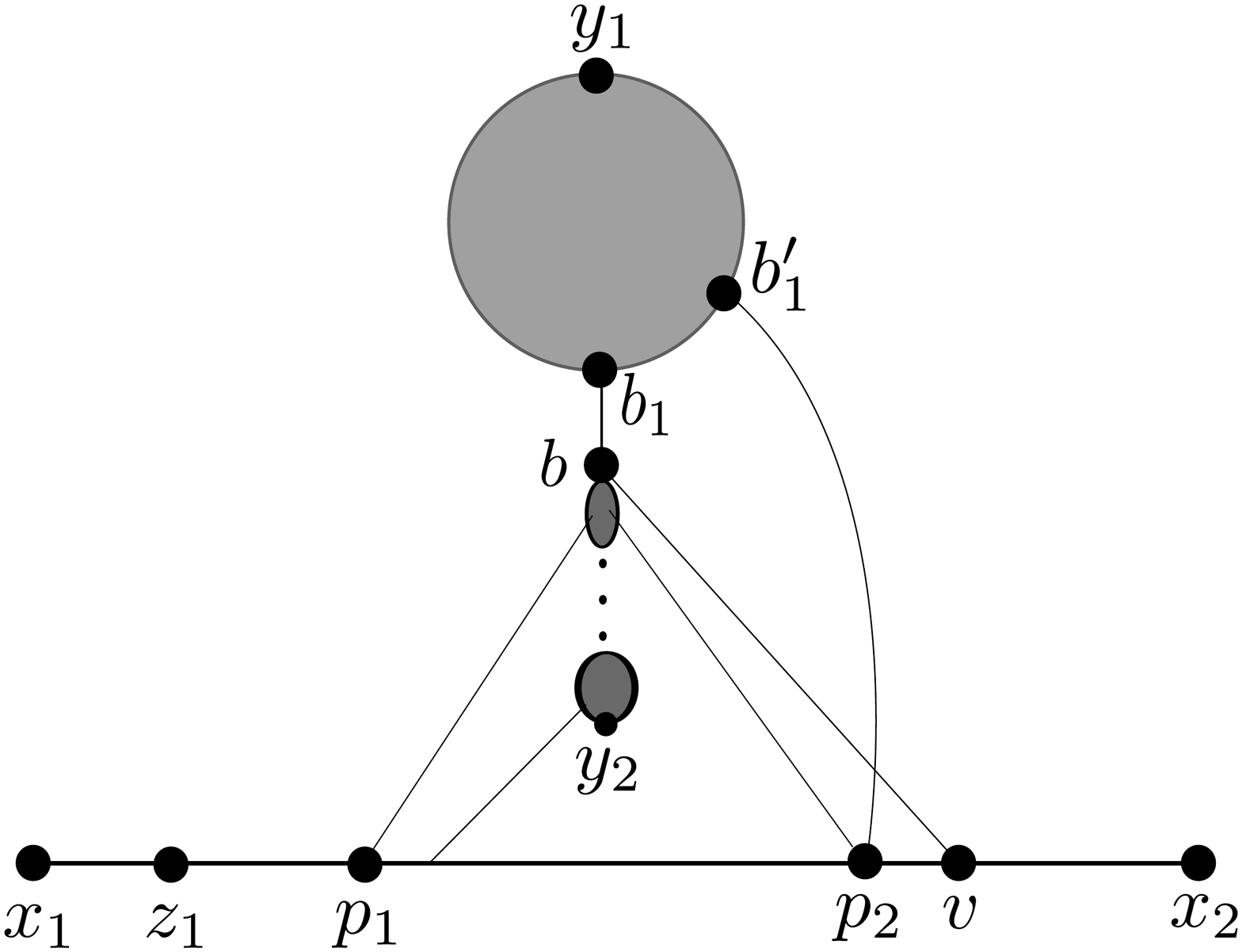}
\caption{\label{blocks} Structure of $G$ in $(iv)$ of
  Lemma~\ref{classify2} (with $j=1$). }
\end{center}
\end{figure}

We have results from \cite{HWY15I, HWY15II, MY13} that can be used to deal with $(i)$ or $(ii)$ of Lemma~\ref{classify1}. 
In this paper, we deal with $(iii)$ of Lemma~\ref{classify1}.  Parts
$(iii)$ and $(iv)$ of the next lemma give more detailed structure of $G$ 
when $(iii)$ of Lemma~\ref{classify1} occurs. We refer the reader to
Figure~\ref{path} for $(iii)$ of Lemma~\ref{classify2}, and
Figure~\ref{blocks} for $(iv)$ of Lemma~\ref{classify2}.

For a graph $H$ and a subgraph $L$ of $H$, an $L$-bridge of $H$ is a subgraph of $H$ 
that is induced by an edge in $E(H)-E(L)$ with both incident vertices in $V(L)$, or is induced by the edges in a 
component of $H-L$ as well as edges from that component to $L$.

\begin{lem}\label{classify2}
Let $G$ be \textsl{a} $5$-connected nonplanar graph and $x_1,x_2,y_1,y_2 \in V(G)$ be distinct such that 
$G[\{x_1 , x_2 , y_1 , y_2 \}]\cong K_4^-$ and $y_1y_2\notin E(G)$. Let  $z_0, z_1\in N(x_1)-\{x_2,y_1,y_2\}$ be distinct and let $G':=
G-\{x_1x: x\notin \{x_2,y_1,y_2,z_0, z_1\}\}$. Then one of the following  holds:
\begin{itemize}
 \item [$(i)$] $G'$ contains $TK_5$, or $G$ contains \textsl{a} $TK_5$ in which $x_1$ is not  \textsl{a} branch vertex. 
\item [$(ii)$] $G-x_1$ contains $K_4^-$, or $G$ contains  \textsl{a} $K_4^-$ in which $x_1$ is of degree 2. 
\item [$(iii)$] The notation of $y_1,y_2,z_0, z_1$ 
may be chosen so that $(G-x_1)-x_2y_2$ has an induced path $X$ from $z_1$ to $x_2$ such that $z_0, y_1\notin V(X)$,  
and $(G-x_1)-X$ is 2-connected. 
\item [$(iv)$]  The notation of $z_0, z_1$ may be chosen so that there exists an induced path $X$ in $G-x_1$ from $z_1$ to $x_2$ such that 
$z_0\notin V(X)$,  $(G-x_1)-X$ is a chain of blocks $B_1, \ldots, B_k$ from $y_1$ to $y_2$ with $B_1$ nontrivial, 
$z_0\in V(B_1)$ when $z_1$ has at least two neighbors in $B_1$, and $(G-x_1)-x_2y_2$ has a 3-separation $(Y_1, Y_2)$ 
such that $V(Y_1\cap Y_2)=\{b, p_1, p_2\}$, $z_1, p_1, p_2, x_2$ occur on $X$ in this order, 
$Y_1=G[B_1\cup z_1Xp_1\cup p_2Xx_2+b]$, $p_1Xp_2+y_2\subseteq Y_2$, and $p_1, p_2$ each have at least two neighbors in $Y_2-B_1$. 
Moreover, if  $b\notin V(B_1)$ then  $V(B_2)=\{b_1,b\}$ with $b_1\in V(B_1)$, and 
there exists some $j\in [2]$ such that $p_{3-j}$ has a unique neighbor $b_1'$ in $B_1$, 
$b$ has  \textsl{a} unique neighbor $v$ in $X-p_1Xp_2$ such that $vp_{3-j}\in E(X)-E(p_1Xp_2)$, $vb_1\notin E(G)$ and $p_jb\notin E(G)$.
\end{itemize}
\end{lem}

\pf We begin our proof by applying Lemma~\ref{classify1} to
$G,x_1,x_2,y_1,y_2$. If $(i)$ or $(ii)$ of Lemma~\ref{classify1} holds
then assertion  $(i)$ or $(ii)$ of this lemma holds. 
 So we may assume that $(iii)$ of Lemma~\ref{classify1} holds. Then we may assume $(G-x_1)-x_2y_2$ has an induced path $X$ from $z_1$ to $x_2$ such that 
$z_0,y_1\notin V(X)$, $(G-x_1)-X$ has a nontrivial block $B_1$ containing $y_1$, and $y_1$ is not a cut vertex of $(G-x_1)-X$. 
(Note that we are not requiring the stronger condition that  $(G-x_1)-X$ be a chain of blocks from $y_1$ to $y_2$.) 
We choose such a path $X$ that 
\begin{itemize}
 \item[(1)] $B_1$ is maximal,
\item [(2)] subject to (1), whenever possible, $(G-x_1)-X$ has a chain of blocks from $y_1$ to $y_2$ and containing $B_1$, and 
 \item[(3)] subject to (2), the component $H$ of $(G-x_1)-X$
   containing $B_1$ is maximal.
\end{itemize}

Let ${\cal C}$ be the set of all components of $(G-x_1)-X$ different
from 
$H$. Then

\begin{itemize}
\item [(4)] ${\cal C} = \emptyset$, $H=(G-x_1)-X$, and if $y_2\notin V(X)$ then  $H$ is a chain of blocks from $y_1$ to $y_2$ and containing $B_1$. 
\end{itemize}
First, suppose ${\cal C} = \emptyset$. Then  $H=(G-x_1)-X$. 
Suppose  $y_2\notin V(X)$. Then $H$ has a chain of blocks, say $B$, from $y_1$
to $y_2$ and containing $B_1$. By applying  Lemma~\ref{4A} to $G-x_1,z_1,x_2,y_1,y_2$, $X$ can be chosen so that $(G-x_1)-X$ is a chain of blocks from $y_1$ 
to $y_2$, or 
$G$ has a 5-separation $(G_1,G_2)$ such that $x_1\in V(G_1\cap G_2)$,
$B+\{x_1,x_2,z_1\}\subseteq G_1$, $|V(G_2)|\ge 7$ and $(G_2-x_1,V(G_1\cap
G_2)-\{x_1\})$ is planar. We may assume the latter as otherwise (4) holds. Since $B_1$ is nontrivial and 
$y_1y_2\notin E(G)$, $|V(B)|\ge 4$. 
%So $|V(G_1)|\ge 6$. If $|V(G_1)|=6$ then, since $y_1y_2\notin E(G)$ and $G$ is 5-connected,
%$y_1,y_2,z_1\in V(G_1\cap G_2)$ and there exists $v\in V(G_1)-V(G_2)$ such that $N(v)=V(G_1\cap G_2)$; now
%$G[\{v,x_1,x_2,y_1\}]-x_1y_1$ is a $K_4^-$ in which $x_1$ is of degree 2, and $(ii)$ holds. 
So $|V(G_1)|\ge 7$; and  $(i)$ or $(ii)$ follows from Lemma~\ref{apex}.

Now suppose ${\cal C}\ne \emptyset$. For each $D\in {\cal C}$, let
$u_D,v_D\in V(X)$ be the neighbors of $D$ 
in $G-x_2y_2$ with $u_DXv_D$ maximal, and assume  that 
$z_1,u_D,v_D,x_2$ occur on $X$ in this order. 
 Define a new graph $G_{\cal C}$ such that $V(G_{\cal C}) = {\cal C}$, and 
two components $C, D\in {\cal C}$ are adjacent in $G_{\cal C}$ if $u_CXv_C-\{u_C,v_C\}$ 
contains a neighbor of $D$ or  $u_DXv_D-\{u_D,v_D\}$ contains a neighbor of $C$. 

Note that, for any component ${\cal D}$ of $G_{\cal C}$, $\bigcup_{D\in V({\cal D})}u_DXv_D$ is a subpath of $X$. Since $G$ is $5$-connected, 
there exist $y\in V(H)$ and $C\in V({\cal D})$ with $N(y)\cap V(u_CXv_C - \{u_C, v_C\})\ne \emptyset$. 

If $y\ne y_1$ then let $Q$ be an induced path in $G[C+\{u_C,
v_C\}]-x_2y_2$ from $u_C$ to $v_C$, and let $X'$ be obtained from $X$
by replacing $u_CXv_C$ with $Q$. Then $B_1$ is contained in a block of
$(G-x_1)-X'$, and $y_1$ is not a cut vertex of $(G-x_1)-X'$. Moreover,
if $(G-x_1)-X$ has a chain of blocks from $y_1$ to $y_2$ then so does
$(G-x_1)-X'$. However, the component of $(G-x_1)-X'$ containing $B_1$
is larger than $H$, contradicting (3). 

So we may assume that $y = y_1$ for all choices of $y$ and $C$. Let $uXv : = \bigcup_{D\in V({\cal D})} u_DXv_D$. Since $G$ is $5$-connected, 
$y_2\in V(\bigcup_{D\in V({\cal D})} D)\cup V(uXv-\{u, v\})$ and $G$ has a separation $(G_1, G_2)$ such that $V(G_1\cap G_2)=\{u,v,x_1, x_2, y_1\}$,
$G_1:=G[\bigcup_{D\in V({\cal D})} D\cup uXv +\{x_1, x_2, y_1\}]$, and $B_1\cup z_1Xu\cup vXx_2\subseteq G_2$. Clearly, $|V(G_i)|\ge 7$ for $i\in [2]$. 
Since $G[\{x_1, x_2, y_1\}]\cong K_3$,  $(i)$ or $(ii)$ or $(iii)$ follows from
Lemma~\ref{5cut_triangle}. This completes the proof of (4).

\medskip

Let ${\cal B}$ be the set of all $B_1$-bridges of $H$. For each $D\in
{\cal B}$, let $b_D\in V(D)\cap V(B_1)$ and $u_D, v_D \in V(X)$ 
be the neighbors of $D-b_D$ in $G-x_2y_2$ with $u_DXv_D$ maximal. Define a new graph $G_{\cal B}$ such that $V(G_{\cal B}) = {\cal B}$, and 
two $B_1$-bridges $C, D\in {\cal B}$ are adjacent in $G_{\cal B}$ if $u_CXv_C-\{u_C,v_C\}$ 
contains a neighbor of $D-b_D$ or  $u_DXv_D-\{u_D,v_D\}$ contains a neighbor of $C-b_C$. 
Note that, for any component ${\cal D}$ of $G_{\cal B}$,
$\bigcup_{D\in V({\cal D})}u_DXv_D$ is a subpath of $X$, whose ends
are denoted by $u_{{\cal D}}, v_{{\cal D}}$. 
We let $S_{{\cal D}}:= \{b_D: D\in V({\cal D})\}\cup (N(u_{{\cal D}}Xv_{{\cal D}}-\{u_{{\cal D}}, v_{{\cal D}}\})\cap V(B_1))$.
We may assume that 

\begin{itemize}
\item [(5)] for any component ${\cal D}$ of $G_{\cal B}$, $|S_{\cal D}|\le 2$ and 
$y_2\in \left(\bigcup_{D\in V({\cal D})}V(D)-S_{{\cal D}}\right)\cup V(u_{{\cal D}}Xv_{{\cal D}}-\{u_{{\cal D}}, v_{{\cal D}}\})$. 
\end{itemize}
First,  we may assume $|S_{\cal D}|\le 2$. For, suppose $|S_{\cal D}|\ge 3$. 
Then there exist $D\in V({\cal D})$, $r_1,r_2\in V(u_DXv_D)-\{u_D,v_D\}$, and distinct $r_1',r_2'\in V(B_1)$ such that 
for $i\in [2]$, $r_ir_i'\in E(G)$ or $r_i'=B_{D_i}$ for some $D_i\in
V({\cal D})-\{D\}$. (To see this, we choose $D\in V({\cal D})$ such
that there is a maximum  number of vertices in $B_1$ from which $G$ has
a path to $u_DXv_D-\{u_D,v_D\}$ and internally disjoint from $B_1\cup
D\cup X$. If this number is at most $1$, we can show that $|S_{\cal D}|\le 2$.)
Let $R_i=r_ir_i'$ if $r_ir_i'\in E(G)$; and otherwise let $R_i$ be a path in $G[D_i+r_i]$ 
from $r_i$ to $r_i'$ and internally disjoint from $X$. 
Let $Q$ denote an induced path in $G[D+\{u_{D},v_{D}\}]-b_D-x_2y_2$ 
between $u_{D}$ and $v_{D}$, and let $X'$ be obtained from $X$ by
replacing $u_{D}Xv_{D}$ with $Q$. Clearly, 
the block of $(G-x_1)-X'$ containing $y_1$  contains $B_1$ as well as the path  $R_1\cup r_1Xr_2\cup R_2$. 
Note that $y_1\ne b_{D}$ (as $y_1$ is not a cut vertex in $H$). 
Moreover, if $y_1=r_i'$ for some $i\in [2]$ then $D_i$ is not defined and $r_ir_i'\in E(G)$. So $y_1$ is not a cut vertex of $(G-x_1)-X'$.
Thus, $X'$ contradicts the choice of $X$,  because of (1). 

Now assume $y_2\notin \bigcup_{D\in V({\cal D})}V(D)\cup V(u_{{\cal D}}Xv_{{\cal D}})-(\{u_{{\cal D}}, v_{{\cal D}}\}\cup S_{{\cal D}})$. Then $S_{\cal D}
\cup \{u_{{\cal D}},v_{{\cal D}},x_1\}$ is a cut in $G$; so $|S_{\cal
  D}|=2$ (as $G$ is 5-connected). Let $S_{\cal D}=\{p,q\}$. 
Then $G$ has a 5-separation $(G_1,G_2)$ such that $V(G_1\cap G_2)=\{p,q,u_{{\cal D}},v_{{\cal D}},x_1\}$,  
$B_1\cup z_1Xu_{{\cal D}}\cup v_{{\cal D}}Xx_2\subseteq G_1$, and $G_2$ contains $u_{{\cal D}}Xv_{{\cal D}}$ and the $B_1$-bridges of $H$ contained in ${\cal D}$. 
If $(G_2-x_1,u_{{\cal D}},p,v_{{\cal D}},q)$ is planar then, since
$|V(G_i)|\ge 7$ for $i\in [2]$, the assertion of this lemma  follows from Lemma~\ref{apex}. 
So we may assume that $(G_2-x_1,u_{{\cal D}},p,v_{{\cal D}},q)$ is not
planar. Then by Lemma~\ref{2path}, $G_2-x_1$ contains disjoint paths
$S,T$ from $u_{{\cal D}},p$ to $v_{{\cal D}},q$, respectively. 

We  apply Lemma~\ref{4A} to $G_2-x_1$ and $\{u_{{\cal D}},v_{{\cal D}},p,q\}$. If $(i)$ of Lemma~\ref{4A}
holds then from the separation in $G_2-x_1$, we derive a 5-separation
$(G_1',G_2')$ in $G$ such that $x_1\in V(G_1'\cap G_2')$, $B_1\cup
T+x_1\subseteq G_1'$, $|V(G_2')|\ge 7$, and $(G_2'-x_1,V(G_1'\cap
G_2')-\{x_1\})$ is planar. So  $(i)$ or $(ii)$ follows from Lemma~\ref{apex}.  
We may thus assume that $(ii)$ of Lemma~\ref{4A} holds.  
Thus, there is an induced path $S'$ in $G_2-x_1$ from $u_{{\cal D}}$ to $v_{{\cal D}}$ such that $(G_2-x_1)-S'$ is a 
chain of blocks from $p$ to $q$. 
Now let $X'$ be obtained from $X$ by replacing $u_{{\cal D}}Xv_{{\cal D}}$ with $S'$. Then $y_1$ is not a cut vertex of $(G-x_1)-X'$, and 
the block of $(G-x_1)-X'$ containing $y_1$ contains $B_1$ and
$(G_2-x_1)-S'$, contradicting (1). This completes the proof of  (5).

\medskip

We may also assume that

\begin{itemize}
\item[(6)] for any $B_1$-bridge $D$ of $H$, $y_2\notin V(u_DXv_D)-\{u_D,v_D\}$. 
\end{itemize}
For, suppose  $y_2\in V(u_DXv_D)-\{u_D,v_D\}$ for some $B_1$-bridge $D$ of $H$. Choose $X$ and $D$ so that, subject to (1)-(3), $u_DXv_D$ is maximal. 

We claim that $\{D\}$ is a component of $G_{\cal B}$. For, otherwise,
by the maximality of $u_DXv_D$, 
there exists a $B_1$-bridge $C$ of $H$ such that $N(C)\cap V(u_DXv_D-\{u_D,v_D\})\ne\emptyset$.
Let $T$ be an induced path in $G[D+\{u_D,v_D\}]-b_D-x_2y_2$ from $u_D$ to $v_D$. 
By replacing $u_DXv_D$ with $T$ we obtain a path $X'$ from $X$ such
that $y_1$ is not a cut vertex in $(G-x_1)-X'$, $B_1$ is contained in
a block of $(G-x_1)-X'$, and 
$(G-x_1)-X'$ has a chain of blocks from $y_1$ to $y_2$ and containing $B_1$, contradicting 
the choice of $X$ (in (2) as $y_2\in V(X)$).

Hence, by (5), $V(G_{\cal B})=\{D\}$. If $G$ has an edge from $u_DXv_D-\{u_D,v_D\}$ to $B_1-y_1$ or if $y_1$ has two neighbors, one on 
$u_DXy_2-u_D$ and one on $v_DXy_2-v_D$, 
then let $X'$ be obtained from $X$ by replacing $u_DXv_D$ with an induced path in $G[D+\{u_D,v_D\}]-b_D-x_2y_2$ from $u_D$ to $v_D$. 
In the former case,  $(G-x_1)-X'$ has a chain of blocks 
from $y_1$ to $y_2$ and containing $B_1$, contradicting (2). In the
latter case, $(G-x_1)-X'$ has a cycle containing $\{y_1,y_2\}$. 
So by Lemmas~\ref{4A} and 
\ref{apex}, $(i)$ or $(ii)$ holds, or there is an induced path $X^*$ in $G-x_1$ from $z_1$ to $x_2$ such that $y_1,y_2\notin V(X^*)$ 
and $(G-x_1)-X^*$ is 2-connected, and $(iii)$ holds. 

 Therefore, we may assume $N(u_DXv_D-\{u_D,v_D\})\cap V(B_1)=\{y_1\}$, and $N(y_1)\cap V(u_DXv_D-\{u_D,v_D\})\subseteq V(u_DXy_2)$ or  
$N(y_1)\cap V(u_DXv_D-\{u_D,v_D\})\subseteq V(v_DXy_2)$.  
 Let $L=G[D\cup u_DXv_D]$ and let $L'=G[L+y_1]$. 

Suppose $L$ has disjoint paths from $u_D,b_D$ to $v_D,y_2$,
respectively. We may apply Lemma~\ref{4A} to $L$ and $\{u_D,v_D,b_D,y_2\}$. If
$L$ has an induced path $S$ from $u_D$ to $v_D$ such that $L-S$ is a
chain of blocks from $b_D$ to $y_2$ then  let $X'$ be 
obtained from $X$ by replacing $u_DXv_D$ with $S$; now $(G-x_1)-X'$ is a chain of blocks 
from $y_1$ to $y_2$  and containing $B_1$, contradicting (2). So we
may assume that $L$ has a 4-separation as given in $(i)$ of
Lemma~\ref{4A}. Thus $G$ has a 5-separation $(G_1,G_2)$ such that
$x_1\in V(G_1\cap G_2)$, $|V(G_i)|\ge 7$ for $i\in [2]$, and
$(G_2-x_1,V(G_1\cap G_2)-\{x_1\})$ is planar. Hence, $(i)$ or $(ii)$
follows from Lemma~\ref{apex}.

Thus, we may assume that such disjoint paths do not exist in $L$. 
By Lemma~\ref{2path}, there exists a collection ${\cal A}$ of subsets of $V(L)-\{b_D,u_D,v_D,y_2\}$ such that 
$(L,{\cal A}, u_D,b_D,v_D,y_2)$ is 3-planar. 

We now show that $(L'-y_1v_D, u_D,b_D,v_D,y_2,y_1)$ is planar (when $N(y_1)\cap V(u_DXv_D-\{u_D,v_D\})\subseteq V(u_DXy_2)$), 
or  $(L'-y_1u_D, u_D,b_D,v_D,y_1,y_2)$ is planar (when $N(y_1)\cap V(u_DXv_D-\{u_D,v_D\})\subseteq V(v_DXy_2)$). 
Since the arguments for these two cases are the same, we consider only
the case  when $N(y_1)\cap V(u_DXv_D-\{u_D,v_D\})\subseteq V(u_DXy_2)$.  
Since $G$ is 5-connected, for each $A\in {\cal A}$,  $\{x_1,y_1\}\subseteq N(A)$ and
$|N_L(A)|=3$; and since 
$N(y_1)\cap V(L-b_D)\subseteq V(u_DXy_2)$ and $G$ is 5-connected, $|N_L(A)\cap V(X)|=2$. For each such $A$, let $a_1,a_2\in N_L(A)\cap V(X)$ and 
let $a\in N_L(A)-V(X)$. If for each $A\in{\cal A}$, $(G[A\cup \{a,a_1,a_2,y_1\}],
a_1,a,a_2,y_1)$ is planar, then $(L'-y_1v_D, u_D,b_D,v_D,y_2,y_1)$ 
is planar. So we may assume that, for some choice of $A$, $(G[A\cup \{a,a_1,a_2,y_1\}],
a_1,a,a_2,y_1)$ is not planar. (Note that
$G[A\cup \{a,a_1,a_2, y_1\}]$ is $(4, \{a,a_1,a_2,y_1\})$-connected.)
Hence, by Lemma~\ref{2path}, $G[A\cup \{a,a_1,a_2, y_1\}]$ contains
disjoint paths from $a_1,a$ to $a_2,y_1$, respectively.  So we can apply Lemma~\ref{4A} to $G[A\cup \{
a,a_1,a_2,y_1\}]$ and $\{a, a_1,a_2,y_1\}$. If $(i)$ of Lemma~\ref{4A} occurs then $G$
has a 5-separation $(G_1,G_2)$ such that $x_1\in V(G_1\cap G_2)$,
$|V(G_i)|\ge 7$ for $i\in [2]$, and $(G_2-x_1,V(G_1\cap G_2)-\{x_1\})$
is planar; so $(i)$ or $(ii)$ follows from Lemma~\ref{apex}. Hence, we may
assume that $(ii)$ of Lemma~\ref{4A} occurs. Then 
$G[A\cup  \{a,a_1,a_2, y_1\}]$ has an induced path $S$ from $a_1$ to
$a_2$ such that $G[A\cup \{a,a_1,a_2,y_1\}]-S$ 
is a chain of blocks from $y_1$ to $a$. Let $X'$ be obtained from $X$ by replacing $a_1Xa_2$ with $S$. Then the block of $(G-x_1)-X'$ 
containing $y_1$ contains $B_1$ and $G[A\cup N_L(A)\cup \{y_1\}]-S$, and $y_1$ is not a cut vertex in $(G-x_1)-X'$, contradicting (1).

Hence, $G$ has a 6-separation $(G_1, G_2)$ with $V(G_1\cap G_2)=\{b_D,
u_D,v_D,x_1,y_1,y_2\}$ and $G_2-x_1=L'-y_1v_D$ (or $G_2-x_1=L'-y_1u_D$). Since 
$(L'-y_1v_D, u_D,b_D,v_D,y_2,y_1)$ (or $(L'-y_1u_D,u_D,b_D,v_D,y_1,y_2)$)
is planar and $|V(G_2)|\ge 8$,  the assertion follows from
Lemma~\ref{6-cut2} (and then Lemma~\ref{apex}). This completes the proof of (6).

\medskip

If $y_2\in V(X)$ then by (4), (5) and (6), $H$ is 2-connected; so $(iii)$ holds. Thus we may assume $y_2\notin V(X)$. 
Then by (4),   $H$ is a chain of blocks from $y_1$ to $y_2$ and containing $B_1$, which we denote as
$B_1, \ldots,  B_k$.  We may assume  $k\ge 2$; as otherwise, $(iii)$
holds. Let $y_1\in V(B_1)-V(B_2)$, $y_2\in V(B_k)-V(B_{k-1})$, and $b_i\in V(B_i)\cap 
V(B_{i+1})$ for $i\in [k-1]$. Note that
\begin{itemize}
\item  if $z_1$ has at least two neighbors in $B_1$ then $z_0\in V(B_1)$.
\end{itemize}
For,  suppose $z_1$ has at least two neighbors in $B_1$ and $z_0\notin V(B_1)$. Let $w\in V(X)$ with $wXx_2$ minimal such that 
$w$ is a neighbor of $\bigcup_{i=2}^kB_i-b_1$ in $G-x_2y_2$. Recall
that $z_0\notin V(X)$.  Let $W$
be an induced  path in  $G[(\bigcup_{i=2}^kB_i)+w-b_1]-x_2y_2$ from
$z_0$ to $w$, and let $X'=W\cup wXx_2$. Then, since $y_1$ is not a cut
vertex of $H$, $y_1$ is not a cut vertex of $(G-x_1)-X'$. However,  the block of 
$(G-x_1)-X'$ containing $y_1$ contains $B_1+z_1$, contradicting (1). 

\medskip

We further choose $X$ so that, subject to (1), (2) and (3),  
\begin{itemize}
\item [(7)] $B_k$ is maximal. 
\end{itemize}

Let $q_1, q_2\in V(X)$ be the neighbors of $\bigcup_{i=2}^kB_i-b_1$ in $G-x_2y_2$ with $q_1Xq_2$ maximal, and assume that $z_1,q_1,q_2,x_2$ occur on $X$ in this order. 
We may assume that 
\begin{itemize}
\item [(8)] there exists $b_1'\in V(B_1-b_1)$ such that $N(q_1Xq_2-\{q_1,q_2\})\cap V(B_1-b_1)=\{b_1'\}$. 
\end{itemize}
For, otherwise, by (5),  $N(q_1Xq_2-\{q_1,q_2\})\cap
V(B_1-b_1)=\emptyset$. Hence,  $(iv)$ holds with $b=b_1$, $p_1=q_1$, and $p_2=q_2$. 
%Now suppose  $|N(q_1Xq_2-\{q_1,q_2\})\cap V(B_1-b_1)|\ge 2$. Let $X'$ be obtained from $X$ by replacing $q_1Xq_2$ with an induced path 
%in $G[\bigcup_{i=2}^kB_i+\{q_1,q_2\}]-b_1$ from $q_1$ to $q_2$.  Now $B_1$ and part of $q_1Xq_2$ are contained in a block of 
%$(G-x_1)-X'$, and $y_1$ is not a cut vertex of $(G-x_1)-X'$. This contradicts (1) and completes the proof of (8). 

\medskip
 
Thus $G$ has a separation $(G_1,G_2)$ such that $V(G_1\cap
G_2)=\{b_1,b_1', q_1,q_2,x_1,y_2\}$, $G_1=G[(B_1\cup z_1Xq_1\cup q_2Xx_2)+\{x_1,y_2\}]$ and $G_2$ contains $\bigcup _{i=2}^kB_i$ and $q_1Xq_2$. 
Note that $xy\notin E(G_2)$ for all $\{x, y\}\subseteq V(G_1\cap G_2)$. 
We may assume that 
\begin{itemize}
\item [(9)] there exists a collection  $\mathcal{A}$ of subsets of $V(G_2-x_1)-\{b_1, b_1',q_1, q_2\}$ such that 
$(G_2-x_1, \mathcal{A}, b_1, q_1, b_1', q_2)$ is $3$-planar.
\end{itemize}
For, otherwise, by Lemma~\ref{2path}, $G_2-x_1$ has disjoint paths $S, S'$ from $b_1, q_1$ to $b_1', q_2$, respectively.  
We may choose $S'$ to be induced and let $X'$ be obtained from $X$ 
by replacing $q_1Xq_2$ with $S'$. Then $B_1\cup S$ is contained in a block of $(G-x_1)-X'$. 
Thus, by (1),  $y_1=b_1'$ and $y_1$ is a cut vertex of $(G-x_1)-X'$.

%Suppose $G_2-x_1$ is $(4,\{b_1,b_1',q_1,q_2\})$-connected.  
%Applying Lemma~\ref{4A} (and then Lemma~\ref{apex}) to $G_2-x_1$ and
%$\{q_1,q_2,b_1,b_1'\}$,  we may assume that there is  
%an induced path $S^*$ in $G_2-x_1$ from $q_1$ to $q_2$
%such that $(G_2-x_1)-S^*$ is a chain of blocks. Let $X^*$ be obtained from $X$ by replacing $q_1Xq_2$ with $S^*$. Then  $B_1$ is properly 
%contained in a  block of $(G-x_1)-X^*$, and $y_1$ is not a cut vertex of $(G-x_1)-X^*$. This contradicts (1).

If $G_2-x_1$ is $(4,\{b_1,b_1',q_1,q_2\})$-connected, let $G_2'=G_2$ and $J=\emptyset=T$. 
Now suppose  $G_2-x_1$ is not $(4,\{b_1,b_1',q_1,q_2\})$-connected. Since $G$ is 5-connected and $y_2$ is the only vertex in $V(G_2)-\{b_1,b_1',q_1,q_2,x_1\}$
adjacent to $x_2$, $G_2-x_1$ has a 3-cut $T$ separating $y_2$ from $\{b_1,b_1',q_1,q_2\}$. Choose $T$ so that the component $J$ of 
$(G_2-x_1)-T$ containing $y_2$ is maximal. Let $G_2'$ be obtained from $G_2-J$ by adding an edge between every pair of vertices in $T$. 

Then $G_2'-x_1$ is $(4,\{b_1,b_1',q_1,q_2\})$-connected, and the paths $S,S'$ also give rise to disjoint paths in $G_2'-x_1$ 
from $b_1, q_1$ to $b_1', q_2$, respectively.  Hence by applying
Lemma~\ref{4A}  (and then Lemma~\ref{apex}) to $G_2'-x_1$ and $\{q_1,q_2,b_1,b_1'\}$, we find an induced path $S''$ 
in $G_2'-x_1$ from $q_1$ to $q_2$ such that $(G_2'-x_1)-S''$ is a chain of blocks 
from $b_1$ to $b_1'$. Note that $S''$ gives rise to an induced path $S^*$ in $G_2$ by replacing $S''\cap G_2'[T]$ with an induced path in $G_2[J+T]$. 
Let $X^*$ be obtained from $X$ by replacing $q_1Xq_2$ with $S^*$.  Then $B_1$ is properly contained in a block of $(G-x_1)-X^*$. 
Since $y_2\notin V(X)$, $b_1'\notin T\cup V(J)$. Hence, $y_1$ is not a cut vertex in $(G-x_1)-X^*$. 
Thus, we have a contradiction to (1) which completes the proof of  (9). 

\medskip

We may assume that, for any choice of ${\cal A}$ in (9), 
\begin{itemize}
\item [(10)] ${\cal A}\ne \emptyset$.
\end{itemize}
For, otherwise, $G_2-x_1$ has no cut of size at most 3 separating $y_2$ from $\{b_1,b_1',q_1,q_2\}$. Hence, 
$G_2$ is $(5,\{b_1, b_1',q_1,q_2,x_1\})$-connected and $(G_2-x_1, b_1, q_1, b_1', q_2)$ is planar. 
We may assume that $G_2-x_1$ is a plane graph with $b_1, q_1, b_1', q_2$ incident with its outer face. 

If $y_2$ is also incident with the outer face of $G_2-x_1$ then $(i)$
or $(ii)$ holds by applying Lemma~\ref{6-cut2} (and then
Lemma~\ref{apex}) to $G_2-x_1$ and $\{b_1,b_1',q_1,q_2,x_1,y_2\}$. 
So  assume that $y_2$ is not incident with the outer face of $G_2-x_1$. Then by Lemma~\ref{apex1}, 
the vertices of $G_2-x_1$ cofacial with $y_2$ induce a cycle $C_{y_2}$ in $G_2-x_1$, and
$G_2-x_1$ contains paths $P_1,P_2, P_3$ from $y_2$ to $\{b_1,b_1', q_1,q_2\}$  such that $V(P_i\cap P_j)=\{y_2\}$
for $1\le i<j\le 3$, and $|V(P_i\cap C_{y_2})|=|V(P_i)\cap \{b_1,b_1',q_1,q_2\}|=1$ for $i\in [3]$. Let $K=C_{y_2}\cup P_1\cup P_2\cup P_3$. 

If $P_1,P_2,P_3$ end at $q_1$, $b_1$ (or $b_1'$), $q_2$, respectively,
then let $Q$ be a path in $B_1$ from $y_1$ to $b_1$ (or $b_1'$); 
%and if $P_1,P_2,P_3$ end at $q_1,b_1',q_2$, respectively, then let $Q$ be a path in $B_1$ from $y_1$ to $b_1'$. 
now $K\cup (x_1z_1\cup z_1Xq_1)\cup (x_1x_2\cup x_2Xq_2)\cup  (x_1y_1\cup Q)\cup x_1y_2$ is a $TK_5$ in $G'$.
For the remaining cases, let $Q_1,Q_2$ be independent paths in $B_1$ from $y_1$ to $b_1', b_1$, respectively. 
If $P_1,P_2,P_3$ end at $b_1, q_1,b_1'$, respectively, then $K\cup Q_1\cup Q_2\cup (y_1x_1z_1\cup z_1Xq_1)\cup y_1x_2y_2$ is a $TK_5$ in $G'$. If 
$P_1,P_2,P_3$ end at $b_1, q_2,b_1'$, respectively then $K\cup Q_1\cup Q_2\cup (y_1x_2\cup x_2Xq_2)\cup y_1x_1y_2$ is a $TK_5$ in $G'$.  This proves (10). 

\medskip

By (10) and the 5-connectedness of $G$,  we may let
$\mathcal{A}=\{A\}$ and $y_2\in A$. Moreover,  $|N(A)-\{x_1,x_2\}|=
3$.  Choose  ${\cal A}$ so that 
\begin{itemize}
\item [(11)] $A$ is maximal.
\end{itemize}

Then 
\begin{itemize}
\item [(12)]  $b_1'\notin N(A)$, and we may assume that $N(b')\cap
  V(B_k-b_{k-1})=\emptyset$ for any $b'\in N(b_1')\cap V(q_1Xq_2)$, and  $|N(A)\cap V(q_1Xq_2)|= 2$. 
\end{itemize}
Suppose $b_1'\in N(A)$. Then $A\cap V(q_1Xq_2-\{q_1, q_2\})\neq \emptyset$. 
Hence, $|N(A)\cap V(q_1Xq_2)|\geq 2$. Since $y_2\in A$ and $y_2\notin V(X)$, $|N(A)\cap V(B_i)|\geq 1$ for some $2\le i\le k$, 
a contradiction as $|N(A)-\{x_1,x_2\}|=3$.

Now suppose there exist $b'\in N(b_1')\cap V(q_1Xq_2)$ and $b''\in N(b')\cap V(B_k-b_{k-1})$. Then $B_k$ has independent paths $P_2, P_2'$ from $y_2$ 
to $b_{k-1}$, $b''$, respectively. 
Let $P_1, P_1'$ be independent paths in $B_1$ from $y_1$ to $b_1, b_1'$, respectively, 
and let $P$ be a path  in $\bigcup_{j=2}^{k-1}B_j$ from $b_1$ to $b_{k-1}$. 
Then $G[\{x_1, x_2, y_1, y_2\}]\cup (b'Xz_1\cup z_1x_1)\cup b'Xx_2\cup (b'b_1'\cup P_1')\cup (b'b''\cup P_2')\cup (P_1\cup P\cup P_2)$ is a 
$TK_5$ in $G'$ with branch vertices $b',x_1,x_2,y_1,y_2$.

Finally, assume $|N(A)\cap V(q_1Xq_2)|\leq 1$. Then, since $B_k-b_{k-1}$ has at least two neighbors on $q_1Xq_2$ (as $G$ is 5-connected), $B_k$ is 
2-connected and $V(B_k-b_{k-1})\not\subseteq A$. Hence, $|N(A)\cap V(B_k)|\geq 2$. 
Let  $q_1',q_2'\in N(B_k-b_{k-1})\cap V(X)$ such that $q_1'Xq_2'$ is maximal. 
Then there exists $b'\in N(b_1')\cap V(q_1'Xq_2'-\{q_1', q_2'\})$; otherwise $V(B_k\cup q_1'Xq_2')-\{b_{k-1}, q_1', q_2'\}$ 
contradicts the choice of $A$ in (11). 
Since $G$ is 5-connected and $(G_2-x_1, \mathcal{A}, b_1, q_1, b_1', q_2)$ is $3$-planar, 
$b'$ has a neighbor $b''$ in $B_k-b_{k-1}$, a contradiction. 

So $|N(A)\cap V(q_1Xq_2)|\ge 2$.  Indeed $|N(A)\cap V(q_1Xq_2)|= 2$,
since $(G-x_1)-X$ is connected, $y_2\notin V(X)$, and $|N(A)-\{x_1,x_2\}|=3$.
This concludes the proof of (12). 

\medskip

Since $|N(A)\cap V(q_1Xq_2)|=2$ (by (12)), there exists  $2\le l\le k-1$ such that
$b_l\in N(A)$ and $\bigcup_{j=l+1}^kV(B_j)\subseteq A$. (Here $l\ne 1$ since $\bigcup_{j=2}^kV(B_j)\not \subseteq A$.)
Note that $N(A)\cap V(q_1Xq_2)\ne  \{q_1, q_2\}$, as $b_1'$ has a neighbor in $q_1Xq_2-\{q_1,q_2\}$.  We may assume that 
%\medskip

\begin{itemize}
\item [(13)] there exists $i\in  [2]$ such that $q_i\in N(A)$ and
$N(q_i)\cap V(G_2-x_1)\subseteq A\cup N(A)$. 
\end{itemize}
For, suppose otherwise. Then for $i\in  [2]$, $q_i\notin N(A)$ or $N(q_i)\cap V(G_2-x_1)\nsubseteq A\cup N(A)$. 
Hence, $G_2[\bigcup_{j=2}^l B_j+ \{q_1, q_2\}-b_1]$ contains an induced
path $P$ from $q_1$ to $q_2$. 

We may assume $b_1'\ne y_1$.  For, suppose $b_1'=y_1$. Since $G$ is
5-connected, there exists $t\in [2]$ such that $G[\bigcup_{j=l+1}^kV(B_j)\cup q_1Xq_2+y_1]-\{b_l,q_{3-t}\}$ has independent paths $P_1,P_2$
from $y_2$ to $y_1, q_t$, respectively. If $q_t$ has a neighbor $s\in V(B_1)$ then let $S$ be a path in $B_1$ from $s$ to $y_1$; 
now $G[\{x_1,x_2,y_1,y_2\}]\cup (x_1z_1\cup z_1Xq_1\cup P\cup q_2Xx_2)\cup  (q_ts\cup S)\cup 
P_2\cup P_1$ is a $TK_5$ in $G'$ with branch vertices $q_t,x_1,x_2,y_1,y_2$. So assume that $q_t$ has no neighbor in $B_1$. 
Then we may assume $q_t\notin \{z_1,x_2\}$ and $q_tx_2\notin E(X)$; for otherwise, 
$\{b_1,q_{3-t},x_1,x_2,y_1\}$ is a 5-cut in $G$ containing the triangle $x_1x_2y_1x_1$, 
and the assertion follows from Lemma~\ref{5cut_triangle}. 
Now let $vq_t\in E(X)-E(q_1Xq_2)$. Then $G[B_1+v]$ has independent paths $R_1,R_2$ from $v$ to $y_1,b_1$, respectively. 
Let $R$ be a path in $G[\bigcup_{j=2}^{l}B_j+q_{3-t}]$ from $b_1$ to $q_{3-t}$.  Then 
$G[\{x_1,x_2,y_1,y_2\}]\cup  R_1\cup (vq_t\cup P_2)\cup
(R_2\cup R\cup (X-(q_1Xq_2-q_{3-t}))\cup x_1z_1)\cup P_1$ 
is a $TK_5$ in $G'$ with branch vertices $v,x_1,x_2,y_1,y_2$.

Let $t_1,t_2\in V(X-x_2)\cap N(B_k-b_{k-1})$ with $t_1Xt_2$
maximal. Then $t_1\ne t_2$ as $\{b_{k-1},t_1,x_1,x_2\}$ is not a cut in $G$. 
We claim that $G[B_k\cup t_1Xt_2]-b_{k-1}$ is
2-connected. For, suppose not. Then $G[B_k\cup t_1Xt_2]$ has a
2-separation $(L_1,L_2)$ such that $b_{k-1}\in V(L_1\cap L_2)$ and
$t_1Xt_2\subseteq L_1$. Now $V(L_1\cap L_2)\cup \{x_1,x_2\}$ is a
cut in $G$, a contradiction. 

Let $X'$ be obtained from $X$ by replacing $q_1Xq_2$ with
$P$. Then $(G-x_1)-X'$ has a chain of blocks from $y_1$ to $y_2$, in
which $B_1$ is a block containing $y_1$, and the block containing
$y_2$ contains $(B_k-b_{k-1})\cup t_1Xt_2$ (whose size is larger than
$B_k$ as $t_1\ne t_2$). Since $b_1'\ne y_1$,  $y_1$ is
not a cut vertex in $(G-x_1)-X'$. This contradicts the choice of
$X$ for (7) (subject to (1), (2) and (3)).  So we have (13).

\medskip

Then $q_{3-i}\notin N(A)$ (as $N(A)\cap V(q_1Xq_2)\ne \{q_1,q_2\}$), and $x_2\ne q_i$ (otherwise $N(A)\cup
\{x_1\}$ would be a cut in $G$ of size at most 4).  Let $a\in
N(A)-\{x_1,x_2,q_i,b_l\}$. Then $a\in V(X)$ and 
$\{ a, b_1, b_1', b_l,q_{3-i}, x_1\}$ is a cut in $G$. So $G$ has a
6-separation $(G_1',G_2')$ such that $V(G_1'\cap G_2')=\{ a, b_1,
b_1', b_l,q_{3-i}, x_1\}$ and $G_2':=G_2-(A\cup \{q_i\})$. 
Note that $(G_2'-x_1, b_1, b_l, a, b_1', q_{3-i})$ is planar.

If $|V(G_2')|\ge 8$ then we may apply Lemma~\ref{6-cut2} to
$(G_1',G_2')$ and  conclude, with help from Lemma~\ref{apex}, that $(i)$ or $(ii)$ holds. 
So assume  $|V(G_2')|=6$ or   $|V(G_2')|=7$. Note that $G-x_1$ has a
separation $(Y_1,Y_2)$ such that  $V(Y_1\cap Y_2)=  \{a,b_l,q_i\}$,  $aXq_i+y_2\subseteq Y_2$, 
$Y_1=G[B_1\cup G_2'\cup X-(q_iXa-\{a,q_i\})]$. 

\medskip

{\it Case}  1. $|V(G_2')|=6$.  

Then $l=2$ and $b_2q_{3-i},aq_{3-i}, ab_1'\in E(G)$. We claim that $b_2q_i\notin E(G)$. For, suppose 
$b_2q_i\in E(G)$. Let $P$ be a path  in $\bigcup_{j=3}^{k-1}B_j$
from $b_2$ to $b_{k-1}$. Since $G$ is $5$-connected, $B_k - b_{k-1}$ has at least two neighbors 
on $q_iXa$. Thus, we may choose $a_1a_2\in E(G)$ with $a_1\in V(q_iXa-q_i)$ and $a_2\in V(B_k-b_{k-1})$.  Let $Q_1,Q_2$ be independent paths 
in $B_k$ from $y_2$ to $b_{k-1}, a_2$, respectively, and $P_1,P_2$ be independent paths in $Y_1$ from $y_1$ to $b_1,b_1'$, respectively. Now 
$G[\{x_1,x_2,y_1,y_2\}]\cup (b_2q_1\cup q_1Xz_1\cup z_1x_1)\cup (b_2q_2\cup q_2Xx_2) \cup (P\cup Q_1)\cup  (b_2b_1\cup P_1)\cup (P_2\cup b_1'a\cup aXa_1\cup a_1a_2\cup Q_2)
$ is a $TK_5$ in $G'$ with branch vertices $b_2, x_1, x_2,y_1,y_2$. 

We also claim that $ab_1\notin E(G)$. For, otherwise, let $P$ be an
induced path in $G[\bigcup_{j=3}^kB_j + q_i]$ from $q_i$ to $b_2$. 
Let $X'$ be obtained from $X$ by replacing $q_iXq_{3-i}$ with
$P\cup b_2q_{3-i}$. Then, in $(G-x_1)-X'$, there is a block containing
both $B_1$ and $a$, and $y_1$ is not a cut vertex. This contradicts (1).

If  $q_{3-i}b_1\notin E(G)$ then $(iv)$ holds with $b=b_2$, $p_j=q_i$,
$p_{3-j}=a$, and $v=q_{3-i}$. So we may assume  $q_{3-i}b_1\in
E(G)$.  Note that $q_{3-i}\ne x_1$ as $x_1\notin V(X)$ and $q_{3-i}\in V(X)$.  
 We consider two cases: $x_2\ne q_{3-i}$ and $x_2=q_{3-i}$. 

First, suppose $x_2\ne q_{3-i}$. Since
$G$ is $5$-connected, $x_2$ has at least one neighbor in $B_1-b_1'$. Thus, $G[B_1+x_2]$ has independent paths $P_1, P_2$ 
from $b_1$ to $x_2, b_1'$, respectively. If $G[Y_2+x_2]$ contains a
path $P$ from $q_i$ to $x_2$ and containing $\{a, b_2\}$ then 
$G[\{b_1, b_2, q_{3-i}\}]\cup P\cup P_1\cup (ab_1'\cup P_2)\cup aq_{3-i}\cup
 (x_2x_1z_1\cup z_1Xq_1)\cup x_2Xq_2$ is a $TK_5$ in $G'$ 
with branch vertices $a, b_1, b_2, q_{3-i}, x_2$. 
Thus, it remains to prove the existence of $P$. Note that $G[Y_2+x_2]$ is $(4, \{a, b_2,q_i, x_2\})$-connected.
First, consider the case when $G[Y_2+x_2]$ has disjoint paths from $b_2, x_2$ to $a, q_i$, 
respectively. Then by Lemma~\ref{4A} and then  Lemma~\ref{apex}, $(i)$ or $(ii)$ holds, or there is a path $S$ in $G[Y_2+x_2]$ from $a$ to $b_2$ 
such that $G[Y_2+x_2]-S$ is a chain of blocks from $q_i$ to $x_2$. Now the existence of $P$ follows from the fact that $Y_2$ is 2-connected.  
So assume $G[Y_2+x_2]$ has no disjoint paths from $b_2, x_2$ to $a, q_i$, 
respectively. Then by Lemma~\ref{2path}, $(G[Y_2+x_2], b_2, x_2, a, q_i)$ is planar. 
If $|V(G[Y_2+x_2])|\geq 6$ then the assertion of the lemma follows from Lemma~\ref{apex}. So assume $|V(G[Y_2+x_2])|=5$. 
If $ab_2\in E(G)$ then $G[\{q_i, a, b_2, y_2\}]\cong K_4^-$ (as $b_2q_i\notin E(G)$); and if
$ab_2\notin E(G)$ then $ax_1\in E(G)$ (as $G$ is 5-connected), and $G[\{q_i, a, x_1, y_2\}]$ contains a $K_4^-$ in
which $x_1$ is of degree 2. So $(ii)$ holds.

Now suppose $x_2=q_{3-i}$. Then we may assume that $b_1'\ne y_1$, for
otherwise $G[\{a,x_1,x_2,y_1\}]$ contains a $K_4^-$ in which $x_1$ is
of degree 2, and $(ii)$ holds. Thus $B_1$ has independent paths
$P_1,P_2$ from $b_1$ to $y_1,b_1'$, respectively. If $Y_2$ has a cycle
$C$ containing $\{a,b_2,y_2\}$, then $C\cup
G[\{a,b_1,b_2,q_{3-i}\}]\cup y_2x_2\cup (P_2\cup b_1'a)\cup (P_1\cup
y_1x_1y_2)$ is a $TK_5$ in $G'$ with
branch vertices $a,b_1,b_2,q_{3-i},y_2$ (as we assume $b_1q_{3-i}\in E(G)$). So we may assume that the
cycle $C$ in $Y_2$ does not exist. Since $Y_2$ is 2-connected, it
follows from Lemma~\ref{Watkins} that $Y_2$ has 2-cuts $S_u$, for $u\in
\{a,b_2,y_2\}$,
separating $u$ from $\{a,b_2,y_2\}-\{u\}$. Since $G$ is 5-connected,
we see that $S_{y_2}$ separates $\{q_i,y_2\}$ from
$\{a,b_2\}$. Hence, $d_G(b_2)=5$ and
$x_1b_2\in E(G)$. Now $G[\{b_1,b_2,x_1,x_2\}]$ contains a $K_4^-$ in
which $x_1$ is of degree 2(as we assume $b_1q_{3-i}\in E(G)$), and $(ii)$ holds.

\medskip

{\it Case} 2. $|V(G_2')|=7$.

Let $z\in V(G_2')-\{a, b_1', b_1, b_l,  q_{3-i},x_1\}$. 
Suppose $z\notin V(X)$. Then $b_1'a\in E(G)$. Since $G$ is 5-connected and $B_1$ is a block of $H$,  $zb_1'\notin E(G)$ and 
$za,zb_1,zb_l, zq_{3-i},zx_1\in E(G)$. We may assume $b_1'q_{3-i}\notin E(G)$, as 
otherwise, $G[\{a,b_1',q_{3-i},z\}]$ contains $K_4^-$ and $(ii)$ holds. Thus, $G[B_1+q_{3-i}]$ has independent paths $P_1,P_2$ from 
$b_1$ to $b_1',q_{3-i}$, respectively. Note $b_1b_l\in E(G)$ by the maximality of $A$ in (11). 
In $G[A\cup \{a,b_l,q_i\}]$ we find independent paths $Q_1,Q_2$ from $b_l$ to $q_i,a$, respectively. 
Now $G[\{a,b_1,b_l,q_{3-i},z\}]\cup (P_1\cup b_1'a)\cup P_2\cup Q_2\cup (q_2Xx_2\cup x_2x_1z_1\cup z_1Xq_1\cup Q_1)$ is a  
$TK_5$ in $G'$ with branch vertices  $a,b_1,b_l,q_{3-i},z$.  

So we may assume  $z\in V(X)$.  Then $b_1b_l, q_{3-i}b_l\in E(G)$. By (9),  $b_1a,b_1z\notin E(G)$. 
%For, suppose $b_1a\in E(G)$ or $b_1z\in E(G)$.  
%Let $X'$ be obtained from $X$ by replacing $q_1Xq_2$ with $b_lq_{3-i}$ and a path in $Y_2-a$ from $b_l$ to $q_i$. 
%Then, $B_1+a$ or $B_1+z$ is contained in a block of $(G-x_1)-X'$, and $y_1$ is not a cut vertex of $(G-x_1)-X'$, contradicting (1).
Hence, since $G$ is 5-connected, $zb_1'$, $zb_l$, $zx_1\in E(G)$ and $q_{3-i}\ne x_1$. We may assume
$x_1q_{3-i}\notin E(G)$; as otherwise, $G[\{b_l,q_{3-i},x_1,z\}]$
contains a $K_4^-$ in which $x_1$ is of degree 2, and $(ii)$ holds. 
Note that $b_1'a\in E(G)$ by the maximality of $A$ in (11). 
Let $q\in N(q_{3-i})\cap V(B_1-b_1)$, and let $P_1,P_2$ be independent paths in $B_1$ from $b_1'$ to $b_1,q$, respectively. 
Let $Q_1,Q_2$ be independent paths in $Y_2$ from $a$ to $b_l,q_i$, respectively. Then $G[\{a,b_l,b_1',q_{3-i},z\}]\cup (P_1\cup b_1b_l)\cup  
(P_2\cup qq_{3-i})\cup Q_1\cup (Q_2\cup q_1Xz_1\cup z_1x_1x_2 \cup x_2Xq_2)$ is a $TK_5$ in $G'$ with branch vertices $a, b_1',b_l,q_{3-i},z$. \qed

\section{Two special cases}
To prove Theorem~\ref{x1a},
we need to take care of the conclusions $(iii)$ and $(iv)$ of Lemma~\ref{classify2}. 
Results from \cite{MY10} can be used to deal with 
$(iii)$ of Lemma~\ref{classify2} when $y_2\notin V(X)$. So it remains to consider  $(iii)$ of Lemma~\ref{classify2} 
when $y_2\in V(X)$ and $(iv)$ of Lemma~\ref{classify2}. 

We will use
the notation in  the statement of Lemma~\ref{classify2}. See Figures~\ref{path} and \ref{blocks}. In particular, 
$X$ is an induced path in 
$(G-x_1)-x_2y_2$ from $z_1$ to $x_2$ and $G':=G-\{x_1x : x\notin
\{x_2, y_1, y_2, z_0, z_1\}\}.$ Also recall from   $(iv)$ of
Lemma~\ref{classify2} 
the separation $(Y_1,Y_2)$ and the vertices $p_j, p_{3-j}, v, b, b_1,
b_1'$. Note that $x_2\ne p_2$; as otherwise, $\{b,p_1,x_1,x_2\}$ would a cut in $G$.
Let $z_2$ be the neighbor of $x_2$ on $X$. For any  $x\in V(G)$ and $S\subseteq G$, we use $e(x,S)$ to denote the number of edges in $G$ from $x$ to $S$.

In this section, we deal with two special cases of Theorem~\ref{x1a}.
First, we need some structural information on $Y_2$.

\begin{lem}\label{Cycle}
Suppose $(iv)$ of Lemma~\ref{classify2} holds. 
%$y_2\notin V(X)$. 
Then $Y_2$ has independent paths from $y_2$
to $b, p_1, p_2$, respectively, and, for $i\in [2]$, $Y_2$ has \textsl{a} path
from $b$ to $p_{3-i}$ and containing $\{y_2, p_i\}$. 
Moreover, one of the following holds:
\begin{itemize}
\item [$(i)$] $G'$ contains $TK_5$, or $G$ contains \textsl{a} $TK_5$
  in which $x_1$ is not \textsl{a} branch vertex. 
\item [$(ii)$] $G-x_1$ contains $K_4^-$, or $G$ contains \textsl{a} $K_4^-$ in which $x_1$ is of degree 2.
\item [$(iii)$] %Let $i\in [2]$. Then $Y_2$ has a path from $b$ to $p_{3-i}$ and containing $\{y_2, p_i\}$. 
If  $e(p_i, B_1-b_1)\ge 1$ for some $i\in [2]$ then $Y_2$ has a path through $b, p_i,y_2,p_{3-i}$ 
in order, and $Y_2-b_1$ has a cycle containing $\{p_1, p_2, y_2\}$. If $b\ne b_1$, $p_2v\in E(X)$, and $vb, vx_1\in E(G)$ then 
$Y_2$ has a cycle containing $\{b,p_2,y_2\}$. 
\end{itemize}
\end{lem}

\pf Since $G$ is 5-connected, $Y_2$ is $(3,\{b,p_1,p_2\})$-connected. So by Menger's theorem, 
$Y_2$ has independent paths from $y_2$ to $b,p_1, p_2$,
respectively. 

\medskip

Next, let $i\in [2]$. We claim that $Y_2$ has a path $Q_i$ from $b$ to $p_{3-i}$ and containing $\{p_i,y_2\}$. 
To see this, let $Y_2':=Y_2+\{t,tb,tp_{3-i}\}$, which is
2-connected. If $Y_2'$ has a cycle $C$ containing $\{b,t,y_2\}$ then
$Q_i:=C-t$ is as desired.
%is a path in $Y_2$ from $b$ to $p_{3-i}$ and containing $\{y_2, p_i\}$. 
So suppose such a
cycle $C$ does not exist. Then by Lemma~\ref{Watkins}, $Y_2'$ has a
2-cut $T$ separating $y_2$ from $\{p_i,t\}$ and $\{p_i,t\}\cap T=\emptyset$. However, $T\cup
\{x_1,x_2\}$ is a cut in $G$, a contradiction.

\medskip
We now show that $(i)$ holds or  the first part of $(iii)$ holds. Suppose $e(p_i,B_1-b_1)\ge
1$ for some $i\in [2]$. 
%Let $S$ denote a path in $Y_2$ from $b$ to $p_{3-i}$ and containing $\{p_i,y_2\}$. 

First, we may assume that $Q_i$ must go through $b,p_i,y_2,p_{3-i}$ in
order. For, suppose $Q_i$ goes through $b,y_2,p_i,p_{3-i}$ in this
order. Since
$e(p_i,B_1-b_1)\ge 1$,  $G[B_1+p_i]$ has independent paths $P_1,P_2$
from $y_1$ to $b_1, p_i$, respectively. Then $G[\{x_1, x_2, y_1, y_2\}]\cup Q_i
\cup P_2\cup ((X-(p_1Xp_2-\{p_1,p_2\}))\cup x_1z_1)\cup (P_1\cup b_1b)$ is a $TK_5$ in $G'$ with branch vertices 
$p_i, x_1,x_2,y_1,y_2$, and $(i)$ holds.

Next, note that $Y_2-b_1$ is 2-connected. For, suppose not. Then $b=b_1$ and 
$Y_2-b_1$ has a 1-separation $(Y_{21},Y_{22})$, and we may assume $|V(Y_{21}-Y_{22})\cap \{p_1,p_2,y_2\}|\le 1$. 
Since each of $\{p_1,p_2,y_2\}$ has at least two neighbors in $Y_2-b_1$,  $(V(Y_{21}-Y_{22})\cap \{p_1,p_2,y_2\})\cup \{b,x_1\}\cup V(Y_{21}\cap Y_{22})$ 
is a cut in $G$ of size at most 4, a contradiction. 
%Thus $Y_2-b_1$ is 2-connected.

Now suppose no cycle in  $Y_2-b_1$ contains  $\{p_1, p_2, y_2\}$. 
Then, $(i)$ or $(ii)$ or $(iii)$ of Lemma~\ref{Watkins} holds. We use the notation 
in Lemma~\ref{Watkins} (with  $p_1, p_2, y_2$ playing the roles of  $y_1,y_2,y_3$ there). 
If $(i)$ of Lemma~\ref{Watkins} occurs then let $S=\{a_1, a_1'\}$, $a_2=a_3=a_1$, and $a_2'=a_3'=a_1'$;
if $(ii)$ or $(iii)$ of Lemma~\ref{Watkins} occurs let $S_{p_j}=\{a_j,
a_j'\}$ for $j\in [2]$ and let $S_{y_2}=\{a_3, a_3'\}$.  
Let $A, A'$ denote the components of $(Y_2-b_1)-(D_{p_1}\cup D_{p_2}\cup D_{y_2})$ such that $a_j\in V(A)$ and $a_j'\in V(A')$ for $j\in [3]$. 
Note that if $(ii)$ of Lemma~\ref{Watkins} occurs and $A\ne A'$, then either $A= a_3$ and $\{a_1',a_2',a_3'\}\subseteq V(A')$, or 
$A'=a'_3$ and $\{a_1,a_2,a_3\}\subseteq V(A)$.

Since $Y_2-b_1$ is 2-connected,  there exist paths $S_1, S_2, S_3$ in $D_{p_1}, D_{p_2}, D_{y_2}$, respectively, 
with $S_j$ from $a_j$ to $a_j'$ for $j\in [3]$, $p_j\in V(S_j)$ for $j\in [2]$, and $y_2\in V(S_3)$. 
Since $G$ is 5-connected, $b\in V(D_{y_2})$ or $b=b_1$ has a neighbor in $D_{y_2}$. 
Hence, $G[D_{y_2}+b]$ contains a path $T_3$ from $b$ to some $t\in
V(S_3)-\{a_3,a'_3\}$ and internally disjoint from $S_3$. By symmetry, we may assume $t\in V(y_2S_3a_3)$. 
Let $T_1$ be a path in $A$ from $a_i$ to $a_{3-i}$, and $T_2$ be a path in $A'$ from $a_i'$ to $a_3'$.  
Then $T_3\cup tS_3a_3'\cup T_2\cup S_i\cup T_1\cup
a_{3-i}S_{3-i}p_{3-i}$ is a path from $b$ to $p_{3-i}$ through $y_2,
p_i$ in order. 
This is a contradiction as we have assumed that such a path $Q_i$ does not exist.

\medskip

Next, we prove that  $(i)$ or $(ii)$ holds or the second part of
$(iii)$ holds. Suppose $b\ne b_1$,
$p_2v\in E(p_2Xx_2)$, and $vb, vx_1\in E(G)$. Suppose $Y_2$ has no
cycle containing  $\{b,p_2,y_2\}$. Then $(i)$ or $(ii)$ or $(iii)$ of Lemma~\ref{Watkins} holds. 
In particular, $\{b,p_2,y_2\}$ is independent in $G$. We use the notation 
in Lemma~\ref{Watkins} (with  $b,p_2,y_2$ playing the roles of
$y_1,y_2,y_3$ there, respectively). So there is a 2-cut
$S_{y_2}=\{a_3, a_3'\}$ in $Y_2$ such that $Y_2-S_{y_2}$ has a
component $D_{y_2}$ with $y_2\in V(D_{y_2})$ and $b,p_2\notin
V(D_{y_2})\cup S_{y_2}$. Since $G$ is 5-connected, $p_1\in V(D_{y_2})$. Note that $Y_2-D_{y_2}$ is $(4, \{a_3, a_3',b, p_2\})$-connected. 

Suppose $(Y_2-D_{y_2}, a_3, b, a_3', p_2)$ is not planar. Then by
Lemma~\ref{2path}, $Y_2-D_{y_2}$ contains disjoint paths from $a_3, b$
to $a_3', p_i$, respectively. By Lemma~\ref{4A}, we may assume that
$Y_2-D_{y_2}$ has an induced path $S$ from $b$ to $p_2$ 
such that $(Y_2-D_{y_2})-S$ is a chain of blocks from $a_3$ to $a_3'$;
for otherwise, we may apply Lemma~\ref{apex} to show that $(i)$ or
$(ii)$ holds. Thus $Y_2-D_{y_2}$ has a path $S_1$ from $a_3$ to $a_3'$
and containing $\{b, p_2\}$ (as $Y_2$ is 2-connected). Let $S_2$ be a path in $G[D_{y_2}+\{a_3, a_3'\}]$ from $a_3$ to
 $a_3'$ through $y_2$. Then $S_1\cup S_2$ is a cycle in $Y_2$ containing $\{b,
 p_2, y_2\}$, a contradiction. 

So we may assume $(Y_2-D_{y_2}, a_3, b, a_3', p_2)$ is planar. 
%Hence, $bp_2 \notin E(G)$. 
If $|V(Y_2-D_{y_2})|\ge 6$ then $(i)$ or $(ii)$ follows from
Lemma~\ref{apex} (by considering the 5-cut $\{a_3,a_3',b,p_i,x_1\}$).  

Now suppose $|V(Y_2-D_{y_2})|=5$. Let $t\in
V(Y_2-D_{y_2})-\{a_3,a_3',b,p_2\}$. Since $G$ is 5-connected,
$ta_3,ta_3', tb, tp_2,tx_1\in E(G)$. By symmetry between $a_3$ and
$a_3'$, we may assume $a_3'\in V(X)$.  Then $a_3'p_2\in E(G)$. If
$ba_3'\in E(G)$ then $G[\{a_3',b,p_2,t\}]\cong K_4^-$, and $(ii)$
holds. So assume $ba_3'\notin E(G)$. Then, since $G$ is 5-connected,
$ba_3,bx_1\in E(G)$. Now $G[\{a_3,b,t,x_1\}]$ contains a $K_4^-$
in which $x_1$ is of degree 2, and $(ii)$ holds.

So $|V(Y_2-D_{y_2})|=4$ and, hence,  $(i)$ of Lemma~\ref{Watkins} occurs, with $V(D_{b})=\{b\}$ and $V(D_{p_2}) = \{p_2\}$.
We claim that $D:=G[D_{y_2}+\{a_3,a_3',x_1\}]+\{c,cx_1,cy_2\}$ has a cycle $C$ containing $\{c, a_3, a_3'\}$; for otherwise, 
by Lemma~\ref{Watkins}, $D-c$ has a 2-cut either separating 
$a_3$ from $\{x_1, y_2, a_3', p_1\}$ or separating $a_3'$ from $\{x_1,
y_2, a_3, p_1\}$, contradicting the 5-connectedness of $G$. 
Let $Q$ be a path in $G[B_1+\{b, p_2\}]$ from $b$ to $p_2$. Now $G[\{a_3,a_3',b,p_2\}]\cup Q\cup (C-c)\cup vx_1\cup (vXx_2\cup 
x_2y_2)\cup vb\cup vp_2$ is a $TK_5$ in $G$ with branch vertices $a_3,a_3',b,p_2,v$. \qed

\medskip

The next two results provide information on $e(z_i,B_1)$ for $i\in
[2]$ in the case when $y_2\notin V(X)$.

\begin{lem}\label{(I)'}
Suppose $(iv)$ of Lemma~\ref{classify2} holds with $b\ne b_1$. Then one of the following holds:
\begin{itemize}
\item [$(i)$] $G'$ contains  $TK_5$, or $G$ contains \textsl{a} $TK_5$
  in which $x_1$ is not \textsl{a} branch vertex.
\item [$(ii)$] $G-x_1$ contains $K_4^-$, or $G$ contains \textsl{a} $K_4^-$ in which $x_1$ is of degree 2.
\item [$(iii)$] $e(z_i, B_1)\ge 2$ for $i\in [2]$.
\end{itemize} 
\end{lem}

\pf Recall the notation from $(iv)$ of Lemma~\ref{classify2}. In
particular, $v\in V(X)-V(p_1Xp_2)$. Suppose $e(z_i, B_1)\le 1$ for some $i\in [2]$.  

\medskip

{\it Case} 1.  $v\in V(z_1Xp_1-p_1)$; so $p_1v\in E(X)$. 

In this case,  $e(z_1, Y_2)\le 2$ (with equality only if $z_1=v$). Hence,  $e(z_1, B_1)\ge 2$, since $G$ is 5-connected. 
Thus, $e(z_2, B_1)\le 1$.  Hence, $z_2=p_2$ and $e(z_2,B_1)=1$, since $\{x_1, x_2, p_1, b\}$
cannot be a cut in $G$. 
By Lemma~\ref{Cycle}, $Y_2$ has a path $Q$ from $b$ to $p_1$ and
containing $\{y_2, z_2\}$. 

Suppose $b, z_2, y_2, p_1$ occur on $Q$ in this order. If $b_1'\in N(z_2)$ then let $P_1, P_2$ be independent paths in 
$G[B_1+x_2]$ from $b_1'$ to $y_1, x_2$, respectively; now
$G[\{x_1,x_2,y_2\}]\cup z_2x_2\cup (z_2Qb\cup bv\cup vXz_1\cup z_1x_1) \cup z_2Qy_2\cup 
b_1'z_2\cup (b_1'p_1\cup p_1Qy_2)\cup (P_1\cup y_1x_1) \cup P_2$ is a $TK_5$ in $G'$ with branch vertices $b_1',x_1,x_2,y_2, z_2$. 
So assume $b_1'\notin N(z_2)$. Let $P_1, P_2$ be independent paths in $G[B_1+z_2]$ from $y_1$ to $ b_1', z_2$, respectively. 
Then $ G[\{x_1,x_2, y_1, y_2\}]\cup z_2x_2\cup (z_2Qb\cup bv\cup
vXz_1\cup z_1x_1)\cup z_2Qy_2\cup P_2\cup (y_2Qp_1\cup p_1b_1'\cup P_1)$ is a $TK_5$ in $G'$ with branch vertices $x_1,x_2, y_1, y_2,z_2$. 

So assume that $b, y_2, z_2, p_1$ must occur on $Q$ in this order. 
Then, by Lemma~\ref{Cycle}, we may assume $e(z_2,B_1-b_1)=0$; 
%Since $G$ is 5-connected and $p_2=z_2$, 
so $b_1z_2\in E(G)$. 
%as otherwise, $\{b,p_1,x_1,x_2\}$ would be a cut in $G$. 
Let $P_1, P_2$ be independent paths in $G[B_1+x_2]$ from $b_1$ to $y_1, x_2$, respectively.
Then $G[\{x_1,x_2,y_2\}]\cup  z_2x_2\cup (z_2Qp_1\cup p_1Xz_1\cup
z_1x_1)\cup z_2Qy_2\cup (b_1b\cup bQy_2)\cup b_1z_2\cup (P_1\cup
y_1x_1)\cup P_2$ is a $TK_5$ in $G'$ with branch vertices $b_1, x_1,x_2,y_2,z_2$.

\medskip

{\it Case} 2. $v\in V(p_2Xx_2-p_2)$; so $p_2v\in E(X)$.

Since $\{b,p_2,x_1, x_2\}$ cannot be a cut in $G$, $e(z_1, B_1) \ge  1$. We consider two cases. 

\medskip
{\it Subcase} 2.1.  $e(z_1, B_1) =  1$. 

Then $z_1=p_1$. By Lemma~\ref{Cycle}, $Y_2$ has a path $Q$ from $b$ to
$p_2$ and containing $\{z_1, y_2\}$. 

Suppose $b, z_1, y_2, p_2$ occur on $Q$ in this order.    If $b_1'\in
N(z_1)$ then  $x_2\ne v$ as $\{x_1, x_2, b_1, b_1'\}$ is not a cut in
$G$; so $e(x_2, B_1-y_1)\ge 1$. 
 Let $P_1, P_2$ be independent paths in $G[B_1+x_2]$ from $b_1'$ to $y_1, x_2$, respectively. 
Then $G[\{x_1,x_2,y_2\}] \cup z_1x_1\cup (z_1Qb\cup bv\cup vXx_2)\cup
z_1Qy_2\cup b_1'z_1\cup  (b_1'p_2\cup p_2Qy_2)\cup (P_1\cup y_1x_1)\cup P_2
$ is a $TK_5$ in $G'$ with branch vertices
$b_1',x_1,x_2,y_2, z_1$. Hence, we may assume  $b_1'\notin N(z_1)$. Then let $P_1, P_2$ be independent paths in 
$G[B_1+z_1]$ from $y_1$ to $ b_1', z_1$, respectively; now 
$G[\{x_1,x_2, y_1, y_2\}] \cup z_1x_1\cup (z_1Qb\cup bv\cup vXx_2)\cup
z_1Qy_2\cup P_2\cup (y_2Qp_2\cup p_2b_1'\cup P_1)$ is a $TK_5$ in $G'$ with branch vertices $x_1,x_2, y_1, y_2, z_1$.

So we may assume $b, y_2, z_1, p_2$ must occur on $Q$ in this
order. Hence, by Lemma~\ref{Cycle}, we may assume $e(p_1,B_1-b_1)=0$;
so $b_1\in N(z_1)$ as $\{b,p_2,x_1,x_2\}$ is not a cut in $G$. Then 
$e(x_2, B_1-y_1)\ge 1$;  otherwise, $x_2 = v$, and $\{b_1, b_1', x_1, x_2\}$ would be a cut in $G$.
Let $P_1, P_2$ be independent paths in $G[B_1+x_2]$ from $b_1$ to $y_1, x_2$, respectively. Then 
$G[\{x_1,x_2,y_2\}] \cup z_1x_1\cup (z_1Qp_2\cup p_2Xx_2)\cup z_1Qy_2\cup b_1z_1\cup (b_1b\cup bQy_2)\cup (P_1\cup y_1x_1)\cup P_2$ is a $TK_5$ in $G'$ with branch vertices $b_1,z_1,x_1,x_2,y_2$.

\medskip

{\it Subcase} 2.2. $e(z_1, B_1)\ge 2$. 

Then  $e(z_2, B_1)\le 1$. Hence, $z_2=p_2$ or $z_2 = v$. Suppose $z_2
= p_2$. Then $x_2= v$; so $x_1v\in E(G)$. Hence, 
by $(iii)$ of Lemma~\ref{Cycle}, $Y_2$ has a cycle $C$ containing
$\{b, y_2, z_2\}$. Let $P_1, P_2$ be independent paths in $B_1$ from
$y_1$ to $b_1, b_1'$, respectively. Now $C\cup x_2y_2\cup x_2z_2\cup
x_2b\cup y_1x_2\cup y_1x_1y_2\cup (P_1\cup b_1b)\cup (P_2\cup b_1'z_2)$ is a $TK_5$ in $G'$ 
with branch vertices $b, x_2, y_1, y_2, z_2$.

So we may assume $z_2 = v$. Since $e(z_2, B_1)=1$, $x_1v\in E(G)$.
Hence, by $(iii)$ of Lemma~\ref{Cycle}, $Y_2$ has a cycle $C$ 
containing $\{b, p_2, y_2\}$. Let $P_1, P_2$ be independent paths in $G[B_1+x_2]$ from $x_2$ to $b_1, b_1'$, 
respectively. Note that $P_1, P_2$ exist since $x_2$ has at least two
neighbors in $B_1$. Then $C\cup z_2b\cup z_2p_2\cup z_2x_1y_2\cup x_2y_2\cup x_2z_2\cup (P_1\cup b_1b)\cup (P_2\cup b_1'p_2)$ is a $TK_5$ in 
$G'$ with branch vertices $b, p_2,x_2,y_2, z_2$. \qed

\begin{lem}\label{(I)}
Suppose $y_2\notin V(X)$. Then one of the following holds: 
\begin{itemize}
\item [$(i)$] $G'$ contains  $TK_5$, or $G$ contains  \textsl{a} $TK_5$ in which
  $x_1$ is not \textsl{a} branch vertex.
\item [$(ii)$] $G-x_1$ contains $K_4^-$, or $G$ contains $K_4^-$ in which $x_1$ is of degree 2.
\item [$(iii)$] There exists  $i\in [2]$ such that $e(z_i, B_1-b_1)\ge 2$ and $e(z_{3-i}, B_1-b_1)\ge 1$.
\end{itemize} 
\end{lem}

\pf Suppose $(iii)$ fails. First, assume $b\ne b_1$; so $(iv)$ of
Lemma~\ref{classify2} occurs. Then by
Lemma~\ref{(I)'}, we have, for $i\in [2]$, $e(z_i,
B_1-b_1)=1$ and $b_1z_i \in E(G)$. Let $P_1, P_2$ be independent paths
in $B_1$ from $y_1$ to $b_1, b_1'$, respectively. 
Recall, from $(iv)$ of Lemma~\ref{classify2}, the role of $j\in [2]$
and the vertices $p_{3-j},v$. Since $b_1'$ is the only neighbor of
$p_{3-j}$ in $B_1$, $p_{3-j}\notin \{z_1, z_2\}$. Let $Q$ be a path in
$Y_2-\{z_1, z_2\}$ from $b$ to $p_{3-j}$ and through $y_2$. Then
$G[\{x_1, x_2, y_1, y_2\}]\cup b_1z_1x_1\cup b_1z_2x_2\cup (b_1b\cup bQy_2)\cup P_1
\cup (y_2Qp_{3-j}\cup p_{3-j}b_1'\cup P_2)$ is a $TK_5$ in $G'$ with branch vertices $b_1,x_1, x_2, y_1, y_2$. 

So we may  assume $b = b_1$. Then, for $i\in [2]$, $e(z_i, B_1-b_1)\ge 1$ as 
$\{b,p_{3-i},x_1,x_2\}$ is not a cut in $G$. Hence, since $(iii)$
fails,   $e(z_i, B_1-b_1)= 1$ for $i\in [2]$.  For $i\in [2]$, let
$z_i'\in N(z_i)\cap V(B_1)$. 
Since $G$ is 5-connected, $z_1=p_1$.

\medskip
{\it Case} 1. $z_2\ne p_2$.

Then, since $G$ is 5-connected, $z_2x_1,z_2b\in E(G)$. First, assume
that  there is no
edge from $p_2Xz_2-z_2$ to $B_1-b$. Then  $G$ has a separation $(G_1,G_2)$ such that $V(G_1\cap
G_2)=\{b,x_1,x_2,z_1,z_2\}$, $B_1\subseteq G_1$, and $Y_2\subseteq
G_2$. Clearly, $|V(G_i)|\ge 7$ for $i\in [2]$. Since $x_1x_2z_2x_1$ is a triangle in
$G$, the assertion of the lemma follows from
Lemma~\ref{5cut_triangle}. 

Hence, we may assume that there exists $uu'\in E(G)$ with $u\in
V(p_2Xz_2-z_2)$ and $u'\in V(B_1-b)$. Suppose, for some choice of $uu'$,
$u'\ne z_1'$ and  $B_1-b$ contains independent paths $P_1,P_2$ from $y_1$ to $z_1',u'$,
respectively. By Lemma~\ref{Cycle} (since $e(p_1,B_1-b_1)=1$), $Y_2$
contains a path $Q$ from $b$ to $p_2$ through $p_1,y_2$ in order. Now
$G[\{x_1,x_2,y_1,y_2\}]\cup z_1x_1\cup (z_1Qb\cup bz_2x_2)\cup
(z_1z_1'\cup P_1)\cup z_1Qy_2\cup (P_2\cup u'u\cup uXp_2\cup p_2Qy_2)$
is a $TK_5$ in $G'$ with branch vertices $x_1,x_2,y_1,y_2,z_1$.

Therefore, we may assume that for any choice of $uu'$, $u'=z_1'$ or the paths
$P_1,P_2$ do not exist.  If $u'=z_1'$ for all $u'\in N(p_2Xz_2-z_2)$ we let
$B'=B_1$ and $B''=\{b,z_1'\}$; otherwise, 
since $B_1$ is 2-connected, $B_1$ has a
2-separation $(B',B'')$ such that $b\in V(B'\cap B'')$, $y_1\in V(B')$
and $z_1',u'\in V(B'')$ for all $u'\in N(p_2Xz_2-z_2)$. Thus $G$ has a
5-separation $(G_1,G_2)$ such that $V(G_1\cap G_2)=V(B'\cap B'')\cup
\{x_1,x_2,z_2\}$, $B'\subseteq G_1$ and $B''\cup Y_2\subseteq
G_2$. Clearly, $|V(G_2)|\ge 7$. 

If $|V(G_1)|\ge 7$ then the assertion
of the lemma follows from Lemma~\ref{5cut_triangle} (as $x_1x_2z_2x_1$
is a triangle in $G$). So assume $|V(G_1)|\le 6$. Then, since $G$ is 5-connected, $z_2y_1\in E(G)$. So
$G[\{x_1,x_2,y_1,z_2\}]-x_1y_1\cong K_4^-$ in which $x_1$ is of degree
2, and $(ii)$ holds.

\medskip 
{\it Case} 2. $z_2=p_2$. 

We may assume $z_i'\ne y_1$ for $i\in [2]$. For, otherwise, $G$ has a 5-separation $(G_1,G_2)$ such that 
$V(G_1\cap G_2)=\{b,p_{3-i},x_1, x_2,y_1\}$, $B_1\subseteq G_1$ and $Y_2\subseteq G_2$. Clearly, $|V(G_2)|\ge 7$. 
If $|V(G_1)|\ge 7$ then, since $G[\{x_1, x_2, y_1\}]\cong K_3$, the assertion of the lemma follows from
Lemma~\ref{5cut_triangle}. So we may assume $|V(G_1)=6$. Then $|V(B_1)|=3$. Let $z\in V(B_1)-\{b_1,y_1\}$. Then, since $G$ is 
5-connected, $zx_1,zx_2,zy_1\in E(G)$; so $G[\{x_1,x_2,y_1,z\}]-x_1x_2\cong K_4^-$ in which $x_1$ is of degree 2, and $(ii)$ holds. 

Note that $z_1'\ne z_2'$ as otherwise $\{b_1,x_1,x_2,z_1'\}$ would be a cut in $G$. 
Let $K=G[B_1 +\{x_2, z_1, z_2\}]$.
Suppose $K$ contains disjoint paths $Z_1, Z_2$ from $z_1, z_2$ to $x_2,
y_1$, respectively. By $(iii)$ of Lemma~\ref{Cycle}, let $C$ be a cycle in $Y_2-b_1$  
containing $\{y_2, z_1, z_2\}$. Then 
$G[\{x_1, x_2, y_2\}]\cup C\cup z_1x_1\cup z_2x_2\cup (Z_2\cup y_1x_1)\cup Z_1$ is a $TK_5$ in $G'$ with branch vertices $x_1,x_2,y_2,z_1, z_2$.

So we may assume that such $Z_1, Z_2$ do not exist. Then by
Lemma~\ref{2path}, there exists a collection ${\cal A}$ of pairwise disjoint subsets of $V(K)-\{x_2,y_1,z_1, z_2\}$ 
such that $(K, {\cal A}, z_1,z_2,x_2,y_1)$ is 3-planar.
Since $G$ is 5-connected, either ${\cal A}=\emptyset$ or $|{\cal
  A}|=1$. When $|{\cal A}|=1$  let ${\cal A}=\{A\}$; then  $b_1\in A$. We choose ${\cal A}$ 
so that $|{\cal A}|$ is minimal and, subject to this, $|A|$ is minimal
when ${\cal A}=\{A\}$. 
Note that if $A$ exists then $|A|\ge 2$ (by the minimality of $|{\cal A}|$ and $|A|$). 
Moreover,  $|N_K(A)|=3$ as $N_K(A)\cup \{b_1,x_1\}$ is not a cut in $G$.

We may assume that  if ${\cal A}\ne \emptyset$ then $\{x_2,z_1,z_2\}\cap N_K(A)=\emptyset$. For, suppose
there exists $u\in \{x_2,z_1,z_2\}\cap N_K(A)$. Let $S:=(N_K(A)\cup
\{x_1,x_2,z_1,z_2\})-\{u\}$ if $u\in \{z_1,z_2\}$ and let $S:=N_K(A)\cup
\{x_1,x_2,z_1,z_2\}$ if $u=x_2$. Then $S$ is a cut in $G$ separating $B_1-A$ from $Y_2$.
Since $G$ is 5-connected,
$|S|=5$ if $u\in \{z_1,z_2\}$ and $|S|\in \{5,6\}$ if $u=x_2$. Therefore, $G$ has a
separation $(G_1,G_2)$ such that $V(G_1\cap G_2)=S$, $B_1-A\subseteq
G_1$, and $Y_2\subseteq G_2$. Note that $(G_1-x_1,S-\{x_1\})$ is
planar. Since $|V(Y_2)|\ge 4$, $|V(G_2)|\ge 7$ if $|S|=5$, and $|V(G_2)|\ge 8$ is $|S|=6$. Since
$b_1,y_1\notin \{z_1',z_2'\}$, $|V(G_1)|\ge 7$ if $|S|=5$ and $|V(G_1)|\ge 8$ if $|S|=6$. Thus,  if $|S|=5$ then the assertion of the
lemma follows from Lemma~\ref{apex}, and if $|S|=6$ then the assertion of the
lemma follows from Lemma~\ref{6-cut2} and Lemma~\ref{apex}.

If ${\cal A}=\emptyset$ let $K^*=K$; and if ${\cal A}\ne \emptyset$ let 
$K^*$ be the graph obtained from $K$ by deleting $A$ and adding new
edges 
joining every pair of distinct vertices in $N_{K}(A)$. Since $B_1$ is
2-connected and $G$ is 5-connected, 
$K':= K^*-\{x_2,z_1, z_2\}$ is a 2-connected planar graph. 
Take a plane embedding of $K'$ and let $D$ denote its outer
cycle. Let $t\in V(D)$ such that $t\in N(x_2)$
and $tDz_2'$ is minimal.

When ${\cal A}\ne \emptyset$,  $N_K(A)\not\subseteq V(D)$; as otherwise, if we write
$N_K(A)=\{s_1,s_2,s_3\}\subseteq V(D)$ with  $s_2\in V(s_1Ds_3)$, then $\{b_1,s_1,s_3,x_1\}$ is a cut in $G$, a
contradiction.  Further, if  ${\cal A}\ne \emptyset$ let $N_K(A)=\{a,a_1,a_2\}$ with $a\in N_{K}(A)-V(tDz_1')$; 
then,  by the minimality of ${\cal A}$ and $A$, $G[A\cup N_K(A)]$ contains disjoint  paths $P_1,P_2$ from $a,a_2$ to $b_1,a_1$, respectively. 
If ${\cal A}=\emptyset$ let $Q=tDz_1'$, $P_1=a=a_1=a_2=b_1$ and $P_2=\emptyset$. If ${\cal A}\ne \emptyset$ 
let $Q=tDz_1'$ if $a_1a_2\notin E(tDz_1')$;  and otherwise let $Q=(tDz_1'-a_1a_2)\cup P_2$. Note that $Q$ is a path in $B_1$. 

Suppose $K'-(tDz_1'-z_2')$ has independent paths $S_1, S_2$ from $y_1$
to $z_2', \{a,a_1,a_2\}$, respectively, and internally disjoint from
$\{a,a_1,a_2\}$. We may assume the notation is chosen so that $a\in
V(S_2)$. For $i\in [2]$, let $S_i'=S_i$ if $a_1a_2\notin E(S_i)$; and
otherwise let  $S_i'$ be obtained from
$S_i$ by replacing $a_1a_2$ with $P_2$. 
By Lemma~\ref{Cycle}, let $Q_1, Q_2$ be independent paths in $Y_2$ from $y_2$ to $z_2, b_1$, respectively.
Then $G[\{x_1, x_2, y_1, y_2\}] \cup (z_2'Qz_1'\cup z_1'z_1x_1)\cup
(z_2'Qt\cup tx_2)\cup (z_2'z_2\cup Q_1)\cup S_1' \cup (S_2'\cup P_2\cup Q_2)$ is a  
$TK_5$ in $G'$ with branch vertices $x_1,x_2,y_1,y_2, z_2'$. 

So we may assume that such $S_1,S_2$ do not exist. Then by planarity, $K'$ has a cut $\{s_1, s_2, s_3\}$ separating $y_1$ from $\{a,a_1,a_2, z_2'\}$,
with $s_1\in V(z_2'Dz_1')$ and $s_3\in V(tDz_2')$. Clearly, $\{s_1, s_2, s_3\}$ is also a cut in $B_1$ separating $y_1$ from $\{z_2'\}\cup A$. 
Denote by $M$ the $\{s_1, s_2, s_3\}$-bridge of $B_1$ containing $y_1$. 
If $V(M)-\{s_1, s_2, s_3\}=\{y_1\}$ then $s_1=z_1'$ and $s_3=t$; now  $G[\{t,x_1, x_2,y_1\}]$ contains a  $K_4^-$ in which 
$x_1$ is of degree 2, and $(ii)$ holds. So assume  $|V(M)-\{s_1, s_2, s_3\}|\geq
2$. Then $G$ has a 6-separation $(G_1,G_2)$ such that 
$V(G_1\cap G_2)=\{s_1,s_2,s_3, x_1,x_2,z_1\}$, 
$G_2=G[M+\{x_1,x_2,z_1\}]$, and $(G_2-x_1, z_1, s_1, s_2, s_3, x_2)$
is planar. It is easy to see that $|V(G_i)|\ge 8$ for $i\in [2]$; so the assertion
follows from Lemma~\ref{6-cut2} and  then Lemma~\ref{apex}.\qed

\section{Substructure}

In this section, we derive a substructure in $G$ by finding 
five paths $A,B,C,Y,Z$ in $H:=G[B_1+\{z_1,z_2\}]$. The paths $Y,Z$ are found in the following lemma.

 \begin{figure}
\begin{center}
\includegraphics[scale=0.25]{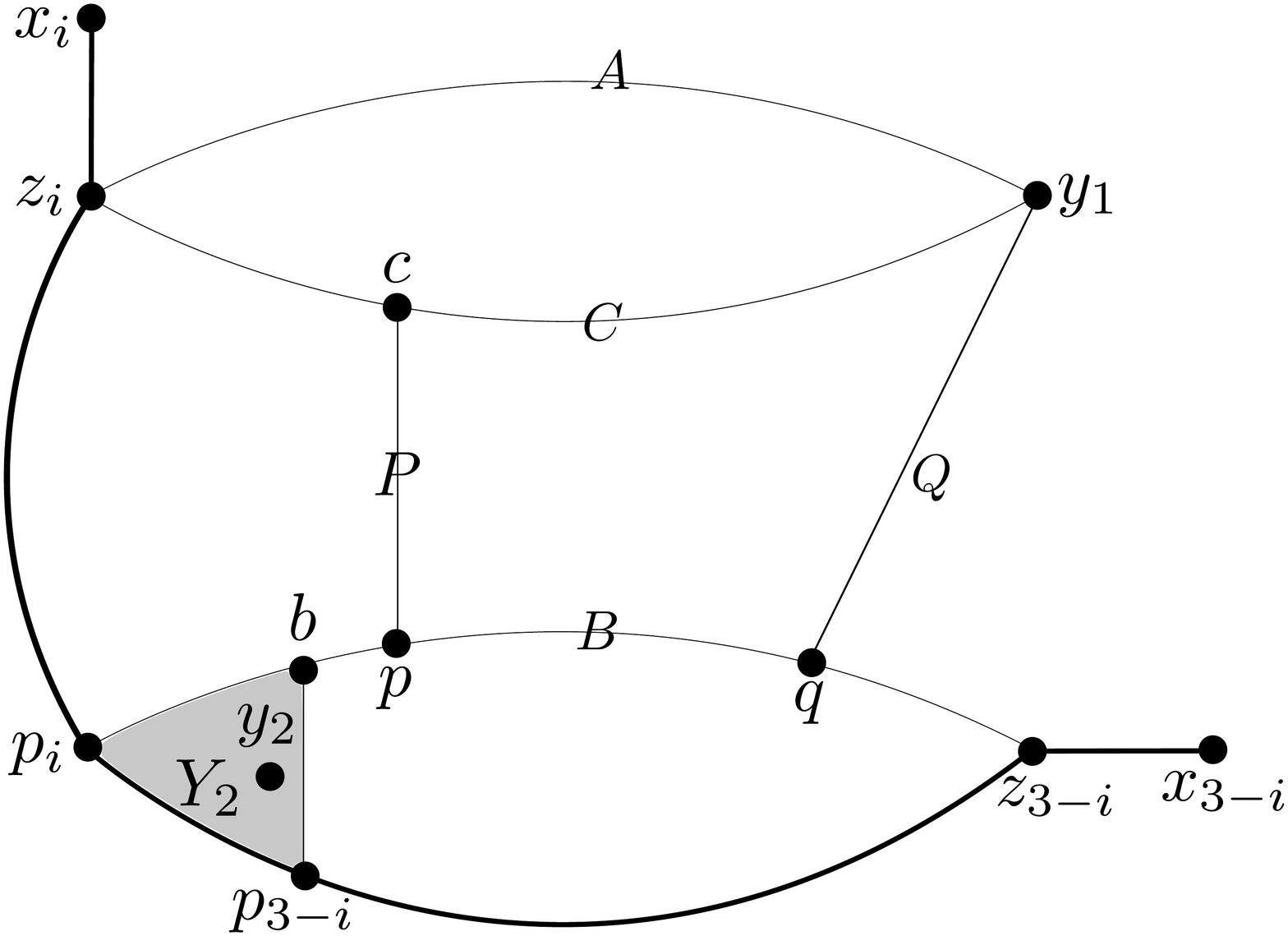}

\caption{\label{structure} An intermediate structure }
\end{center}
\end{figure}

\begin{lem}\label{YZ}
Suppose $y_2\in V(X)$ (see $(iii)$ of Lemma~\ref{classify2}), or $y_2\notin V(X)$ and $e(z_i, B_1)\ge 2$ for some $i\in [2]$ (see $(iv)$ of Lemma~\ref{classify2}).  Let 
$b_1\in N(y_2)\cap V(B_1)$ if $y_2\in V(X)$, and let $\{b_1\}=V(B_1)\cap V(B_2)$ if $y_2\notin V(X)$. Then one of the following holds:
\begin{itemize}
\item [$(i)$]  $G'$ contains $TK_5$ or $G$ contains \textsl{a} $TK_5$ in which $x_1$ is not \textsl{a} branch vertex.
\item [$(ii)$] $G-x_1$ contains $K_4^-$, or $G$ contains \textsl{a} $K_4^-$ in which $x_1$ is of degree 2. 
\item [$(iii)$] $H$ contains disjoint paths $Y,Z$ from $y_1, z_1$ to $b_1, z_2$, respectively.
\end{itemize}
\end{lem}

\pf Suppose $(iii)$ fails.  Then by Lemma~\ref{2path}, there exists a collection ${\cal A}$ of subsets of $V(H)-\{b_1,y_1, z_1, z_2\}$ 
such that $(H,{\cal A}, b_1,z_1, y_1, z_2)$ is 3-planar.

Since $B_1$ is 2-connected, $|N_H(A)\cap \{z_1, z_2\}|\le 1$ for all $A\in {\cal A}$. 
Let ${\cal A}'=\{A\in {\cal A} : |\{z_1,z_2\}\cap N_H(A)|=0\}$ and ${\cal A}''=\{A\in {\cal A} : |\{z_1,z_2\}\cap N_H(A)|=1\}$. 
Let $p(H,{\cal A})$ be the graph obtained from $H$ by deleting $A$  (for each $A\in \cal{A}$) and adding new edges joining 
every pair of distinct vertices in $N_{H}(A)$. Since $G$ is 5-connected and $B_1$ is
2-connected, $p(H,{\cal A})-\{z_1,z_2\}$ is 2-connected and we may
assume that  it is drawn in the plane with  outer cycle $D$, such that for each $A\in {\cal A}''$, the edges 
between the vertices in $N_H(A)-\{z_1, z_2\}$ occur on $D$. 

For each $j\in [2]$, let $t_j\in V(D)$ such that $H$ has a path $R_j$
from $z_j$ to $t_j$  and internally disjoint from $p(H,{\cal A})$, and, subject to this,  $t_2,b_1,t_1,y_1$ occur on $D$ in clockwise order, and 
$t_2Dt_1$ is  maximal. When $e(z_1,B_1)\ge 2$, let $t_1'\in V(b_1Dt_1)$ with $t_1'Dt_1$ maximal such that $H$ has independent paths 
$R_1,R_1'$ from $z_1$ to $t_1,t_1'$, respectively, and  
internally disjoint from $p(H,{\cal A})$. (Note that, for convenience, we use the same $R_j$ in both cases.)
When $e(z_2,B_2)\ge 2$,  let $t_2'\in V(t_2Db_1)$ with $t_2Dt_2'$
maximal such that $H$ has independent paths $R_2,R_2'$ from $z_2$ to $t_2,t_2'$, respectively, and  
internally disjoint from $p(H,{\cal A})$.

Next we define vertices $y_{21},y_{22}$ and paths $Q_1,Q_2,Q_3$. 
If $y_2\in V(X)$, then let $p_1 = p_2= b = y_2$, let $Q_j:= y_2$ for $j\in [3]$, and let 
$y_{21}, y_{22}\in N(y_2)\cap V(D)$ such that $t_2', y_{22}, y_{21}, t_1'$ occur on $D$ in clockwise order and $y_{22}Dy_{21}$ is maximal.
Suppose $y_2\notin V(X)$. By Lemma~\ref{(I)}, we may assume that for some $i\in [2]$, $e(z_i,B_1-b_1)\ge 2$ and $e(z_{3-i},B_1-b_1)\ge 1$. 
If both $e(z_1, B_1)\ge 2$ and $e(z_2,B_2)\ge 2$,
then let $y_{21}=y_{22}=b_1$ and, by  Lemma~\ref{Cycle}, 
let $Q_1,Q_2,Q_3$ be independent 
paths in $Y_2$ from $y_2$ to $p_1,p_2,b$, respectively. 
%Now assume $y_2\notin V(X)$ and
If  $e(z_{3-i}, B_1) =  1$ then
$z_{3-i}=p_{3-i}$ and, by  Lemma~\ref{Cycle}, $Y_2$ has a path $Q_{3-i}^*$ through $b, z_{3-i}, y_2, p_i$ in order;
let $R_{3-i}' := z_{3-i}Q_{3-i}^*b\cup bb_1$, $t_{3-i}' :=b_1$,
$Q_{3-i}:=y_2Q_{3-i}^*z_{3-i}$, and $Q_i:= p_iQ_{3-i}^*y_2$. 
% Let $R_{3-i}$ be a path in $H$ from $z_{3-i}$ to $t_{3-i}$ and internally disjoint from $p(H,{\cal A})$.  
Note that in this final case, $R_{3-i}$ and $R_{3-i}'$ are independent, and 
$Q_3$, $y_{21}$ and $y_{22}$ are  not defined.

Let ${\cal A}_1=\{A\in {\cal A} : z_1\in N_H(A)$ or
$N_{H}(A)\subseteq V(b_1Dy_1)\}$,  ${\cal A}_2=\{A\in {\cal A} :
z_2\in N_H(A)$ or $N_{H}(A)\subseteq V(y_1Db_1)\}$, and  $A_j=\bigcup_{A\in {\cal A}_j} A$ for $j\in[2]$.
Let $F_1:= G'[V(z_1Xp_1)\cup A_1\cup V(b_1Dy_1)]$ and $F_2:= 
G'[V(x_2Xp_2)\cup A_2\cup V(y_1Db_1)]$. Write $b_1Dy_1=v_1\ldots v_m$ 
and $z_1Xp_1=v_{m+1}\ldots v_n$ with $v_1=b_1$, $v_m=y_1$, $v_{m+1}=z_1$, and $v_n=p_1$. 
Write $y_1Db_1=u_1\ldots u_k$ and $p_2Xx_2=u_{k+1}\ldots u_l$ such that $u_1=y_1$, 
$u_k=b_1$, $u_{k+1}=p_2$ and $u_l=x_2$. 
We may assume that

\begin{itemize}
\item  [(1)]   $(F_1, v_1, \ldots, v_n)$ and $(F_2, u_1,\ldots, u_l)$ are planar.
\end{itemize}
Note that $F_1$ is $(4,\{v_1,\ldots, v_k\})$-connected, and $F_2$ is $(4,\{u_1,\ldots, u_l\})$-connected. 
We only prove that $(F_1, v_1, \ldots, v_n)$ is planar; the argument for $(F_2, u_1,\ldots, u_l)$ is similar. 
Suppose  $(F_1, v_1, \ldots, v_n)$ is not planar. Then by Lemma~\ref{society}, there exist $1\le q<r<s<t\le n$ such that 
$F_1$ contains disjoint paths $S_1, S_2$ from $v_q, v_r$ to $v_s, v_t$, respectively. 
By the definition of $F_1$ (and since $X$ is induced), we see that
$r\le m$ and $s\ge m+1$. 
Let  $T_1,T_2,T_3$ be independent paths in $B_1$ corresponding to $y_1Dt_2$, $t_2'Dv_q$, $v_rDy_1$, respectively. 
(By this, we mean that $T_1$ is between $y_1$ and $t_2$, $T_2$ is between $t_2'$ and $v_q$, and $T_3$ is between $v_r$ and $y_1$.)
Hence, $G[\{x_1,x_2,y_1,y_2\}] \cup z_2x_2\cup 
(z_2Xp_2\cup Q_2)\cup (R_2\cup T_1)\cup (R_2'\cup T_2\cup S_1\cup v_sXz_1\cup z_1x_1)\cup (T_3\cup S_2 
\cup v_tXp_1\cup Q_1)$ is a $TK_5$ in $G'$ with branch vertices $x_1,x_2,y_1,y_2, z_2$. 
This completes the proof of (1). 

\medskip
We may also assume that 
\begin{itemize}
\item [(2)]  $N_H(x_2)\subseteq V(F_2+x_1)$. 
\end{itemize}  
For, suppose there exists $a\in N_H(x_2)-V(F_2+x_1)$. 
If $a\notin A$ for all $A\in {\cal A}$ let $a'=a$ and $S=a$; and if $a\in A$ for some $A\in {\cal  A}$ then let $a'\in N_H(A)-V(F_2)$ and $S$ be a path in 
$G[A+a']$ from $a$ to $a'$. 

First, we may choose $a$ and $a'$ so that $a'\notin V(t_1Dy_1-y_1)$ and no 2-cut of $B_1$
separating $a$ from $y_1Dt_2$ is contained in $t_1Dy_1$. For,
otherwise, we may assume $a'\in V(t_1Dy_1-y_1)$ (by modifying $A$ if necessary). Let $T_1,T_2,T_3$ be independent paths in $B_1$ corresponding to $t_2'Dt_1',  t_1Da', y_1Dt_2$,
respectively. Then 
$G[\{x_1,x_2,y_2\}] \cup z_1x_1\cup z_2x_2\cup (R_1'\cup T_1\cup R_2')
\cup (z_1Xp_1\cup Q_1) \cup (z_2Xp_2\cup Q_2)\cup (R_1\cup T_2\cup
S\cup ax_2)\cup (R_2\cup T_3\cup y_1x_1)$ is a $TK_5$ in $G'$ with branch vertices $x_1, x_2, y_2, z_1, z_2$.

Suppose that  $p(H,{\cal A})-\{z_1, z_2\}-t_1Dt_2$ has a path $T$ from
$a'$ to $t_1'$. Let $T_1,T_2$ be independent paths in $B_1$ corresponding to  $T, t_1Dt_2$, respectively. 
So $G[\{x_1,x_2,y_1,y_2\}] \cup z_1x_1\cup (z_1Xp_1\cup Q_1) 
\cup (R_1\cup t_1T_2y_1)\cup (R_1'\cup T_1\cup S\cup ax_2)\cup 
(y_1T_2t_2\cup R_2\cup z_2Xp_2\cup Q_2)$  is a $TK_5$ in $G'$ with branch vertices $x_1,x_2,y_1,y_2,z_1$.

So we may assume that such $T$ does not exist. By planarity, there is a cut $\{s_1,s_2\}$
in $B_1$ separating $t_1'$ from $N_H(x_2)-V(F_2+x_1)$, with
$s_1,s_2\in V(t_1Dt_2)$. Since $\{s_1,s_2\}\not\subseteq V(t_1Dy_1)$
and $a\notin V(F_2+x_1)$, we may let  $s_1\in V(t_1Dy_1-y_1)$ and $s_2\in V(y_1Dt_2-y_1)$. Let $M$ be the $\{s_1, s_2\}$-bridge of 
$B_1$ containing $y_1$. We choose $\{s_1,s_2\}$ so that $M$ is minimal (subject to the property that 
$s_1\in V(t_1Dy_1-y_1)$ and $s_2\in V(y_1Dt_2-y_1)$). 

Since $\{s_1,s_2,x_1, x_2\}$ cannot be a cut in $G$, there 
exists  $vv'\in E(G)$ with $v'\in V(M)-\{s_1, s_2\}$ and $v\in V(z_jXp_j-z_j)$ for some $j\in [2]$.
By minimality, $M-S_j$ has independent paths $P_1, P_2$ from $y_1$ to
$s_{3-j}, v'$, respectively. Let $T_1$ be a path in $B_1-(M-s_j)$ corresponding to $t_2'Dt_1'$, and 
$T_2$ be a  path in $B_1-(M-s_j)$ corresponding to $t_1Ds_1$ (when $j=2$) or $s_2Dt_2$ (when 
$j=1$).  Then $G[\{x_1, x_2, y_1, y_2\}] \cup z_{3-j}x_{3-j}\cup (z_{3-j}Xp_{3-j}\cup Q_{3-j})\cup (R_{3-j}'\cup T_1\cup R_j'\cup z_jx_j)\cup 
(R_{3-j}\cup T_2\cup P_1)\cup (P_2\cup v'v\cup vXp_j\cup Q_j)$ is a $TK_5$ in $G'$ with branch vertices $x_1,x_2,y_1,y_2,z_{3-i}$. 

\medskip

We may assume 
\begin{itemize}
\item [(3)] $N(z_1Xp_1-z_1)\cap V(B_1)\nsubseteq V(F_1)$ or $N(z_2Xp_2-z_2)\cap V(B_1)\nsubseteq V(F_2)$. 
\end{itemize}
For, suppose $N(z_jXp_j-z_j)\cap V(B_1)\subseteq V(F_j)$ for $j\in [2]$. 
%and $N(z_2Xp_2-z_2)\cap V(B_1)\subseteq V(F_2)$. 
If $y_2\in V(X)$ then by (1) and (2), $G-x_1$ is planar; so the assertion of this lemma follows from  Lemma~\ref{apexvertex}.
Hence, we may assume $y_2\notin V(X)$. By (1) and (2), $b=b_1$, and
$(G[B_1\cup z_1Xp_1\cup p_2Xx_2], p_1, b, p_2, x_2)$ is planar. 
So $G$ has a separation $(G_1,G_2)$ such that $V(G_1\cap G_2)=\{b,p_1,p_2,x_1,x_2\}$ and 
$G_2=G[(B_1\cup z_1Xp_1\cup x_2Xp_2)+x_1]$. Clearly, $|V(G_j)| \ge
7$ for $j\in [2]$. Hence, the assertion of this lemma follows from  Lemma~\ref{apex}. 

\medskip

Since the rest of the  argument is the same  for the two cases in
(3), we will  assume 
\begin{itemize}
\item [(4)] $N(z_2Xp_2-z_2)\cap V(B_1)\nsubseteq V(F_2)$ (and, hence, $e(z_2,B_1)\ge 2$). 
\end{itemize}

Let  $vv'\in E(G)$ with $v\in V(B_1-F_2)$ and $v'\in V(z_2Xp_2-z_2)$. Let $v''=v$ and $S=v$ if $v\notin A$ for all $A\in {\cal A}$; otherwise,
let $v\in A\in {\cal A}$ and $v''\in N_H(A)$ such that $v''\notin
V(F_2)$, and let $S$ be a path in $G[A + v'']$ from $v$ to $v''$.

Suppose $(p(H,{\cal A})-\{z_1,z_2\})-t_2'Dt_1'$ has independent paths $P_1,P_2$ from $y_1$ to
$t_1, v''$, respectively. Let $P_1',P_2',T$ be independent paths in $B_1$ corresponding to $P_1,P_2,t_2'Dt_1'$, 
respectively.  
Now $G[\{x_1,x_2,y_1,y_2\}] \cup z_1x_1\cup (R_1\cup P_1') \cup (z_1Xp_1\cup Q_1)\cup (R_1'\cup T\cup R_2'\cup z_2x_2)\cup 
(P_2'\cup S\cup vv'\cup v'Xp_2\cup Q_2)$ is a $TK_5$ in $G'$ with branch vertices $x_1,x_2,y_1,y_2,z_1$. 

So we may assume that such $P_1,P_2$ do not exist in $p(H,{\cal A})$. 
Then by planarity and the existence of $t_1Dy_1$, $p(H,{\cal A})-\{z_1,z_2\}$ has a cut $\{s_1,s_2\}$, with $s_1\in V(t_2'Dt_1')$ and $s_2\in V(t_1Dy_1)$,
separating $y_1$ from $\{v'',t_1\}$. Clearly, $\{s_1,s_2\}$ is also a cut in $B_1$. 
Denote by $M_v, M_y$ the $\{s_1, s_2\}$-bridges of $B_1$ containing $\{v'', t_1\}$, $y_1$, respectively. 
We choose $\{s_1,s_2\}$ so that $M_y$ is minimal. Since $v$ is arbitrary, we have  
$N(z_2Xp_2-z_2)\cap V(B_1-F_2)\subseteq V(M_v)$. We further choose $vv'$
with $v'Xx_2$ minimal.  

\medskip
 Recall that $y_{22}$ is defined only when $y_2\in
V(X)$, or when $y_2\notin V(X)$ and both $e(z_1,B_1)\ge 2$ and
$e(z_2,B_2)\ge 2$.
We may assume
\begin{itemize}
\item [(5)] $y_{22}\in V(M_v)$ (when defined) and, for any $q\in V(p_2Xv'-v')$, $N(q)\cap V(M_y-\{s_1, s_2\})= \emptyset$.
\end{itemize} 
Suppose (5) fails. If $y_{22}$ is defined and $y_{22}\notin V(M_v)$
let $q=b$, $q'= y_{22}$, and $Q'=q'q\cup Q_3$;  and if $y_{22}$ is
defined, $y_{22}\in V(M_v)$, and there exist  $q\in V(p_2Xv'-v')$ and  $q'\in N(q)\cap V(M_y-\{s_1, s_2\})$, then let
$Q'=q'q\cup qXp_2\cup Q_2$.  

Since $B_1$ is 2-connected, there exists $j\in [2]$ such that $M_v-s_{3-j}$ contains disjoint paths $T_1, T_2$ from $\{t_1, t_1'\}$ to $\{v'', s_j\}$. 
Note that $R_1\cup R_1'\cup T_1\cup T_2$ contains independent paths $T_1', T_2'$ from $z_1$ to $v'', s_j$, respectively. 
If $M_y$ contains independent paths $S_1, S_2$ from $y_1$ to $q',
s_j$, then $G[\{x_1,x_2,y_1,y_2\}] \cup z_1x_1\cup (z_1Xp_1\cup Q_1)\cup (T_1'\cup S\cup vv'\cup v'Xx_2)\cup (T_2'\cup S_2) \cup 
(Q'\cup S_1)$ is a $TK_5$ in $G'$ with branch vertices $x_1,x_2,y_1,y_2,z_1$. 
So we may assume $S_1, S_2$ do not exist in $M_y$; hence $M_y$ has a cut vertex $c$ that separates $y_1$ from $\{q', s_j\}$.  
%Thus,we must have $j=1$; so 
%$M_v-s_1$ has no disjoint paths from $\{t_1,t_1'\}$ to $\{s_2,v''\}$.

By the minimality of $M_y$ and the existence of $y_1Ds_1$, $c\in
V(y_1Dt_2'-t_2')$; so we must have $j=1$. 
Denote by $C_q$, $C_y$ the $c$-bridges of $M_y$ containing  $\{q', s_1\}$, $y_1$, respectively, and choose $c$ with $C_y$ minimal. 
Then $N(p_2Xv'-v')\cap V(C_y-\{c,s_2\})=\emptyset$. 

We may assume that  there exist $uu'\in E(G)$ with $u'\in V(z_1Xp_1-z_1)$ and $u\in V(C_y)-\{c,s_2\}$. For, otherwise, by (1) and (2), there exists
$z\in V(v'Xx_2)$ such that $\{c,s_2,x_1,x_2,z\}$ is a cut in $G$, and $G$ has a separation $(G_1,G_2)$ 
such that $V(G_1\cap G_2)=\{c,s_2,x_1,x_2,z\}$, $M_v\cup z_1Xz\cup Y_2\subseteq G_1$, 
$M_y\subseteq G_2$, and $(G_2-x_1,\{c,s_2,x_2,z\})$ is planar. Clearly, $|V(G_1)|\ge 7$ and $|V(G_2)|\ge 6$. 
If $|V(G_2)|\ge 7$ then the assertion of the lemma follows from Lemma~\ref{apex}. 
So assume $|V(G_2)|=6$. Then $z=z_2$ and $y_1z_2\in E(G)$; now $G[\{x_1,x_2,y_1,z_2\}]-x_1z_2\cong K_4^-$ in which $x_1$ is of degree 2, and $(ii)$ holds. 

By the minimality of $M_y$ and $C_y$, $C_y-s_2$ has  independent paths $U_1, U_2$ from $y_1$ to $c, u$, respectively.
In $M_v-s_1$, we find a path $T$ from $t_1$ to $v''$. Let $X^*$ be an induced path 
in $G-x_1$ from $z_1$ to $x_2$ such that $V(X')\subseteq V(R_1\cup T\cup S\cup vv'\cup v'Xx_2)$.
Now $U_1\cup U_2\cup (C_q-s_1)\cup uu'\cup u'Xp_1\cup Q_1\cup Q_2\cup p_2Xq\cup qq'$ is a subgraph of $(G-x_1)-X^*$ and  
has a cycle containing $\{y_1,y_2\}$. Hence by Lemma~\ref{4A} and Lemma~\ref{apex}, 
we may assume that $G-x_1$ contains an induced  path $X'$ from $z_1$
to $x_2$ such that $y_1,y_2\notin V(X')$ and $(G-x_1)-X'$ is
2-connected. So the assertion of this lemma follows from Lemma~\ref{2-connected}. 
This proves (5).

\medskip
We may assume $N(z_1Xp_1-z_1)\cap V(M_y-\{s_1, s_2\})\ne \emptyset$. For, otherwise, by (5),  
$G$ has a 5-separation $(G_1, G_2)$ such that $V(G_1\cap
G_2)=\{s_1,s_2,v',x_1,x_2\}$, $G_2:=G[v'Xx_2\cup M_y+x_1]$ and
$(G_2-x_1, s_1, s_2, x_2, v')$ is planar (by (1) and (2)). Clearly, $|V(G_1)|\ge 7$ and $|V(G_2)|\ge 6$. If $|V(G_2)|\ge 7$ then 
the assertion of this lemma follows from Lemma~\ref{apex}. So assume $|V(G_2)|=6$. Then $v'=z_2$ and $y_1z_2\in E(G)$. 
So $G[\{x_1,x_2,y_1,z_2\}]-x_1z_2\cong K_4^-$ in which $x_1$ is of degree 2, and $(ii)$ holds.

So there exists  $uu'\in E(G)$ with $u'\in V(z_1Xp_1-z_1)$ and $u\in
V(M_y)-\{s_1, s_2\}$. Hence, $e(z_1,B_1)\ge 2$; so $y_{21},y_{22},Q_3$ are defined.  Let $P_u$ be a path in $M_y$ 
from $u$ to some $u_D\in V(s_2Ds_1)-\{s_1,s_2\}$ and 
internally disjoint from $V(D)$ (which exists 
by minimality of $M_y$), and $P_v$ be a path in $M_v$ from $v''$ to some $v_D\in V(s_1Ds_2)$ and internally disjoint from $V(D)$. 
By the definition of $F_2$, we may choose $v_D$ so that $v_D\notin V(s_1Dy_{22})$. 

We may assume $v_D\in V(t_1'Dy_1-t_1')$. For, suppose $v_D\in
V(y_{22}Dt_1'-y_{22})$. Let $T_1,T_2,T_3$ be independent paths in $B_1$ corresponding to 
$t_1Dy_1$, $v_DDt_1'$, $y_1Dy_{22}$, respectively. Then 
 $G[\{x_1,x_2,y_1,y_2\}] \cup z_1x_1\cup (z_1Xp_1\cup Q_1)\cup (R_1\cup T_1)\cup (R_1'\cup T_2\cup P_v\cup S\cup vv'\cup v'Xx_2)\cup (T_3\cup y_{22}b\cup Q_3)$ 
is a $TK_5$ in $G'$ with branch vertices $x_1,x_2,y_1,y_2,z_1$.  

Next, we consider the location of $u_D$. 
Suppose $u_D\in V(t_2'Ds_1-s_1)$.  Let $T_1,T_2,T_3$ be independent paths in $B_1$ corresponding to $y_1Dt_2$, $t_2'Du_D$, $y_{21}Dy_1$, respectively. 
Then $G[\{x_1,x_2,y_1,y_2\}] \cup z_2x_2\cup (z_2Xp_2\cup Q_2)\cup (R_2\cup T_1)\cup 
(R_2'\cup T_2\cup P_u\cup uu'\cup u'Xz_1\cup z_1x_1)\cup (T_3\cup y_{21}b\cup Q_3)$ is a $TK_5$ in $G'$ with branch vertices $x_1,x_2,y_1,y_2,z_2$. 

Now suppose $u_D\in V(s_2Dy_1)$.  Let $T_1,T_2,T_3$ be independent paths in $B_1$ corresponding to $y_1Dt_2$, $t_2'Dt_1'$, $u_DDy_1$, respectively. 
Then $G[\{x_1,x_2,y_1,y_2\}] \cup z_2x_2\cup (z_2Xp_2\cup Q_2)\cup (R_2\cup T_1)\cup (R_2'\cup T_2\cup R_1'\cup z_1x_1)\cup 
(T_3\cup P_u\cup uu'\cup u'Xp_1\cup Q_1)$ is a $TK_5$ in $G'$ with branch vertices $x_1,x_2,y_1,y_2,z_2$. 

So we may assume $u_D\in V(y_1Dt_2'-t_2')$.  Let $T_1,T_2,T_3$ be
independent  paths in $B_1$ corresponding to $y_1Du_D$, $t_2'Dt_1'$, $v_DDy_1$, respectively. 
Thus, $(G-x_1)-(R_1'\cup T_2\cup R_2'\cup z_2x_2)$ contains the cycle $T_1\cup P_u\cup uu'\cup u'Xp_1\cup Q_1\cup Q_2\cup p_2Xv'\cup v'v\cup S\cup P_v\cup T_3$. 
Hence, by Lemma~\ref{4A} and Lemma~\ref{apex}, we may assume that $G-x_1$ contains a path $X'$ from $z_1$ to $x_2$ such that $y_1,y_2\notin V(X')$ and $(G-x_1)-X'$ is 2-connected. 
So the assertion of this lemma follows from Lemma~\ref{2-connected}. \qed

\medskip

We now prove the existence of three paths $A,B,C$  in $H:=G[B_1+\{z_1,z_2\}]$.
%such that two of these paths are from $z_i$ (for some $i\in [2]$) to $y_1$ and one from $b$ to $z_{3-i}$. 

\begin{lem}\label{ABC} 
Let $b_1\in N(y_2)\cap V(B_1)$ when $y_2\in V(X)$, and let $\{b_1\}=V(B_1)\cap V(B_2)$ when $y_2\notin V(X)$. Then one of the following holds: 
\begin{itemize}
\item [$(i)$] $G'$ contains $TK_5$, or $G$ contains \textsl{a} $TK_5$ in which $x_1$ is not \textsl{a} branch vertex.
\item [$(ii)$] $G-x_1$ contains $K_4^-$, or $G$ contains \textsl{a} $K_4^-$ in which $x_1$ is of degree 2.
\item  [$(iii)$] There exists $i\in [2]$ such that $H$ contains independent paths $A,B,C$, with $A$ and $C$ from $z_i$ to $y_1$ 
and $B$ from $b_1$ to $z_{3-i}$. 
\end{itemize}
\end{lem}

\pf If $y_2\notin V(X)$ then by Lemma~\ref{Cycle}, let $Q_1,Q_2,Q_3$ be independent paths in
$Y_2$ from $y_2$ to $p_1,p_2,b$, respectively. When $y_2\in
V(X)$ let $Q_1=Q_2=Q_3=y_2$.  We may assume that 
\begin{itemize}
\item [(1)] for $i\in [2]$, $H$ has no path through $z_{3-i}, z_i, y_1, b_1$ in order.
\end{itemize}
For, if $H$ has a path $S$ through $z_{3-i}, z_i, y_1, b_1$ in
order. Then $G[\{x_1, x_2, y_1, y_2\}] \cup z_ix_i\cup (z_iXp_i\cup Q_i)\cup z_iSy_1\cup 
(z_iSz_{3-i}\cup z_{3-i}x_{3-i})\cup (y_1Sb_1\cup b_1b\cup Q_3)$ is a $TK_5$ in $G'$ with branch vertices $x_1,x_2,y_1,y_2,z_i$.

\medskip

We may also assume that 
\begin{itemize}
\item [(2)] for $i\in [2]$ with $e(z_i,B_1-b_1)\ge 2$,  $H$ has a
  2-separation $(F_i',F_i'')$ such that  $b_1\in V(F_i')$, $z_i\in V(F_i'-F_i'')$ and $\{y_1, z_{3-i}\}\subseteq V(F_i''-F_i')$. 
\end{itemize}
Suppose  $i\in [2]$ and $e(z_i,B_1-b_1)\ge 2$. Let $K$ be obtained from $H$ by duplicating $z_i$ and
$y_1$ with copies $z_i'$ and $y_1'$, respectively. So in $K$, $y_1$ and $y_1'$ 
are not adjacent, but have the same set of neighbors, namely $N_H(y_1)$; and the same holds for $z_i$ and $z_i'$.

Suppose $K$ contains disjoint paths $A',B',C'$ from $\{z_i,z_i',
b_1\}$ to $\{y_1,y_1', z_{3-i}\}$, with $z_i\in V(A'), z_i'\in V(C')$
and 
$b_1\in V(B')$. If $z_{3-i}\notin V(B')$ then, after identifying $y_1$ with $y_1'$ and
$z_i$ with $z_i'$, we obtain from $A'\cup B'\cup C'$ a path in $H$
from $z_{3-i}$ to $b_1$ through $z_i,y_1$ in order, contradicting (1). Hence $z_{3-i}\in V(B')$, 
and we get the desired paths for $(iii)$ from $A'\cup B'\cup C'$, by identifying $y_1$ with $y_1'$ and $z_i$
with $z_i'$. 

So we may assume that such $A', B', C'$ do not exist. Then $K$ has a separation $(K',K'')$ such that $|V(K'\cap K'')|\le 2$, 
$\{b_1,z_i, z_i'\}\subseteq V(K')$ and $\{y_1,y_1', z_{3-i}\}\subseteq V(K'')$. Since $H-z_{3-i}$ is 
2-connected, $z_{3-i}\notin  V(K'\cap K'')$.  

We claim that $z_i,z_i'\notin  V(K'\cap K'')$. For, if  exactly one of $z_i,z_i'$ is in $V(K'\cap
K'')$ then, since $z_i,z_i'$ have  the same set of neighbors in $K$, $V(K'\cap K'')-\{z_i,z_i'\}$ is a cut in $H$
separating $\{z_{3-i}, y_1\}$ from $\{z_i, b_1\}$, a contradiction. Now assume 
$\{z_i,z_i'\}= V(K'\cap K'')$.  Then $z_i$ is a cut vertex in $H$ separating $b_1$ from $\{y_1, z_{3-i}\}$, a contradiction.

We may assume that $y_1,y_1'\notin V(K'\cap K'')$. First, suppose exactly  one of $y_1,y_1'$ is in $V(K'\cap
K'')$. Then, since $y_1,y_1'$ have the same set of neighbors in $K$, $V(K'\cap K'')-\{y_1,y_1'\}$ is a cut  in $H$
separating $\{z_{3-i}, y_1\}$ from $\{z_i, b_1\}$, a contradiction. Now assume 
$\{y_1,y_1'\}= V(K'\cap K'')$.  Then $y_1$ is a cut vertex in $H$ separating $z_{3-i}$ from $\{b_1, z_i\}$. 
This implies that $N(z_{3-i})\cap V(B_1)=\{y_1\}$; so $y_2\notin V(X)$
and $z_{3-i}=p_{3-i}$. We may assume $i=2$; for otherwise,
$G[\{x_1,x_2,y_1,z_2\}]-x_1z_2\cong K_4^-$ in which $x_1$ is of degree
2, and $(ii)$ holds.  
Then $z_1=p_1$,  and $G$ has a 5-separation $(G_1,G_2)$ such that  
$V(G_1\cap G_2)=\{b,p_2,x_1,x_2,y_1\}$ and $G_2=G[B_1\cup x_2Xp_2+\{b,x_1\}]$. 
Note that $x_1x_2y_1x_1$ is a triangle and $|V(G_j)|\ge 7$ for $j\in
[2]$.  So the assertion of this lemma follows from Lemma~\ref{5cut_triangle}.

Thus, since $B_1$ is 2-connected, $|V(K'\cap K'')|=2$.  Let $V(K'\cap K'')=\{s,t\}$, and let $F_i'$ (respectively, $F_i''$) be 
obtained from $K'$ (respectively, $K''$) by identifying $z_i'$ with
$z_i$ (respectively, $y_1'$ with $y_1$).  
Then $(F_i',F_i'')$ gives the desired  2-separation in $H$, completing the proof of (2).

\medskip
By Lemma~\ref{YZ}, we may assume that 
\begin{itemize}
\item [(3)] $H$ has disjoint paths $Y,Z$ from $y_1,z_1$ to $b_1,z_2$, respectively. 
\end{itemize}

We now consider three cases.

\medskip
{\it Case} 1. $e(z_i,B_1-b_1)\ge 2$ for $i\in [2]$. 

For $i\in [2]$, let $V(F_i'\cap F_i'')=\{s_i,t_i\}$ as in (2). 
Let $Z_1,B_1'$ denote the $\{s_1,t_1\}$-bridges of $F_1'$ containing
$z_1,b_1$, respectively, and let $Y_1,Z_2$ denote the $\{s_1,t_1\}$-bridges of $F_1''$
containing $y_1,z_2$, respectively. 

Suppose $Y_1\ne Z_2$, and suppose $Z_1\ne B_1'$ or $b_1\in \{s_1,t_1\}$.  Let $b_1=s_1$ if $b_1\in \{s_1,t_1\}$. 
Then $Z_1$ has independent paths $S_1,T_1$ from $z_1$ to $s_1,t_1$,
respectively. Moreover, 
$Z_2$ has independent paths $S_2,T_2$ from $z_2$ to $s_1,t_1$, respectively, $B_1'-t_1$ has a path $P$ from $s_1$ to $b_1$, and $Y_1$ has 
independent paths $S_3,T_3$ from $y_1$  to $s_1,t_1$, respectively. 
So $x_1z_1\cup (z_1Xp_1\cup Q_1) \cup x_1y_2\cup (z_2Xp_2\cup Q_2)\cup z_2x_2x_1 \cup (T_2\cup T_1)
\cup S_1\cup S_2\cup (S_3\cup y_1x_1)\cup (P\cup b_1b\cup Q_3)$ is a $TK_5$ in $G'$ with branch vertices $s_1,x_1,y_2,z_1,z_2$. 

Thus, we may assume that $Y_1=Z_2$, or $Z_1=B_1'$ and $b_1\notin \{s_1,t_1\}$. 
First, suppose  $Y_1\ne Z_2$. Then $Z_1=B_1'$ and $b_1\notin \{s_1,t_1\}$, and hence $B_1'-\{s_1,t_1\}$ has 
a path from $z_1$ to $b_1$. Since $H$ is 2-connected, $Y_1\cup Z_2$ has two independent paths 
from $y_1$ to  $z_2$. However, this contradicts the existence of the
separation $(F_2',F_2'')$.

So $Y_1=Z_2$. Thus, by symmetry, we may assume $t_2\in V(Y_1)-\{s_1,t_1\}$.
 Suppose  $b_1\notin \{s_1,t_1\}$ and $B_1'=Z_1$. Then $s_2\in V(B_1')-\{s_1,t_1\}$.  Moreover, $\{s_2,t_2\}$ separates $s_1$ from $t_1$ in $H$; for otherwise, either 
$t_2$ separates $z_2$ from $\{b_1,y_1,z_1\}$ in $H$, or $t_2$ separates $y_1$ from $\{b_1,z_1,z_2\}$ in $H$, a contradiction. Thus, we may 
assume that in $H$, $\{s_2,t_2\}$ separates $\{b_1,s_1,z_2\}$ from $\{t_1,y_1,z_1\}$. However, this contradicts (3). 

Therefore, $B_1'\ne Z_1$ or $b_1\in \{s_1,t_1\}$.  If $b_1\notin
\{s_1,t_1\}$ then $B_1'\ne Z_1$; so  $s_2\in \{s_1,t_1\}$ (because of $(F_2',F_2'')$), 
and we may assume $s_2=s_1$. If $b_1\in \{s_1,t_1\}$ then we may assume that $b_1=s_1$; so $s_2=s_1$ or, in $Z_1$,  
$s_2$ separates $s_1$ from $\{t_1,z_1\}$. 
Let $Y_1',Z_2'$ be the $t_2$-bridges of $Y_1-\{s_1,t_1\}$ containing $y_1,z_2$, respectively. Again, because of the existence of $(F_2',F_2'')$, $t_1$ has no 
neighbor in $Z_2'-t_2$. Hence, by (3),  $s_1$ has a neighbor in $Y_1'-t_2$; and, thus, $s_2=s_1$ and $G[Y_1'+\{s_1,t_1\}]$ has disjoint paths 
$S_1,T_1$ from $s_1,t_1$ to $y_1,t_2$, respectively. Let $S_2,T_2$ be independent paths in $G[Z_2'+s_1]$ from $z_2$ to $s_1,t_2$, respectively, and 
$S,T$ be independent paths in $Z_1$ from $z_1$ to $s_1,t_1$, respectively. Let $P$ be a path in $B_1'-t_1$ from $s_1$ to $b_1$. Then 
$x_1z_1\cup (z_1Xp_1\cup Q_1)\cup x_1y_2\cup (z_2Xp_2\cup Q_2)\cup z_2x_2x_1\cup  (T_2\cup T_1\cup T)\cup S\cup (S_1\cup y_1x_1)\cup S_2\cup 
(P\cup b_1b\cup Q_3)$ is a $TK_5$ in $G'$ with branch vertices $s_1,x_1,y_2,z_1,z_2$. 

\medskip

{\it Case} 2. $e(z_2, B_1-b_1)\ge 2$. 

If $y_2\in V(X)$ then $e(z_1,B_1-b_1)\ge 2$, and if $y_2\notin V(X)$
then, by Lemma~\ref{(I)}, $e(z_1, B_1-b_1)\ge 1$. In view of Case 1,
we may assume $e(z_1, B_1-b_1)=1$; so $z_1=p_1$ and $y_2\notin V(X)$. 
Note that if $b\ne b_1$ then, by Lemma~\ref{(I)'}, we may assume $z_1b_1\in
E(G)$; so $b_1\in V(F_2'\cap F_2'')$.
By Lemma~\ref{Cycle}, we may assume that $Y_2$ has a path $Q$ from $p_2$ to $b$ through $y_2,z_1$ in this order.

For convenience, let $F':=F_2'$, $F'':=F_2''$, $s:=s_2$ and
$t:=t_2$. So $b_1,z_2\in V(F')$ and $y_1,z_1\in V(F'')$. 
We choose $(F',F'')$ so that $F''$ is minimal.
Let $z_1'$ denote the unique neighbor of $z_1$ in $B_1-b_1$.  

\medskip

{\it Subcase} 2.1. $N(z_2Xp_2-z_2)\cap V(F''-\{z_1, s, t\})\not\subseteq \{z_1'\}$. 

Let $uu'\in E(G)$, with $u\in V(F'')-\{z_1, z_1', s, t\}$ and $u'\in V(z_2Xp_2-z_2)$. 
Note that $F'$ contains a path $S$ from $z_2$ to $b_1$ such that
$|V(S)\cap \{s, t\}|\leq 1$. Moreover, if there exists $r\in \{s,t\}$ such that $r\in V(S)$ for all 
such path $S$, then $b_1=r$. 

 If $(F''-z_1)-S$ contains independent paths $T_1, T_2$ 
from $y_1$ to $z_1', u$, respectively, then $G[\{x_1, x_2, y_1,
y_2\}] \cup z_1x_1\cup z_1Qy_2\cup (z_1Qb\cup bb_1\cup S\cup z_2x_2)
\cup (z_1z_1'\cup T_1)\cup (T_2\cup uu'\cup u'Xp_2\cup p_2Qy_2)
$ is a $TK_5$ in $G'$ with branch vertices $x_1,x_2,y_1,y_2,z_1$.

 So we may assume that such $T_1, T_2$ do not exist.
Hence, there is a cut vertex $c$ in $(F''-z_1)-S$ separating $y_1$ from $\{u,z_1'\}$. 
Denote by $M_1, M_2$ the $(\{c\}\cup (V(S)\cap \{s,t\}))$-bridges of $F''-z_1$ containing $y_1$, $\{u, z_1'\}$, respectively. 
We may choose $c$ so that $M_1$ is minimal. Then $N(z_2Xp_2-z_2)\cap V(F'')\subseteq V(M_2)$ (as $uu'$ was chosen arbitrarily).

Since $G$ is $5$-connected, $\{s, t\}\subseteq V(M_1)$ (as otherwise $\{c,x_1,x_2\}\cup (\{s,t\}\cap V(M_1))$ would be a cut in $G$),  
 and $M_1$ contains independent paths $R_1, R_2, R_3$ from $y_1$ to $c, s, t$, respectively. 
Since $B_1$ is 2-connected, $\{s, t\}\cap V(M_2)\neq \emptyset$ and
there exist choices of $u$ and $r\in \{s, t\}\cap V(M_2)$ 
such that $M_2$ contains disjoint paths $R_4, R_5$ from $\{z_1', u\}$
to $\{c, r\}$ and avoiding $\{s,t\}\cap V(M_2)-\{r\}$. Thus, $R_1\cup R_2\cup R_3\cup R_4\cup R_5$ contains 
independent paths from $y_1$ to $z_1', u$, respectively, and avoiding $\{s,t\}\cap V(M_2)-\{r\}$. By the
non-existence of $T_1$ and $T_2$, $r\in V(S)$ for every choice of
$S$. Hence, $b_1=r$, $\{s, t\}\cap V(M_2)= \{r\}$, and $V(S)\cap \{s, t\}=\{r\}$ for every choice of $S$. 
Without loss of generality, we may assume that $r = t$. 

We further choose $uu'$ so that $u'Xp_2$ is maximal. 
Suppose $N(u'Xp_2-u')\cap V(F'-\{s, t\})=\emptyset$. Then $G$ has a 5-separation $(G_1,G_2)$ 
such that $V(G_1\cap G_2)=\{s,t,u',x_1, x_2\}$ and $G_2= G[F'\cup
x_2Xu'+x_1]$.  Clearly, $|V(G_1)|\ge 7$. Since $e(z_2, B_1-b_1)\geq
2$, $|V(G_2)|\geq 7$.  If   $(G_2-x_1, x_2, s, t, u')$ is planar then
the assertion of this lemma follows from Lemma~\ref{apex}. 
Hence, we may assume, by Lemma~\ref{2path}, that $G_2-x_1$ contains disjoint paths $X_1, X_2$ 
from $u', x_2$ to $s, t$, respectively. Let $X_3$ be a path in $M_2-t$
from $z_1'$ to $c$. Then $G[\{x_1, x_2, y_1, y_2\}] \cup z_1x_1\cup
z_1Qy_2\cup (z_1Qb\cup bb_1\cup X_2)\cup 
(z_1z_1'\cup X_3\cup R_1)\cup (R_2\cup X_1\cup u'Xp_2\cup p_2Qy_2)$ is a $TK_5$ in $G'$ with branch vertices $x_1,x_2,y_1,y_2,z_1$. 

So assume that there exists  $ww'\in E(G)$ with $w'\in V(u'Xp_2-u')$ and $w\in V(F'-\{s, t\})$. 
Let $S_1$ be a path in $F'-t$ from $w$ to $s$ and $S_2$ be a path in $M_2-t$ from $z_1'$ to $u$.
Then $G[\{x_1, x_2, y_1, y_2\}] \cup z_1x_1\cup z_1Qy_2\cup (z_1Qb\cup
bb_1\cup R_3)\cup (z_1z_1'\cup S_2\cup uu'\cup u'Xx_2)\cup (R_2\cup S_1\cup ww'\cup w'Xp_2\cup p_2Qy_2)$ is a $TK_5$ in $G'$ 
with branch vertices $x_1,x_2,y_1,y_2,z_1$.

\medskip 

{\it Subcase} 2.2. $N(z_2Xp_2-z_2)\cap V(F''-\{z_1, s, t\})\subseteq \{z_1'\}$.

Then $\{s, t,x_1,x_2, z_1'\}$ is a $5$-cut in $G$ separating $F''-z_1$ from $F'\cup Y_2$. 
Since $G$ is 5-connected,  $F''-z_1$ has independent paths $T_1, T_2, T_3$ from $y_1$ to $s, t, z_1'$, respectively. 
Next, we find a path $R$ in $F''-z_1$ from  $s$ to $t$ and containing
$\{y_1,z_1'\}$. For this, let $F_g:=(F''-z_1)+\{g, gs, gt\}$, where $g$ is a new vertex. Since
$G$ is 5-connected and we are in Subcase 2.2, 
$F_g$ has no $2$-cut separating $y_1$ from $\{g, z_1'\}$. Hence, 
by Lemma~\ref{Watkins}, there is a cycle in $F_g$ containing
$\{g,y_1,z_1'\}$ and, after removing $g$ from this cycle, we get the
desired $R$.

Let $x=p_2$ if $N(z_2Xp_2-z_2)\cap V(F''-\{z_1,s,t\})=\emptyset$ and, otherwise, let $x\in N(z_1')\cap N(z_2Xp_2-z_2)$ with  $xXz_2$ minimal. 

We may assume that $N(xXp_2-x)\cap
V(B_1-\{b_1,z_1'\})=\emptyset$. For, otherwise, there exists  $rr'\in E(G)$ such that  $r\in V(B_1)-\{b_1, z_1'\}$ and $r'\in
V(xXp_2-x)$. Then $r\in V(F')$ and $x\ne p_2$; so $xz_1'\in E(G)$. Note that $F'$ 
has disjoint paths from $\{s,t\}$ to $\{b_1,r\}$, which, combined with $T_1,T_2$, gives independent paths $P_1,P_2$ 
in $B_1-z_1'$ from $y_1$ to $b_1,r$, respectively. Hence, in $(G-x_1)-(z_1z_1'x\cup xXx_2)$, 
$\{y_1,y_2\}$ is contained in the cycle $P_1\cup P_2\cup rr'\cup r'Xp_2\cup Q_2\cup
Q_3\cup bb_1$. Hence, by Lemma~\ref{4A} and Lemma~\ref{apex}, we
may assume that $G-x_1$ has a path $X'$ from $z_1$ to $x_2$ such that 
$y_1,y_2\notin V(X)$, and $(G-x_1)-X'$ is 2-connected. Thus, the
assertion of this lemma follows from Lemma~\ref{2-connected}. 

We may assume  $b=b_1$. For, suppose $b\ne b_1$. Then, using the
notation from $(iv)$ of Lemma~\ref{classify2}, $v\in V(p_2Xx_2-p_2)$
and $b_1'\in V(B_1-b_1)$. Let $P_1,P_2$ be independent paths in $B_1$
from $y_1$ to $b_1, b_1'$, respectively. Then $G[\{x_1, x_2,
y_1,y_2\}]\cup z_1x_1\cup z_1Qy_2\cup (z_1Qb\cup bb_1\cup P_1)\cup
(z_1Qb\cup bv\cup vXx_2)\cup (P_2\cup b_1'p_2\cup p_2Qy_2)$
is a $TK_5$ in $G'$ with branch vertices $x_1,x_2,y_1,y_2,z_1$.

Therefore, $G$ has a separation $(G_1,G_2)$ 
such that $V(G_1\cap G_2)=\{b_1,s,t,x,x_1,x_2\}$ and $G_2=G[F'\cup xXx_2+x_1]$. Let $G_2'=G_2+\{r,rs,rt\}$, where $r$ is a new vertex. 

We may assume that $(G_2'-x_1,{\cal A}, b_1,x,x_2,r)$ is 3-planar for some collection ${\cal A}$  of subsets of $V(G_2'-x_1)-\{b_1,x,x_2,r\}$. 
For, otherwise, by Lemma~\ref{2path}, $G_2'-x_1$ contains disjoint paths $R_1,R_2$ from $b_1,x$ to $x_2,r$, respectively. Let $R=T_2\cup (R_2-r)$ if $R_2-r$ ends at $t$,
and $R=T_1\cup (R-r)$ otherwise.  Then $G[\{x_1, x_2, y_1,
y_2\}]\cup z_1x_1\cup z_1Qy_2\cup (z_1Qb_1\cup R_1)\cup (z_1z_1'\cup T_3)\cup (R\cup xXp_2\cup p_2Qy_2)$
is a $TK_5$ in $G'$ with branch vertices $x_1,x_2,y_1,y_2,z_1$.

We choose ${\cal A}$ to be minimal and define $J, s',t'$ as follows. 
If ${\cal A}=\emptyset$ then after relabeling of $s,t$ (if necessary), we may assume  
$(G_2'-x_1,b_1,x,x_2,s,t)$ is planar and let $J=G_2$, $s'=s$ and $t'=t$.
Now assume ${\cal A}\ne \emptyset$. 
Then, by the minimality of ${\cal A}$ and 5-connectedness of $G$, 
${\cal A}$ has a unique member, say $A$, such that $r\in N(A)$ and $\{s,t\}\subseteq A$ and, moreover, $G'[A\cup \{s',t'\}]$ is connected,
where $N(A)\cap V(F')=\{r,s',t'\}$. Let $J$ denote the $\{s',t',x_1\}$-bridge of $G_2'$ containing $\{b_1,x,x_2\}$. 
We may assume, after suitable labeling
of $s',t'$, $(J-x_1,b_1,x,x_2,s',t')$ is planar.

Suppose $b_1\in \{s',t'\}$.  Then  $G$ has a 5-separation $(L_1,L_2)$ such that $V(L_1
\cap L_2)=\{s',t',x,x_1,x_2\}$ and $L_2=J$. If $|V(J)|\ge 7$ then the assertion of
this lemma follows from Lemma~\ref{apex}. So assume $|V(J)|\le 6$. 
Since $e(z_2,B_1-b_1) \ge 2$, there exists $v\in N(z_2)\cap   V(F'-\{s',t',z_2\})$.  Since $G$ is 5-connected, $vx_1,vx_2\in E(G)$. Hence, 
$G[\{v,x_1,x_2,z_2\}]$ contains a $K_4^-$ in which $x_1$ is of degree 2.

Thus, we may assume that $b_1\notin \{s',t'\}$. Then  $G$ has a 6-separation $(L_1,L_2)$ such that $V(L_1
\cap L_2)=\{b_1,s',t',x,x_1,x_2\}$ and $L_2=J$. 
 If $|V(J)|\ge 8$ then the assertion of this lemma follows from Lemmas~\ref{6-cut2} and \ref{apex}. 

So assume $|V(J)|\le 7$. By planarity of $J$ and 2-connectedness of
$B_1$, $z_2t'\notin E(G)$. Thus, since $e(z_2,B_1-b_1) \ge 2$,
$z_2s'\in E(G)$ and there exists $v\in V(J)-\{b_1,s',t',x,x_2,z_2\}$ such
that $z_2v\in E(G)$. So $|V(J)|=7$ and $z_2=x$. By the minimality of $F'$, $vt'\in E(G)$; and by the
2-connectedness of $B_1$,  $\{vs',vb_1\}\subseteq E(G)$. By planarity
of $J$, $x_2v\notin E(G)$.
% as otherwise $G[\{s',v,x_2,z_2\}]$ (and, hence, $G-x_1$) contains a
% $K_4^-$ and $(ii)$ holds. 
Thus, $vx_1\in E(G)$ as $G$ is 5-connected. Then we may assume
$x_1b_1\notin E(G)$;  for otherwise $G[\{b_1,t',v,x_1\}]-x_1t'\cong
K_4^-$ in which $x_1$ is of degree 2, and $(ii)$ holds. We may also assume $x_1z_2\notin
E(G)$; for otherwise $G[\{s',v,x_1,z_2\}]-x_1s'\cong K_4^-$ in which
$x_1$ is of degree 2, and $(ii)$ holds. So  $z_2=p_2$ as $G$ is
5-connected. 
%otherwise, $z_2x_1\not E(G)$ (as $G$ is 5-connected) and $G[\{s',v,x_1,z_2\}]-x_1s'\cong
%K_4^-$ in which $x_1$ is of degree 2; so $(ii)$ holds. 

If  $L:=G[(F''-z_1)+A\cup \{s',t'\}]$ has independent paths $P_1,P_2$  from $t'$ to $s',z_1'$,
respectively, and if $Y_2$ has a cycle $C$ containing
$\{b, z_1,z_2\}$, then $G[\{b_1,t',v\}]\cup z_2v\cup (z_2s'\cup
P_1)\cup  C\cup (z_1z_1'\cup P_2)\cup z_1x_1v$ is a $TK_5$ in $G$
with branch vertices $b_1,t',v,z_1,z_2$. 
So we may assume $P_1,P_2$ do not exist, or $C$ does not exist. 

Suppose $P_1,P_2$ do not exist in $L$. Then $L$ has
1-separation $(L_1,L_2)$ such that $t'\in V(L_1-L_2)$ and
$\{s',z_1'\}\subseteq V(L_2)$. Since $G$ is 5-connected,
$|V(L_1)|=2$ and $x_1t'\in E(G)$. Now $G[\{b_1,t',v,x_1\}]-x_1b_1\cong
K_4^-$ in which $x_1$ is of degree 2, and $(ii)$ holds. 

Now assume $C$ does not exist. Then by Lemma~\ref{Watkins}, $Y_2$ has
2-cuts $S_b,S_z$ such that $b_1$ is a in component $D_b$ of $Y_2-S_b$,
$p_1=z_1$ is in a component $D_z$ of $Y_2-S_z$, and $V(D_b)\cap
(V(D_z)\cup S_z\cup \{p_2\})=\emptyset=V(D_z)\cap (V(D_b)\cup S_b\cup \{p_2\})$. If $y_2\notin
V(D_b)$ then $S_b\cup\{b,x_1\}$ is a cut in $G$, a contradiction. So
$y_2\in V(D_b)$. Then $y_2\notin V(D_z)$. Then $S_z\cup \{x_1,z_1'\}$ is
a cut in $G$, a contradiction.

\medskip

{\it Case} 3. $e(z_2, B_1-b_1)\le 1$.

If $y_2\in V(X)$ then, since $G$ is 5-connected, $e(z_1,B_1-b_1)\ge 2$
and $e(z_2,B_1-b_1)=1$. If 
$y_2\notin V(X)$ then, by $(iii)$ of Lemma~\ref{(I)},   $e(z_2,
B_1-b_1)=1$ and $e(z_1, B_1-b_1)\ge 2$. 

For convenience, let $F':=F_1'$, $F'':=F_1''$, $s:=s_1$ and
$t:=t_1$. Then $b_1,z_1\in V(F')$ and $y_1,z_2\in V(F'')-V(F')$. 
We choose $(F',F'')$ so that $F''$ is minimal.
Let $z_2'$ denote the unique neighbor of $z_2$ in $B_1-b_1$. 
We may assume $z_2'\ne y_1$; for, otherwise,
$G[\{x_1,x_2,y_1,z_1'\}]-x_1z_2\cong K_4^-$ in which $x_1$ is of
degree 2, and $(ii)$ holds. 
Note that if $z_2\ne p_2$ then $z_2b_1,z_2x_1\in E(G)$. By $(iii)$ of
Lemma~\ref{Cycle}, $G[Y_2+b_1+p_2Xz_2]$
contains a path $Q$ from $p_1$ to $b_1$ through $y_2,p_2$ in order. 

\medskip

{\it Subcase} 3.1.  $N(z_1Xp_1-z_1)\cap V(F''-\{s, t,z_2\})\not\subseteq \{z_2'\}$.

Let $uu'\in E(G)$ with $u'\in V(z_1Xp_1-z_1)$ and $u\in V(F'')-\{s, t, z_2,z_2'\}$. 
Since $B_1$ is 2-connected, $F'$ contains a path $S$ from $z_1$ to $b_1$ such that  $|V(S)\cap \{s, t\}|\leq 1$. 

Suppose $(F''-z_2)-S$ contains independent paths $S_1,S_2$ from $y_1$
to $z_2',u$, respectively. Then $G[\{x_1, x_2, y_1, y_2\}] \cup 
z_2x_2\cup z_2Qy_2\cup (z_2Qb_1\cup S\cup z_1x_1)\cup (z_2z_2'\cup S_1)\cup (S_2\cup uu'\cup u'Xp_1\cup p_1Qy_2)$ is a $TK_5$ in $G'$ with branch vertices $x_1,x_2,y_1,y_2,z_2$. 

So we may assume that such $S_1, S_2$ do not exist in $(F''-z_2)-S$
for any choice of $S$ and any choice of $u$. 
Hence, $(F''-z_2)-S$ has  a cut vertex $c$ which separates $y_1$ from $N(z_1Xp_1-z_1)\cup \{z_2'\}$. Denote by $M_1, M_2$ 
the $(\{c\}\cup (\{s, t\}\cap V(S)))$-bridges of $F''-z_2$ containing
$y_1$, $(N(z_1Xp_1-z_1)\cap V(F''-\{s,t,z_2\}))\cup \{z_2'\}$, respectively. 
Since $G$ is $5$-connected, $\{s, t\}\subseteq V(M_1)$ (to avoid the cut $\{c,x_1,x_2\}\cup (V(S)\cap \{s,t\})$) and $M_1$ contains 
independent paths $R_1, R_2, R_3$ from $y_1$ to $c, s, t$, respectively. Since $B_1$ is 2-connected, 
$\{s, t\}\cap V(M_2)\neq \emptyset$. Note that there exists $r\in
\{s,t\}\cap V(M_2)$ such that 
$M_2$ contains disjoint paths $T_1, T_2$ from $\{z_2', u\}$ to $\{c,
r\}$ and avoiding $\{s,t\}\cap V(M_2)-\{r\}$.
Now $R_1\cup R_3\cup T_1\cup T_2$ contains independent  paths from
$y_1$ to $z_2', u$, respectively,  and avoiding $\{s,t\}\cap V(M_2)-\{r\}$.
So by the nonexistence of $S_1,S_2$, $r\in V(S)$ for every choice of
$S$, which implies $b_1=r$. So we may assume $b_1=t$.

Choose $uu'$ so that $u'Xp_1$ is maximal. Since $\{s,t,u',x_1\}$
cannot be a cut in $G$
separating $F'$ from $F''\cup Y_2\cup p_2Xx_2$, there exists $ww'\in E(G)$ such
that $w\in V(F'-\{s, t,z_1\})$ and  $w'\in V(u'Xp_1-u')\cup V(p_2Xx_2)$. 

Suppose $w'\in V(u'Xp_1-u')$. Let $P_1$ be a path in $F'-\{z_1,t\}$ from $w$ to $s$ and 
$P_2$ be a path in $M_2-t$ from $z_2'$ to $u$. Then $G[\{x_1,
x_2, y_1, y_2\}] \cup z_2x_2\cup z_2Qy_2\cup (z_2Qb_1\cup R_3)
\cup (z_2z_2'\cup P_2\cup uu'\cup u'Xz_1\cup z_1x_1)\cup (R_2\cup P_1\cup ww'\cup w'Xp_1\cup p_1Qy_2)$ is a $TK_5$ in 
$G'$ with branch vertices $x_1,x_2,y_1,y_2,z_2$. 

Now assume $w'\in V(p_2Xx_2)$. Let $W$ be a path in  $F'-t$ from
$z_1$ to $w$. Then $X':=W\cup ww'\cup w'Xx_2$ is a path in $G-x_1$
from $z_1$ to $x_2$ such that in $(G-x_1)-X'$, $\{y_1,y_2\}$ is
contained in a cycle (which is contained in $(Y_2-p_2)\cup p_1Xu'\cup u'u\cup
M_2\cup (M_1-s)$). Hence by Lemma~\ref{4A} and Lemma~\ref{apex}, we may
assume that $X'$ is induced,  $y_1,y_2\notin V(X)$, and $(G-x_1)-X'$ is 2-connected. Thus, the
assertion of this lemma follows from Lemma~\ref{2-connected}. 
\medskip

{\it Subcase} 3.2. $N(z_1Xp_1-z_1)\cap V(F''-\{s, t,z_2\})\subseteq \{z_2'\}$. 

First, we show that $\{s,t,x_1, x_2, z_2'\}$ is a $5$-cut in $G$
separating $F''-z_2$ from $F'\cup Y_2\cup X$.  For, otherwise, there exists
$ww'\in E(G)$ with $w\in V(F'')-\{s,t,z_2'\}$ and $w'\in
V(p_2Xz_2-z_2)$. Let $P_1,P_2$ be independent paths in $F'$ from
$z_1$ to $r,b_1$, respectively, with $r\in \{s,t\}$. Without loss of
generality, we may assume $r=s$. By the minimality of $F''$,  $F''-t$ has independent paths $R_1,R_2$ from
$y_1$ to $s,w$, respectively. Now $G[\{x_1,x_2,y_1,y_2\}] \cup z_1x_1\cup (z_1Xp_1\cup Q_1)\cup
(P_1\cup R_1)\cup (P_2\cup b_1z_2x_2)\cup (R_2\cup ww'\cup w'Xp_2\cup
Q_2)$ is a $TK_5$ in $G$ with branch
vertices $x_1,x_2,y_1,y_2,z_1$.

Hence, since $G$ is 5-connected, $F''-z_2$ contains independent paths
$T_1, T_2, T_3$ from $y_1$ to $s, t, z_2'$, respectively.
%, and $F''-z_2$ has no $2$-cut separating $y_1$ from $\{s, t, z_2'\}$.
Let $F_g:=(F''-z_2)+\{g, gs, gt\}$, where $g$ is a new vertex;  
then  by Lemma~\ref{Watkins}, $F_g$ has a cycle containing $\{g,y_1,z_2'\}$. 
Thus, we may assume by symmetry that  $F''-z_2$ has a
path $S$ from $s$ to $t$ and  through $y_1, z_2'$ in order.

We may assume $N(x_2)\cap V(F'-\{s,t\})=\emptyset$. For, suppose there exists  $x_2^*\in N(x_2)\cap V(F'-\{s,t\})$. Since $B_1$ is 2-connected, 
$F'$ contains independent paths $R_1, R_2$ from $z_1$ to $x_2^*, r$, respectively, for some $r\in \{s, t\}$. (This can be done by considering whether
or not  $z_1$ and $x_2^*$ are contained in the same $\{s,t\}$-bridge of $F'$.)  Let $T=T_1$ if $r=s$, and $T=T_2$ if $r=t$. 
Then $G[\{x_1, x_2, y_1, y_2\}] \cup z_1x_1\cup (z_1Xp_1\cup Q_1)\cup (R_1\cup x_2^*x_2)\cup (R_2\cup T)\cup (Q_2\cup p_2Xz_2\cup z_2z_2'\cup T_3)$ is a $TK_5$ in $G'$ with branch vertices $x_1,x_2,y_1,y_2,z_1$. 

Let $x=p_1$ if $N(z_2')\cap V(z_1Xp_1-z_1)=\emptyset$, and otherwise
let $x\in N(z_2')\cap V(z_1Xp_1-z_1)$ with $z_1Xx$ minimal. 

 Suppose $z_2'x_2\in E(G)$. Then we may assume $x_1z_2\notin E(G)$; for
otherwise, $G[\{x_1,x_2,z_2,z_2'\}]-x_1z_2'\cong K_4^-$ in which
$x_1$ is of degree 2, and $(ii)$ holds. Hence, $z_2=p_2$, and $\{b_1, s, t,x,x_1\}$ is a $5$-cut in $G$ separating 
$F'\cup z_1Xx$ from $F''\cup Y_2$. Since $G$ is 5-connected, $b_1\notin \{s,t\}$. 
Let $(G_1,G_2)$ be a 5-separation in $G$ such that $V(G_1\cap G_2)=\{b_1,s,t,x,x_1\}$ and $G_2=
G[F'\cup z_1Xx+x_1]$. Clearly, $|V(G_2)|\ge 7$. We may assume
$|V(G_1)|\ge 7$; for, if not, $|V(G_1)|=6$, and
$G[\{b_1,s,t,z_1\}]-st\cong K_4^-$ and $(ii)$ holds. If
$(G_2-x_1,b_1,x,s,t)$ is planar then
the assertion of this lemma follows from Lemma~\ref{apex}. So we may
assume that this is not the case. Then by Lemma~\ref{2path}, 
$G_2-x_1$ has disjoint paths $S_1,S_2$ from $s,t$ to $b_1,x$,
respectively. Now $z_2z_2'x_2z_2\cup y_1x_2\cup y_1Sz_2'\cup (y_1Ss\cup S_1\cup
b_1Qz_2)\cup y_2Qz_2\cup (y_2Qp_1\cup p_1Xx\cup S_2\cup tSz_2')\cup y_2x_2\cup y_2x_1y_1$ is a $TK_5$ in $G'$ 
with branch vertices $x_2,y_1,y_2,z_2,z_2'$.

Now assume $z_2'x_2\notin E(G)$. Then $x_2$ has a neighbor in
$F''-\{y_1,z_2'\}$ (as $N(x_2)\cap V(F'-\{s,t\})\ne \emptyset$). Let $r$ be a new vertex. We may assume that $(F''+\{r,rs,rt\})-z_2$ 
has disjoint paths $S_1,S_2$ from $r,z_2'$ to $x_2,y_1$, respectively. 
For, suppose such paths $S_1,S_2$ do not exist. Then by Lemma~\ref{2path}, there exists a collection ${\cal A}$ of disjoint 
subsets of $F''-\{x_2,y_1,z_2\}$ such that $(F''+\{r,rs,rt\})-z_2,
r,y_1,x_2,z_2')$ is 3-planar. Since $G$ is 5-connected and $F''$ is minimal, we may assume 
$(F''-z_2, s,t, y_1,x_2,z_2')$ is planar. Thus, since $z_2'$ is the
only neighbor of $z_2$ in $F''-F'$, $G$ has a 5-separation
$(G_1',G_2')$ such that $V(G_1'\cap G_2')=\{s,t, x_1,x_2,z_2\}$,
$G_2'-x_1=F''$,  and $(G_2'-x_1, s,t, x_2,z_2)$ is planar. Since $|V(G_j')|\ge
7$ for $j\in [2]$, the assertion of this lemma follows from Lemma~\ref{apex}.

Without loss of generality, let $rs\in S_1$. 
If $F'-t$ has independent paths $P_1,P_2$ from $z_1$ to $s,b_1$,
respectively, then $G[\{x_1,x_2,y_2\}] \cup z_1x_1\cup (P_1\cup
(S_1-r))\cup (z_1Xp_1\cup p_1Qy_2)\cup (z_2z_2'\cup S_2\cup
y_1x_1)\cup z_2x_2\cup z_2Qy_2\cup (z_2Qb_1\cup P_2)$ is a $TK_5$ in $G'$ with branch vertices $x_1,x_2,y_2,z_1,z_2$. 
So we may assume that such $P_1,P_2$ do not exist in $F'-t$. 

Thus $F'$ has a 2-separation $(F_1,F_2)$ such that $t\in V(F_1\cap F_2)$, $z_1\in V(F_1-F_2)$ 
and $\{b_1,s\}\subseteq V(F_2-F_1)$.  Choose this separation so that $F_1$ is minimal. 
Let $s'\in V(F_1\cap F_2)-\{t\}$. Since $\{s',t,z_1,x_1\}$ cannot be a
cut in $G$, $V(F_1)=\{s',t,z_1\}$ or 
there exists $zz'\in E(G)$ such that 
$z\in V(z_1Xp_1-z_1)\cup V(p_2Xz_2-z_2)$ and $z'\in
V(F_1)-\{s',t,z_1\}$. 

First, assume $V(F_1)=\{s',t,z_1\}$. Then $z_1=p_1$ as $G$ is 5-connected. By
$(iii)$ of Lemma~\ref{Cycle}, let $Q'$ be a path in $Y_2$ from $p_2$ to $b_1$ and
through $y_2,p_1$ in order, and let $C$ be a cycle in $Y_2-b_1$
containing $\{p_1,p_2,y_2\}$.  Let $C':=Q'\cup
p_2Xz_2\cup z_2b_1$ if $z_2\ne p_2$; and  let $C':=C$ if $z_2=p_2$.
If $F'-\{b_1,t,z_1\}$ has a path $S$ from $s'$ to $s$ then
$x_1x_2y_2x_1\cup z_1x_1\cup z_2x_2\cup C'\cup (z_1s'\cup S\cup
S_1)\cup (z_2z_2'\cup S_2\cup y_1x_1)$ is a $TK_5$ in $G'$ with branch
vertices $x_1,x_2,y_2,z_1,z_2$. So we may assume such $S$ does not
exist. Then $F'$ has a separation $(L',L'')$ such that $V(L'\cap L'')
=\{b_1,t\}$, $\{s',z_1\}\subseteq V(L')$ and $s\in
V(L'')-\{b_1,t\}$. Since $G$ is 5-connected, $\{b_1,t,x_1,z_1\}$ is
not a cut in $G$, and $L'-\{b_1,t,z_1\}$ has a path $S'$ from $s'$
to some $z\in N(p_2Xz_2-z_2)$. Let $z'\in N(z)\cap
V(p_2Xz_2-z_2)$. Let $S$ be a path in $L''-t$ from $s$ to $b_1$. Then
$G[\{x_1,x_2,y_1,y_2\}]\cup z_1x_1\cup Q_1\cup (z_1s'\cup S'\cup
zz'\cup z'Xx_2)\cup (z_1t\cup T_2)\cup (T_1\cup S\cup Q_3)$ is a
$TK_5$ in $G'$ with branch vertices $x_1,x_2,y_1,y_2,z_1$. 

Thus, we may assume that $zz'\in E(G)$ such that 
$z\in V(z_1Xp_1-z_1)\cup V(p_2Xz_2-z_2)$ and $z'\in
V(F_1)-\{z_1,s',t\}$. 

Suppose $z\in V(xXp_1-x)$. Let $X^*=z_1Xx\cup xz_2'z_2x_2$. Then,
$T_1\cup T_2\cup (F'-z_1)\cup zz'\cup zXp_1\cup Y_2$ is contained in 
$G-X^*$ and has a cycle containing $\{y_1,y_2\}$. Hence, by
Lemma~\ref{4A} and then Lemma~\ref{apex}, we may assume that $G-x_1$
has an induced path $X'$ from $z_1$ to $x_2$ such that $y_1,y_2\notin
V(X')$ and $G-X'$ is 2-connected. Then the assertion of this lemma
follows from Lemma~\ref{2-connected}.

Now suppose $z\in V(p_2Xz_2-z_2)$. By the minimality of $F_1$, $F_1-t$ has independent
paths $L_1,L_2$ from $z_1$ to $s',z'$, respectively. In $F_2\cup (F''-z_2)$, we find independent paths $L_1',L_2'$ from $y_1$ to $s',b_1$, respectively. 
Then $G[\{x_1,x_2,y_1, y_2\}] \cup z_1x_1\cup (z_1Xp_1\cup Q_1)\cup (L_1\cup L_1')\cup (L_2\cup z'z\cup zXx_2)\cup (L_2'\cup b_1b\cup Q_3)$ is a $TK_5$ in $G'$ with branch vertices $x_1,x_2,y_1, y_2,z_1$. 

Hence,  we may assume $z\in V(z_1Xx-z_1)$ for all such $zz'$. 
Choose $z$ with $z_1Xz$ is maximal. Since   $\{s',t,x_1,z\}$ cannot be a cut in $G$, there exists $uu'\in E(G)$ such that 
$u\in V(z_1Xz)-\{z_1,z\}$ and $u'\in V(F_2)-\{s',t\}$. Let $P_1$ be a path in $F_1-\{s',z_1\}$ from $z'$ to $t$, and $P_2$ be a path in $F_2-t$ from $u'$ to $b_1$. 
Then $G[\{x_1,x_2,y_1,y_2\}] \cup z_2x_2\cup (z_2z_2'\cup T_3)\cup (z_2Xp_2\cup p_2Qy_2)\cup (z_2Qb_1\cup
P_2\cup u'u\cup uXz_1\cup z_1x_1)\cup (T_2\cup P_1\cup z'z\cup zXp_1\cup
p_1Qy_2)$ is a $TK_5$ in $G'$ with branch vertices $x_1,x_2,y_1,y_2,z_2$.  \qed

\section{Finding $TK_5$}

Recall the notation from Lemma~\ref{classify2} and the previous
section. In particular, $H:=G[B_1+\{z_1, z_2\}]$, $G':=G-\{x_1x:x\notin
\{x_2,y_1,y_2,z_0,z_1\}\}$, $b_1\in N(y_2)\cap V(B_1)$ and $p_1=p_2=b=y_2$ if $y_2\in
V(X)$, and $b_1\in V(B_1\cap B_2)$ and $V(Y_1\cap Y_2)=\{b,p_1,p_2\}$ if $y_2\notin V(X)$.  
 Our objective  is to find $TK_5$ in $G'$ using the structural information on $H$ produced in the previous sections. By Lemma~\ref{Cycle}, 
\begin{itemize}
\item [(A1)] $Y_2$ has independent  paths $Q_1, Q_2, Q_3$  from $y_2$ to $p_1, p_2, b$, respectively.
\end{itemize} 

Note that if $y_2\in V(X)$ then  $e(z_1,B_1-b_1)\ge 2$ and
$e(z_2,B_1-b_1)\ge 1$. Thus,  by Lemma~\ref{(I)},  we may assume that there exists $i\in
[2]$ such that $e(z_i, B_1-b_1)\ge 2$ and $e(z_{3-i}, B_1-b_1)\ge 1$. 
(Hence, by Lemma~\ref{(I)'}, $e(z_{3-i}, B_1) = 1$ only if $b= b_1$ and,
therefore, $z_{3-i}=p_{3-i}$.) Then by Lemma~\ref{Cycle},
\begin{itemize}
\item [(A2)] $Y_2$ has a path $T$ from $b$ to $p_i$ through $p_{3-i}, y_2$ in order, respectively.
\end{itemize} 

By Lemma~\ref{YZ}, we may assume that 
\begin{itemize}
\item [(A3)] $H$ has disjoint paths $Y, Z$ from $y_1, z_1$ to $b_1, z_2$, respectively.
\end{itemize}

By Lemma~\ref{ABC}, we may assume that 

\begin{itemize}
\item [(A4)] $H$ has independent paths $A, B, C$,  with $A, C$ from $z_i$ to $y_1$, and $B$ from $b_1$ to $z_{3-i}$. 
\end{itemize}

Let $J(A,C)$ denote the $(A\cup C)$-bridge of $H$ containing $B$, and $L(A, C)$ denote the union of all $(A\cup C)$-bridges of $H$ with attachments on both $A$ and $C$. We may choose $A, B, C$ such that the following are satisfied in the order listed:
\begin{itemize}
\item [(a)] $A, B,C$ are induced paths in $H$,
\item [(b)] whenever possible, $J(A, C)\subseteq L(A,C)$, 
\item [(c)] $J(A, C)$ is maximal, and
\item [(d)] $L(A, C)$ is maximal.
\end{itemize}
 We refer the reader to Figure~\ref{structure} for an illustration. 
We may assume that 

\begin{itemize}
\item [(A5)] for any $j\in [2]$, $H$ contains no path from $z_j$ to
  $b_1$ and through $z_{3-j},y_1$ in order.
% in particular, $y_1z_1, y_1z_2\notin E(G)$.
\end{itemize}
For, suppose $H$ does contain a path $R$ from $z_j$ to $b_1$ and
through $z_{3-j},y_1$ in order. Then $G[\{x_1, x_2, y_1, y_2\}]
\cup z_{3-j}x_{3-j}\cup (z_{3-j}Xp_{3-j}\cup Q_{3-j})\cup
(z_{3-j}Rz_j\cup z_jx_j)\cup z_{3-j}Ry_1\cup (y_1Rb_1\cup b_1b\cup
Q_3)$ is a $TK_5$ in $G'$ with branch vertices
$x_1,x_2,y_1,y_2,z_{3-j}$. Thus, we may assume (A5). 
\medskip

Since $B_1$ is 2-connected and $e(z_{3-i},B_1-b_1)\ge 1$, $H$ has disjoint paths $P, Q$ from $p,q\in V(B)$ to $c, a\in V(A\cup C)-\{z_i\}$, 
respectively, and internally disjoint from $A\cup B\cup C$.
By symmetry between $A$ and $C$, we may assume that $b_1, p, q, z_{3-i}$ 
occur on $B$ in order. By (A5), $c\ne y_1$. We choose such $P, Q$ that 
the following are satisfied in order listed:

\begin{itemize}
\item [(A6)] $qBz_{3-i}$ is minimal, $pBz_{3-i}$ is maximal, the subpath of $(A\cup C)-z_i$ between $a$ and $y_1$ is minimal, 
and the subpath of $(A\cup C)-z_i$ between $c$ and $y_1$ is maximal.
\end{itemize}

Let $B'$ denote the union of $B$ and the $B$-bridges of $H$ not
containing $A\cup C$. Note that all paths in $H$ from $A\cup C$ to
$B'$ and internally disjoint from $B'$ must have an end in $B$.
We may assume that

\begin{itemize}
\item [(A7)] if $e(z_{3-i},B_1)\ge 2$ then, for any $q^*\in V(B'-q)$, $B'$ has  independent paths from $z_{3-i}$ to $q,q^*$, respectively. 
\end{itemize}
For, suppose $e(z_{3-i},B_1)\ge 2$ and for some $q^*\in V(B'-q)$, $B'$
has  no independent paths from $z_{3-i}$ to $q,q^*$, respectively. Then $q\ne z_{3-i}$, and $B'$ 
has a 1-separation $(B_1',B_2')$ such that $q,q^*\in V(B_2')$ and
$z_{3-i}\in V(B_1')-V(B_2')$. Note that $b_1\in V(B_2')$. 
Choose $(B_1',B_2')$ with $B_1'$ minimal, and let $z\in V(B_1'\cap
B_2')$. 

Since $e(z_{3-i},B_1)\ge 2$, $|V(B_1')|\ge 3$; so 
$H$ has a path $R$ from some $s\in V(B_1'-z)$ to some $t\in V(A\cup C\cup P\cup Q)$ and internally disjoint from $A\cup B\cup C\cup P\cup Q$.
By the choice of $P,Q$ in (A6), we see that $t=z_i$. Let $S$ be a path
in $B_1'$ from $z_{3-i}$ to $s$, respectively. Let $R=A\cup y_1Cc$ if
$c\in V(C)$, and $R=C\cup y_1Ac$ if $c\in V(A)$. Then 
 Now $S\cup R\cup  P\cup pBb_1$ is a path contradicting (A5). 

\medskip 
We will show that we may assume  $a=y_1$
(see (3)), derive structural information about $G'$ and $H$ (see
(4)--(7)), and consider whether
or not $z_i\in V(J(A,C))$ (see Case 1 and Case 2). 
First, we may assume that 

\begin{itemize}
\item [(1)] $N(y_1)\cap V(z_jXp_j-z_j)=\emptyset$ for $j\in [2]$. 
\end{itemize}
For, suppose there exists $s\in N(y_1)\cap V(z_jXp_j-z_j)$ for some
$j\in [2]$. By symmetry, assume $c\in V(C)$. 
If $j = 3-i$ then, using  $Q_1,Q_2,Q_3$ from (A1), we see
that $G[\{x_1, x_2, y_1, y_2\}] \cup z_ix_i\cup (z_iXp_i\cup Q_i)\cup A\cup (z_iCc\cup P\cup pBz_{3-i}\cup z_{3-i}x_{3-i})\cup (y_1s\cup sXp_{3-i}\cup Q_{3-i})$ is a $TK_5$ in $G'$ with branch vertices $x_1,x_2,y_1,y_2,z_i$. 

So assume $j= i$. Suppose $e(z_{3-i},B_1)=1$. Then
$z_{3-i}=p_{3-i}$. Recall the path $T$ from (A2). Note that
$z_{3-i}Tb\cup bb_1\cup A\cup B\cup C\cup P\cup Q$ 
contains independent paths $S_1, S_2$ from $z_{3-i}$ to $z_i, y_1$, respectively. 
Hence $G[\{x_1, x_2, y_1, y_2\}] \cup z_{3-i}x_{3-i}\cup z_{3-i}Ty_2\cup (S_1\cup z_ix_i)\cup S_2\cup (y_1s\cup sXp_i\cup p_iTy_2)$ 
is a $TK_5$ in $G'$ with branch vertices $x_1,x_2,y_1,y_2,z_{3-i}$. 

Now assume $e(z_{3-i},B_1)\ge 2$. Let $P_1,P_2$ be independent paths from (A7) with $q^*=p$. 
Then $P_1\cup P_2\cup A\cup B\cup C\cup P\cup Q$ contains 
independent  paths $S_1, S_2$ from $z_{3-i}$ to $z_i, y_1$,
respectively. Now $G[\{x_1, x_2, y_1, y_2\}]\cup z_{3-i}x_{3-i}\cup (z_{3-i}Xp_{3-i}\cup Q_{3-i}) 
\cup (S_1\cup z_ix_i)\cup S_2\cup (y_1s\cup sXp_i\cup Q_i)$ 
is a $TK_5$ in $G'$ with branch vertices
$x_1,x_2,y_1,y_2,z_{3-i}$. This proves (1).

\medskip
We may also assume 

\begin{itemize} 
\item [(2)] $y_1\in V(J(A, C))$. 
\end{itemize}
For, suppose $y_1\notin V(J(A, C))$. By (1) and 5-connectedness of $G$, 
$y_1\in V(D_1)$ for some $(A\cup C)$-bridge $D_1$ of $H$ with $D_1\ne J(A, C)$. 
Thus, let $D_1, \ldots, D_k$ be a maximal sequence of $(A\cup
C)$-bridges of $H$ with $D_j\neq J(A,C)$ for $j\in [k]$, such that, 
for each $l \in [k -1]$, 
$$D_{l+1} \mbox{  has a vertex not in } \bigcup_{j\in [l]}(c_jCy_1\cup
a_jAy_1) \mbox{  and a vertex not in } \bigcap_{j\in [l]}(z_iCc_j \cup z_iAa_j),  $$
 where for each $j\in [k]$, $a_j\in V(D_j\cap A)$ and $c_j\in V(D_j\cap C)$ such that $a_jAy_1$ and $c_jCy_1$ are maximal.
Let $S_l:=\bigcup_{j\in [l]} (D_j\cup a_jAy_1\cup c_jCy_1)$.

We claim that for any $ l\in [k]$ and for any $r_l\in V(S_l)-\{a_l,c_l\}$, $S_l$ has three independent paths $A_l,C_l,R_l$  
from $y_1$ to $a_l,c_l,r_l$, respectively. This is obvious for $l=1$ (if $a_l=y_1$, or $c_l=y_1$, or
$r_l=y_1$  then $A_l$, or $C_l$, or $R_l$  is a trivial path). 
Now assume $k\ge 2$ and the claim  holds for some $l\in [k-1]$. Let $r_{l+1}\in
V(S_{l+1})-\{a_{l+1},c_{l+1}\}$. When $r_{l+1}\in V(S_l)-\{a_l,c_l\}$ let
$r_l:=r_{l+1}$; otherwise, let $r_l\in V(a_lAy_1-a_l)\cup
V(c_lCy_1-c_l)$ with $r_l\in V(D_{l+1})$. By assumption, $S_l$ has independent paths $A_l,C_l,R_l$ from $y_1$ to
$a_l,c_l,r_l$, respectively. 
If $r_{l+1}\in V(S_l)-\{a_l,c_l\}$ then $A_{l+1}:=A_l\cup
a_lAa_{l+1}, C_{l+1}:=C_l\cup c_lCc_{l+1}, R_{l+1}:=R_l$ are the
desired paths in $S_{l+1}$. 
If $r_{l+1}\in V(D_{l+1})-V(A\cup C)$ then
let $P_{l+1}$ be a path in $D_{l+1}$ from $r_l$ to $r_{l+1}$
internally disjoint from $A\cup C$; we see that $A_{l+1}:=A_l\cup
a_lAa_{l+1}, C_{l+1}:=C_l\cup c_lCc_{l+1}, R_{l+1}:=R_l\cup P_{l+1}$
are the desired paths in $S_{l+1}$. 
So we may assume by symmetry
that $r_{l+1}\in V(a_{l+1}Aa_l-a_{l+1})$. Let $Q_{l+1}$ be a path in
$D_{l+1}$ from $r_l$ to $a_{l+1}$ internally disjoint from $A\cup
C$. Now $R_{l+1}:=A_l\cup a_lAr_{l+1}, C_{l+1}:=C_l\cup c_lCc_{l+1},
A_{l+1}:=R_l\cup Q_{l+1}$ are the desired paths in $S_{l+1}$. 

Hence, by (c), $J(A, C)$ does not intersect $(a_kAy_1\cup c_kCy_1) -
\{a_k, c_k\}$. In particular, $a,c\notin (a_kAy_1\cup c_kCy_1) -\{a_k, c_k\}$.
Since $G$ is 5-connected,  $\{a_k, c_k, x_1,x_2\}$ cannot be a cut in $G$ separating $S_k$ from $X\cup J(A,C)$. 
So there exists $ss'\in E(G)$ such that $s\in V(S_k)-\{a_k, c_k\}$ and $s'\in V(z_1Xp_1\cup z_2Xp_2)$. 
By the above claim,  let $A_k, C_k, R_k$ be independent paths in $S_k$ from $y_1$ to $a_k, c_k, s$, respectively; so 
$s'\notin \{z_1,z_2\}$ by (c). 

Suppose $s'\in V(z_{3-i}Xp_{3-i}-z_{3-i})$. Then 
$G[\{x_1,x_2,y_1,y_2\}] \cup z_ix_i\cup (z_iXp_i\cup Q_i)\cup (z_iCc\cup P\cup pBz_{3-i}\cup z_{3-i}x_{3-i})\cup 
(z_iAa_k\cup A_k)\cup (R_k\cup ss'\cup s'Xp_{3-i}\cup Q_{3-i})$ is a $TK_5$ in $G'$ with branch vertices $x_1,x_2,y_1,y_2,z_i$. 

So we may assume $s'\in V(z_iXp_i-z_i)$.
Suppose  $e(z_{3-i}, B_1)=1$. Then  $z_{3-i}=p_{3-i}$. Recall the 
path $T$ from (A2). Note that 
$z_{3-i}Tb\cup bb_1\cup  z_iAa_k\cup z_iCc_k\cup P\cup Q\cup B$ contains independent paths $S_1, S_2$ from $z_{3-i}$ to $z_i, v$, respectively, for 
some $v\in \{a_k, c_k\}$. Let $S=A_k$ if $v=a_k$, and $S=C_k$ if
$v=c_k$. Then $G[\{x_1, x_2, y_1, y_2\}] \cup z_{3-i}x_{3-i}\cup z_{3-i}Ty_2\cup (S_1\cup z_ix_i)\cup (S_2\cup S)\cup (R_k\cup ss'\cup s'Xp_i\cup p_iTy_2)$ is a $TK_5$ in $G'$ with branch vertices $x_1,x_2,y_1,y_2, z_{3-i}$. 

Hence, we may assume $e(z_{3-i},B_1)\ge 2$. Let $P_1,P_2$ be independent paths from (A7) with $q^*=p$. 
Then, $P_1\cup P_2\cup  z_iAa_k\cup z_iCc_k\cup P\cup Q\cup B$ contains independent paths $S_1, S_2$ from $z_{3-i}$ to $z_i, v$, respectively, for 
some $v\in \{a_k, c_k\}$. Let $S=A_k$ if $v=a_k$, and $S=C_k$ if
$v=c_k$. Then $G[\{x_1, x_2, y_1, y_2\}] \cup z_{3-i}x_{3-i}\cup (z_{3-i}Xp_{3-i}\cup Q_{3-i}) \cup (S_1\cup z_ix_i)\cup (S_2\cup S)\cup (R_k\cup ss'\cup s'Xp_i\cup Q_i)$ is a $TK_5$ in $G'$ with branch vertices $x_1,x_2,y_1,y_2, z_{3-i}$. 
This completes the proof of  (2). 

\medskip

For convenience, we let  $K:=A\cup B\cup C\cup P\cup Q$. We claim that
\begin{itemize}
\item [(3)] $a=y_1$
\end{itemize}
Suppose $a\ne y_1$.  
By (2), $J(A, C)$ has a path $S$ from $y_1$ to some vertex $s\in
V(P\cup Q\cup B)-\{c, a\}$ and internally disjoint from $K$.
By (A6),  $s\notin V(Q\cup qBz_{3-i})$. So $s\in V(P\cup b_1Bq) - \{a,q\}$.
If $a\in V(A)$ let $R=aAz_i$ and $R'=C$; and  if
$a\in V(C)$ let $R=aCz_i$ and $R'=A$. Also, let $S'=S\cup sBb_1$ if $s\in V(B)$, and $S'=S\cup
sPp\cup pBb_1$ if $s\in V(P)$. Then $z_{3-i}Bq\cup Q\cup R\cup R'\cup S'$
is a path contradicting (A5). 

\medskip

By symmetry between $A$ and $C$, we may assume $c\in V(C)$. Before we distinguish cases according to whether or not $z_i\in
V(J(A,C))$, we derive further information about $G'$. 
We may assume that 
\begin{itemize}
\item [(4)] for any path $W$  in $G'$ from $x_i$ to some $w\in V(K)-\{z_i, y_1\}$ and internally disjoint from $K$,
 we have $w\in V(A)-\{z_i, y_1\}$.
\end{itemize}
To see this, 
suppose $w\notin V(A)-\{z_i, y_1\}$.  
First, assume  $e(z_{3-i},B_1) = 1$. Then $z_{3-i}=p_{3-i}$. Recall
the path $T$ from (A2), and note that    
$z_{3-i}Tb_1\cup B\cup (C-z_i)\cup P\cup Q\cup W$ contains independent paths $S_1, S_2$ from $z_{3-i}$ to $x_i, y_1$, respectively.  
Then $G[\{x_1, x_2, y_1, y_2\}] \cup z_{3-i}x_{3-i}\cup z_{3-i}Ty_2\cup S_1\cup S_2\cup  (A\cup z_iXp_i\cup p_iTy_2)$ is a $TK_5$ in $G'$ with branch vertices $x_1, x_2, y_1, y_2, z_{3-i}$. 

Thus, we may assume  $e(z_{3-i},B_1)\ge 2$. Let $P_1,P_2$ be
independent paths in $B'$ from (A7) with $q^*=p$. So $P_1\cup P_2\cup
B\cup (C-z_i)\cup W \cup P\cup Q$ contains independent paths $S_1,
S_2$ from $z_{3-i}$ to $x_i, y_1$, respectively. Then $G[\{x_1,
x_2, y_1, y_2\}] \cup z_{3-i}x_{3-i}\cup (z_{3-i}Xp_{3-i}\cup Q_{3-i})
\cup S_1\cup S_2\cup  (A\cup z_iXp_i\cup Q_i)$ is a $TK_5$ in $G'$
with branch vertices $x_1,x_2,y_1,y_2,z_{3-i}$. This completes the proof of (4).

\medskip

Since $G$ is 5-connected and $z_0\in V(B_1)$ when $e(z_1, B_1)\ge 2$
(see $(iv)$ of Lemma~\ref{classify2}), 
it follows from  (4) that 
$$\mbox{$G'$ has a path $W$ from $x_i$ to $w\in
V(A)-\{y_1,z_i\}$ and internally disjoint from $K$.}$$
 Hence, $|V(A)|\ge 3$. Also,  $|V(C)|\ge 3$ as $c\in V(C)-\{y_1,z_1\}$. Since $A$ and $C$ are induced paths in $H$,
$$y_1z_i\notin E(G).$$ 
We may assume that 
\begin{itemize}
\item [(5)] $G'$ has no path from $z_{3-i}Xp_{3-i}-y_2$ to $(A\cup C)-y_1$ and internally disjoint from 
$K$,  $G'$ has no path from $z_iXp_i-z_i$ to $(A\cup cCy_1)-\{z_i, c\}$ and internally disjoint from 
$K$, and if $i=1$ then $G'$ has no path from $x_2$ to $(A\cup C)-y_1$ and internally disjoint from $K$. 
\end{itemize}
First, suppose $S$ is a path in $G'$ from some $s\in V(z_{3-i}Xp_{3-i}-y_2)$ to some $s'\in V(A\cup C)-\{y_1\}$. Then $A\cup C\cup S$ contains 
independent paths $S_1,S_2$ from $z_i$ to $y_1,s$, respectively.
Hence,  $G[\{x_1,x_2,y_1,y_2\}] \cup z_ix_i\cup (z_iXp_i\cup Q_i)\cup S_1\cup (S_2\cup sXz_{3-i}\cup z_{3-i}x_{3-i})\cup (Q\cup qBb_1\cup b_1b
\cup Q_3)$ is a $TK_5$ in $G'$ with branch vertices $x_1,x_2,y_1,y_2,z_i$.

Now assume that $S$ is a path in $G'$ from some  $s\in V(z_iXp_i-z_i)$ to some $s'\in V(A\cup cCy_1)-\{z_i, c\}$ 
and internally disjoint from $K$. Let $S'=y_1As'$ if $s'\in V(A)$, and
$S'=y_1Cs'$ if $s'\in V(cCy_1)$. If $e(z_{3-i},B_1)= 1$  then
$z_{3-i}=p_{3-i}$ and, using the path  $T$ from  (A2), we see that
 $G[\{x_1, x_2, y_1, y_2\}] \cup z_{3-i}x_{3-i}\cup z_{3-i}Ty_2\cup (z_{3-i}Bq\cup Q)\cup (z_{3-i}Tb_1\cup b_1Bp
\cup P\cup cCz_i\cup z_ix_i)\cup (S'\cup S\cup sXp_i\cup p_iTy_2)$ is a $TK_5$ in $G'$ with branch vertices $x_1,x_2,y_1,y_2,z_{3-i}$. So assume  $e(z_{3-i},B_1)\ge 2$. Let $P_1,P_2$ be independent paths from (A7) with  $q^*=p$. 
Now $G[\{x_1, x_2, y_1, y_2\}] \cup z_{3-i}x_{3-i}\cup (z_{3-i}Xp_{3-i}\cup Q_{3-i})\cup (P_1\cup Q)\cup (P_2
\cup P\cup cCz_i\cup z_ix_i)\cup (S'\cup S\cup sXp_i\cup Q_i)$ is a $TK_5$ in $G'$ with branch vertices $x_1,x_2,y_1,y_2,z_{3-i}$.

Now suppose $i=1$ and $S$ is a path in $G'$ from $x_2$ to some  $s\in V(A\cup C)-\{y_1\}$ and internally disjoint from $K$. 
If $s\in V(A-y_1)$, then $G[\{x_1, x_2, y_1, y_2\}] \cup z_1x_1\cup (z_1Xp_1\cup Q_1)\cup C\cup (z_1As\cup S)\cup (Q\cup qBb_1\cup b_1b\cup Q_3)$ is a $TK_5$ in $G'$ with branch vertices $x_1,x_2,y_1,y_2,z_1$.  
So assume $s\in V(C-y_1)$. Then $G[\{x_1, x_2, y_1, y_2\}] \cup
z_1x_1\cup (z_1Xp_1\cup Q_1)\cup A\cup (z_1Cs\cup S)\cup (Q\cup
qBb_1\cup b_1b \cup Q_3)$ is a $TK_5$ in $G'$ 
with branch vertices $x_1,x_2,y_1,y_2,z_1$. This completes the proof of (5).

\begin{itemize}
\item [(6)] We may assume that 
\begin{itemize}
\item [(6.1)] any path in $J(A, C)$ from $A-\{z_i, y_1\}$ to $(P\cup Q\cup B)-\{c,y_1\}$ and internally disjoint from $K$ must end on $Q$, 
\item [(6.2)] if an $(A\cup C)$-bridge of $H$ contained in $L(A, C)$ intersects $z_iCc-c$ and contains a vertex $z\in V(A-z_i)$ 
then $J(A, C)\cap (z_iAz-\{z_i, z\})=\emptyset$, and 
\item [(6.3)] $J(A, C)\cap (z_iCc-\{z_i, c\})=\emptyset$, and any path in $J(A, C)$ from $z_i$ to 
$(P\cup Q\cup B)-\{c, y_1\}$ and internally disjoint from $K$ must end on $(P-c)\cup b_1Bp$.
\end{itemize}
\end{itemize}
To prove (6.1), let $S$ be a path in $J(A, C)$ from $s\in V(A)-\{z_i, y_1\}$ to $s'\in V(P\cup B)-\{c, q, y_1\}$ and internally disjoint from $K$. 
Note that  $s'\notin V(qBz_{3-i}-q)$ by (A6). 
Suppose  $e(z_{3-i},B_1) = 1$. Then $z_{3-i}=p_{3-i}$ and we use the path $T$ from (A2). Let $S'$ be a path in $(P-c)\cup (b_1Bq-q)$ from $b_1$ to $s'$. 
Then $G[\{x_1, x_2, y_1, y_2\}] \cup z_{3-i}x_{3-i}\cup z_{3-i}Ty_2\cup (z_{3-i}Tb_1\cup S'\cup S\cup 
sAw\cup W)\cup (z_{3-i}Bq\cup Q)\cup (C\cup z_iXp_i\cup p_iTy_2)$ is a $TK_5$ in $G'$  with branch vertices 
$x_1,x_2,y_1,y_2,z_{3-i}$. So we may assume $e(z_{3-i},B_1)\ge 2$. Let
$P_1,P_2$ be the  paths from (A7), with $q^*=p$ when $s'\in V(P)$ and
$q^* = s'$ when $s'\in V(B)$. So $P_1\cup P_2\cup B\cup S\cup Q$
contains independent paths $S_1, S_2$ from $z_{3-i}$ to $s, y_1$,
respectively. Now $G[\{x_1, x_2, y_1, y_2\}] \cup z_{3-i}x_{3-i}
\cup (z_{3-i}Xp_{3-i}\cup Q_{3-i})\cup (S_1\cup sAw\cup W)\cup S_2\cup (C\cup z_iXp_i\cup Q_i)$ is a $TK_5$ in $G'$  with branch vertices 
$x_1,x_2,y_1,y_2,z_{3-i}$.

\medskip

To prove (6.2), let $D$ be a path contained in $L(A, C)$ from $z'\in V(z_iCc-c)$ to $z\in V(A-z_i)$ and internally disjoint from $K$. 
Suppose there exists $s\in V(J(A, C))\cap V(z_iAz-\{z_i, z\})$. By (6.1), $J(A,C)$ has a path $S$ from $s$ to some $s'\in V(Q-y_1)$ 
and internally disjoint from $K$. Then $G[\{x_1, x_2, y_1, y_2\}]
\cup z_ix_i\cup (z_iXp_i\cup Q_i)\cup (z_iAs\cup S\cup s'Qq\cup qBz_{3-i}\cup z_{3-i}x_{3-i})
\cup (z_iCz'\cup D\cup zAy_1)\cup (y_1Cc\cup P\cup pBb_1\cup b_1b\cup Q_3)$ is a $TK_5$ in $G'$ 
with branch vertices $x_1,x_2,y_1,y_2,z_i$.

\medskip
To prove (6.3), let $S$ be a path in $J(A, C)$ from $s\in V(z_iCc-c)$ to $s'\in V(P\cup Q\cup B)-\{c, y_1\}$ and internally disjoint from $K$. 
Suppose  $s'\in V(Q\cup z_{3-i}Bp)-\{p,y_1\}$. Then $(S\cup Q\cup pBz_{3-i})-\{p,y_1\}$ contains a path $S'$  from $s$ to $z_{3-i}$. 
%and internally disjoint from $P\cup pBb_1$. 
So $G[\{x_1, x_2,y_1, y_2\}] \cup z_ix_i\cup (z_iXp_i\cup Q_i)\cup (z_iCs\cup S'\cup z_{3-i}x_{3-i})\cup A\cup (y_1Cc\cup P
\cup pBb_1\cup b_1b\cup Q_3)$ is a $TK_5$ in $G'$ with branch vertices $x_1,x_2,y_1,y_2,z_i$. 
Thus, we may assume  $s'\in V(P-c)\cup V(b_1Bp)$.  By (A6), $s=z_i$. This proves (6).

\medskip 

Denote by $L(A)$ (respectively,  $L(C)$) the union of all $(A\cup C)$-bridges of $H$ whose intersection with $A\cup C$ is contained in 
$A$ (respectively, $C$).

\begin{itemize}
 \item [(7)]  $L(A)=\emptyset$, and $L(C)\cap C \subseteq z_iCc$. 
\end{itemize}
 Suppose $L(A)\ne \emptyset$, and let $R_1$ be an $(A\cup C)$-bridge
 of $H$ contained in $L(A)$. Let $R_1, \ldots, R_m$ be  a maximal sequence of $(A\cup C)$-bridges of $H$ contained in $L(A)$, such that for $2\leq i \leq m$, $R_i$ has a vertex  internal to $\bigcup_{j=1}^{i-1}l_jAr_j$ (which is a path), 
where $l_j, r_j\in V(R_j\cap A)$ with $l_jAr_j$ maximal. 
Let  $a_1, a_2\in V(A)$ such that
$\bigcup_{j=1}^{m}l_jAr_j=a_1Aa_2$. By (c), $J(A, C)\cap
(a_1Aa_2-\{a_1, a_2\})=\emptyset$;
by (d) and the maximality of $R_1,\ldots,R_m$, $L(A, C)$ has no path from $a_1Aa_2-\{a_1, a_2\}$ to
$(A-a_1Aa_2)\cup (C-\{y_1,z_i\})$; and by (5), 
$(z_1Xp_1\cup z_2Xp_2)-\{a_1,a_2,z_i\}$ contains no neighbor of 
$(\bigcup_{j=1}^mR_j\cup a_1Aa_2)-\{a_1, a_2\}$. Hence,  $\{a_1, a_2, x_1, x_2\}$ is a cut in $G$, a contradiction. Therefore, $L(A)=\emptyset$.

Now assume $L(C)\cap C\not \subseteq z_iCc$, and let $R_1$ be an
$(A\cup C)$-bridge of $H$ contained in $L(C)$ such that  $R_1\cap
(cCy_1-c)\neq \emptyset$.  Let $R_1, \ldots, R_m$ be a maximal
sequence of $(A\cup C)$-bridges of $H$ contained in $L(C)$ such that for $2\leq i \leq m$, 
$R_i$ has a vertex  internal to $\bigcup_{j=1}^{i-1}l_jCr_j$ (which is a path),   where $l_j, r_j\in V(R_j\cap C)$ with $l_jCr_j$ maximal. Let 
$c_1, c_2\in V(C)$ such that $\bigcup_{j=1}^{m}l_jCr_j=c_1Cc_2$. By
the existence of $P$ and (c), $c_1, c_2\in V(cCy_1)$; 
by (c), $J(A, C)\cap (c_1Cc_2-\{c_1, c_2\})=\emptyset$; 
by (d), $L(A, C)\cap
(c_1Cc_2-\{c_1, c_2\})=\emptyset$; and by (5) and the maximality of $R_1,\ldots, R_m$,
$z_1Xp_1\cup z_2Xp_2$ contains no neighbor of $(\bigcup_{j=1}^mR_j\cup
c_1Cc_2)-\{c_1, c_2\}$. 
Hence, $\{c_1, c_2, x_1, x_2\}$ is a cut in $G$, a contradiction.
Therefore, $L(C)\cap C\subseteq z_iCc$. This proves (7).

\medskip

Let $F$ be the union of all $(A\cup C)$-bridges of $H$ different from
$J(A,C)$ and intersecting $z_iCc-c$. When $F\ne \emptyset$, let $a^*\in V(F\cap A)$ with $a^*Ay_1$ minimal, and let $r$ be the neighbor of $(F \cup z_iAa^*\cup z_iCc)-\{ a^*,c\}$ on $z_iXp_i$ with $rXp_i$ minimal.

\bigskip

{\it Case} 1. $z_i\in V(J(A, C))$.

By (6.3), $J(A, C)$ contains a path $S$ from $z_i$ to some $s\in V(P-c)\cup V(b_1Bp)$ and internally disjoint from $K$.

\medskip

{\it Subcase} 1.1. $F\neq \emptyset$.

Suppose $r\ne z_i$.  Then by (5) and the definition of $r$, $G'$ has a path $R$ 
from $r$ to $r'\in V(z_iCc)-\{z_i, c\}$ and internally disjoint from
$K\cup X$, and by (6.3), $R$ is disjoint from $J(A,C)$. 
First, assume $e(z_{3-i}, B_1) = 1$. Then $z_{3-i}=p_{3-i}$ and we
use the path $T$ from (A2). Note that $S\cup P\cup b_1Bp$ contains a path $S'$ from $z_i$ to $b_1$.  
Hence,  $G[\{x_1, x_2, y_1, y_2\}] \cup z_{3-i}x_{3-i}\cup z_{3-i}Ty_2\cup 
(z_{3-i}Tb\cup bb_1\cup S'\cup z_ix_i)\cup (z_{3-i}Bq\cup Q)\cup (y_1Cr'\cup R\cup rXp_i\cup p_iTy_2)$ is a $TK_5$ in $G'$
with branch vertices $x_1,x_2,y_1,y_2,z_{3-i}$. So assume $e(z_{3-i}, B_1)\ge 2$. Let $P_1,P_2$ be independent paths  from (A7) with
$q^*=p$. So $P_1\cup P_2\cup B\cup S\cup (P-c)\cup Q$ contains independent
paths $S_1, S_2$ from $z_{3-i}$ to $z_i, y_1$, respectively. Then
$G[\{x_1, x_2, y_1, y_2\}] \cup z_{3-i}x_{3-i}\cup (z_{3-i}Xp_{3-i}\cup Q_{3-i}) \cup (S_1\cup z_ix_i)\cup S_2\cup (y_1Cr'\cup R\cup rXp_i\cup Q_i)$ is a $TK_5$ in $G'$ with branch vertices $x_1,x_2,y_1,y_2,z_{3-i}$. 

So $r=z_i$. By (c) and (d), $G'$ has no path from
$z_iAa^*-\{a^*,z_i\}$ to $(cCy_1-c)\cup (a^*Ay_1-a^*)$ and internally
disjoint from $K$.  Hence, by (5), $\{a^*,c,x_1, x_2, z_i\}$ is a cut in
$G$, and  $i=2$. Let  $F^*:=G[F\cup z_iAa^*\cup z_iCc+\{x_1,x_2\}]$

Suppose $F^*-x_1$ has disjoint paths $S_1, S_2$ from $x_i, z_i$ to $c, a^*$, respectively. 
If $e(z_{3-i}, B_1) = 1$ then $z_{3-i}=p_{3-i}$ and, using the path
$T$ from (A2), we see that 
$G[\{x_1, x_2, y_1, y_2\}] \cup z_{3-i}x_{3-i}\cup z_{3-i}Ty_2\cup
(z_{3-i}Tb\cup bb_1\cup b_1Bp\cup P\cup S_1)\cup (z_{3-i}Bq\cup Q)\cup
(y_1Aa^*\cup S_2\cup z_iXp_i\cup p_iTy_2)$ is a $TK_5$ in $G'$ with
branch vertices $x_1,x_2,y_1,y_2,z_{3-i}$. Now assume $e(z_{3-i},B_1)\ge 2$. 
Let $P_1,P_2$ be independent paths from (A7) with $q^*=p$. Then 
$G[\{x_1, x_2, y_1, y_2\}] \cup z_{3-i}x_{3-i}\cup (z_{3-i}Xp_{3-i}\cup Q_{3-i})\cup (P_1\cup Q)\cup 
(P_2\cup P\cup S_1)\cup (y_1Aa^*\cup S_2\cup z_iXp_i\cup Q_i)$ 
is a $TK_5$ in $G'$ with branch vertices $x_1,x_2,y_1,y_2,z_{3-i}$. 

Thus, we may assume that such $S_1, S_2$ do not exist. Then by Lemma~\ref{2path}, $(F^*-x_1, x_i, z_i, c, a^*)$ is planar. If 
$|V(F^*)|\geq 7$, then the assertion of Theorem~\ref{x1a} follows from
Lemma~\ref{apex}. So assume $|V(F^*)|=6$. 
Let $z\in V(F^*-x_1)-\{x_i, z_i, c, a^*\}$. 
Then $G[\{x_i, z_i, z, c\}]\cong K_4^-$, and $(ii)$ of Theorem~\ref{x1a} holds (as $i=2$ in this case).

\medskip

{\it Subcase} 1.2.  $F=\emptyset$.

Then $L(C)=\emptyset$ by (7). Also,  $L(A)=\emptyset$ by (7); hence,
by (4) and the comment preceding (5), $W=x_iw$ with $w\in V(A)-\{z_i,y_1\}$. 

We may assume that $J(A, C)\cap (A-\{z_i, y_1\})=\emptyset$. For, otherwise, let $t\in V(J(A, C))\cap V(A-\{z_i, y_1\})$. By (6.1), 
$J(A, C)$ contains a path $T$ from $t$ to $t'\in V(Q-y_1)$ and internally disjoint from $K$, and 
$T$ must be internally disjoint from $S$. Note that  $(S\cup P\cup b_1Bp)-c$ contains a path $S'$ from $z_i$ to $b_1$ and
 internally disjoint from $T\cup Q\cup z_{3-i}Bq$. If $e(z_{3-i},B_1)=
 1$ then $z_{3-i}=p_{3-i}$ and, using the path  $T$ from (A2), we see
 that $G[\{x_1, x_2, y_2\}] \cup z_{3-i}x_{3-i}\cup z_{3-i}Ty_2\cup
 z_ix_i\cup (z_iXp_i\cup p_iTy_2)\cup (z_{3-i}Tb\cup bb_1\cup S')
\cup (C\cup y_1x_{3-i})\cup (z_{3-i}Bq\cup qQt'\cup T\cup tAw\cup wx_i)$ is a $TK_5$ in $G'$
with branch vertices $x_1,x_2,y_2,z_1,z_2$. So assume $e(z_{3-i}, B_1)\ge 2$. Let $P_1,P_2$ be independent paths from (A7) with $q^*=p$. Then $P_1\cup P_2\cup B \cup S\cup (P-c)\cup (Q-y_1)\cup T$ contains independent paths $S_1, S_2$ from $z_{3-i}$ to $z_i, t$, respectively. Now
$G[\{x_1, x_2, y_2\}] \cup z_{3-i}x_{3-i}\cup (z_{3-i}Xp_{3-i}\cup Q_{3-i})\cup z_ix_i\cup (z_iXp_i\cup Q_i)\cup S_1\cup (C\cup y_1x_{3-i})\cup 
(S_2\cup tAw\cup wx_i)$ is a $TK_5$ in $G'$ with branch vertices $x_1,x_2,y_2,z_1,z_2$. 

By (A5), $J:=J(A, C)\cup C$ contains no disjoint paths from $z_i,y_1$ to $z_{3-i},b_1$, respectively. 
Hence by Lemma~\ref{2path}, there exists a collection ${\cal L}$ of subsets of $V(J)-\{b_1,y_1,z_1,z_2\}$ such that  $(J, {\cal L}, 
z_i,y_1,z_{3-i},b_1)$ is $3$-planar. We choose ${\cal L}$ so that each $L\in {\cal L}$ is minimal and, subject to this, $|{\cal L}|$ is minimal. 

We claim that for each $L\in {\cal L}$, $L\cap V(L(A,C))=\emptyset$. 
For suppose there exists $L\in {\cal L}$ such that $L\cap V(L(A,C))\ne
\emptyset$. 
Then $|N_J(L)\cap V(C)|\ge 2$. Assume for
the moment that $N_J(L)\subseteq V(C)$. 
Then, since $L(C)=\emptyset$ and   $J(A, C)\cap (A-\{z_i,
y_1\})=\emptyset$, $L\subseteq V(C)$. However, 
since $C$ is an induced path in $G$, we see that $(J, {\cal L}-\{L\}, 
z_i,y_1,z_{3-i}, b_1)$ is 3-planar, contradicting the choice of ${\cal L}$. Thus, let $N_J(L)=\{t_1,t_2,t_3\}$ such that $t_1,t_2\in V(C)$ and $t_3\notin 
V(C)$. Then $J(A,C)$ contains a path $R$ from $t_3$ to $B$ and internally disjoint from $B\cup C$. Let $t\in L\cap V(L(A,C))$. 
By the minimality of $L$, $G[L+\{t_1,t_2,t_3\}]$ contains disjoint paths $T_1,T_2$ from $t_1,t$ to $t_2,t_3$, respectively. We may choose $T_1$ to be induced, and let $C':=z_iCt_1\cup T_1\cup t_2Cy_1$.
Then $A,B,C'$ satisfy (a), but $J(A,C')\subseteq L(A,C')$ (because of $T_2$), contradicting (2) (as $J(A, C)\cap (A-\{z_i, y_1\})=\emptyset$). 

Because of the existence of $Y,Z$ in (A3), there are disjoint paths $R_1,R_2$ in $L(A,C)$ from $r_1,r_2\in V(A)$ to $r_1',r_2'\in V(C)$ such that 
$z_i,r_1,r_2,y_1$ occur on $A$ in order and $z_i,r_2',r_1',y_1$ 
occur on $C$ in order. Let $A'=z_iAr_1\cup R_1\cup r_1'Cy_1$ and $C'=z_iCr_2'\cup R_2\cup r_2Ay_1$. 
Let $t_1,t_2\in V(C-\{z_i,y_1\}) \cap V(J(A,C))$ with $t_1Ct_2$
maximal, and assume that $z_i,t_1,t_2,y_1$ occur on $C$ in this
order. By the planarity of $(J, z_i,y_1,z_{3-i},b_1)$ and by (6.3), $t_1=c$.

Then either $t_1Ct_2\subseteq z_iCr_2'$ for all choices of $R_1$ and $R_2$,  
or  $t_1Ct_2\subseteq r_1'Cy_1$ for all choices of  $R_1$ and $R_2$;
for otherwise, $J(A',C')\subseteq L(A',C')$, and  $A',B,C'$  contradict the choice of $A,B,C$
in (b). Moreover, since $F=\emptyset$,  $t_1Ct_2\subseteq z_iCr_2'$
for all choices of  $R_1$ and $R_2$. Choose $R_1,R_2$ so that $z_iAr_1$ and $z_iCr_2'$ are minimal.
Since $G$ is 5-connected, $\{r_1,r_2',x_1,y_1\}$ cannot be a cut in
$G$. So by (5),  $G'$ has a path $R$ from $x_2$ to some $v\in V(r_1Ay_1-\{r_1,y_1\})\cup V(r_2'Cy_1-\{r_2',y_1\})$ and internally disjoint from $K$. 

First, assume $i=1$. If $v\in V(r_1Ay_1)-\{r_1,y_1\}$ then $G[\{x_1,x_2,y_1,y_2\}] \cup z_ix_i\cup C\cup (z_iXp_i\cup Q_i)\cup 
(zAv\cup R)\cup (Q\cup qBz_{3-i}\cup z_{3-i}Xp_{3-i}\cup Q_{3-i})$ is a $TK_5$ in $G'$ with branch vertices $x_1,x_2,y_1,y_2,z_i$. 
If $v\in  V(r_2'Cy_1)-\{r_2',y_1\}$ then $G[\{x_1,x_2,y_1,y_2\}]
\cup z_ix_i\cup A\cup (z_iXp_i\cup Q_i)\cup 
(z_iCv\cup R)\cup (Q\cup qBz_{3-i}\cup z_{3-i}Xp_{3-i}\cup Q_{3-i})$ is a $TK_5$ in $G'$ with branch vertices $x_1,x_2,y_1,y_2,z_i$. 

Hence, we may assume $i=2$. 
If $e(z_{3-i},B_1) = 1$ then $z_{3-i}=p_{3-i}$ and, using the path
$T$ from (A2), we see that  $G[\{x_1,x_2,y_1,y_2\}] \cup z_{3-i}x_{3-i}\cup z_{3-i}Ty_2\cup (z_{3-i}Bq\cup Q)\cup 
(z_{3-i}Tb_1\cup b_1Bp\cup P\cup cCr_2'\cup R_2
\cup r_2Av\cup R)\cup (y_1Cr_1'\cup R_1\cup r_1Az_i
\cup z_iXp_i\cup p_iTy_2)$ is a $TK_5$ in $G'$ with branch vertices $x_1,x_2,y_1,y_2,z_{3-i}$. 
So assume $e(z_{3-i},B_1)\ge 2$. Let $P_1,P_2$ be independent paths from (A7) with $q^*=p$. 
Now  $G[\{x_1,x_2,y_1,y_2\}] \cup z_{3-i}x_{3-i}\cup (z_{3-i}Xp_{3-i}
\cup Q_{3-i})\cup (P_1\cup Q)\cup (P_2\cup P\cup cCr_2'\cup R_2\cup r_2Av\cup R)\cup (y_1Cr_1'\cup R_1\cup r_1Az_i
\cup z_iXp_i\cup Q_i)$ is a $TK_5$ in $G'$ with branch vertices $x_1,x_2,y_1,y_2,z_{3-i}$.

\medskip

{\it Case} 2.  $z_i\notin V(J(A,C))$.

Then $F\ne \emptyset$ as the degree of $z_i$ in $G'$ is at least 5. So
$a^*$ and $r$ are defined.
\medskip

{\it Subcase} 2.1. $r\ne z_i$,  and $G'$ contains a path $S$ from some $s\in V(z_iXr)-\{z_i, r\}$ to some $s'\in V(P\cup Q\cup B')-\{y_1, c\}$ and internally disjoint from $A\cup B'\cup C\cup P\cup Q\cup X$.

%Note that $s'\in V(B)$ if $s'\in V(B')$. 
First, assume  $s'\in V(Q-y_1)\cup V(pBz_{3-i}-p)$. Then $S\cup
(Q-y_1)\cup (pBz_{3-i}-p)$ has a path $S'$ from $s$ to $z_{3-i}$.  
By (5), let $R$ be a path in $G'$ from $r$ to some $r'\in V(z_iCc)-\{z_i, c\}$ and internally disjoint from 
$A\cup C \cup J(A, C)\cup X$. 
Then $G[\{x_1, x_2, y_1, y_2\}] \cup z_ix_i\cup (z_iXs\cup S'\cup z_{3-i}x_{3-i})\cup A \cup (z_iCr'\cup R \cup
rXp_i\cup Q_i)\cup (y_1Cc\cup P\cup pBb_1\cup b_1b\cup 
Q_3)$ is a $TK_5$ in $G'$ with branch vertices $x_1,x_2,y_1,y_2,z_i$.

Hence, we may assume  $s'\in V(P-c)\cup V(B'-(pBz_{3-i}-p)))$.  Let $P_1,P_2$ be independent paths from (A7)  with $q^*=p$ if $s'\in
P$ and $q^* = s'$ if $s'\in V(B'-(pBz_{3-i}-p))$.
Since $F\ne
\emptyset$ and $B_1:=H-\{z_1,z_2\}$ is 2-connected, $a^*\ne z_i$; so $G'$ has a path $R'$ 
from $r$ to some $r'\in V(z_iAa^*-z_i)$ and internally disjoint from $A\cup cCy_1\cup J(A, C)\cup X$. 

Suppose $e(z_{3-i},B_1) = 1$. Then $z_{3-i}=p_{3-i}$ and we use the
path $T$  from (A2). Note that  $P_1\cup P_2\cup S\cup Q\cup B\cup z_{3-i}Tb\cup bb_1$ 
contains  independent paths $S_1, S_2$ from $z_{3-i}$ to $s', y_1$,
respectively.  So $G[\{x_1, x_2, y_1, y_2\}] \cup z_{3-i}x_{3-i}\cup z_{3-i}Ty_2
\cup (S_1\cup S\cup sXz_i\cup z_ix_i)\cup S_2\cup (y_1Ar'\cup R'\cup rXp_i\cup p_iTy_2)$ 
is a $TK_5$ in $G'$ with branch vertices $x_1,x_2,y_1,y_2,z_{3-i}$. 

Now assume $e(z_{3-i}, B_1)\ge 2$.
 Note that  $P_1\cup P_2\cup B\cup S\cup
P\cup Q$ contains independent paths $S_1, S_2$ from $z_{3-i}$ to $s,
y_1$, respectively. Then $G[\{x_1, x_2, y_1, y_2\}] \cup z_{3-i}x_{3-i}\cup (z_{3-i}Xp_{3-i}\cup Q_{3-i}) \cup S_2\cup (S_1\cup sXz_i\cup z_ix_i)\cup (y_1Ar'\cup R'\cup rXp_i\cup Q_i)$ is a $TK_5$ in $G'$ with branch vertices $x_1,x_2,y_1,y_2,z_{3-i}$.

\medskip 

{\it Subcase} 2.2. $r=z_i$, or $G'$ contains no path from
$z_iXr-\{z_i, r\}$ to $(P\cup Q\cup B')-\{y_1, c\}$ and internally
disjoint from $A\cup B'\cup C\cup P\cup Q\cup X$.

By (c) and (d), $G'$ has no path from $z_iAa^*-\{a*,z_i\}$ to
$(aCy_1-a)\cup (a^*Ay_1-a^*)$ and internally disjoint from $K$. 
Then by (5), (6.2) and (6.3),  $\{a^*, c,r,x_1, x_2\}$ is a cut in
$G$. Hence,  since $G$ is 5-connected,  $i=2$  by (5). Therefore, $G$ has a 5-separation $(G_1,G_2)$ 
such that $V(G_1\cap G_2)=\{a^*, c,r,x_1, x_2\}$ and $G_2=G[F\cup z_2Cc\cup z_2Aa^*\cup x_2Xr+x_1]$. 

Suppose $G_2-x_1$ contains disjoint paths $S_1, S_2$ from $r, x_2$ to $a^*, c$, respectively. 
If $e(z_1,B_1) = 1$ then $z_1=p_1$ and, using the path $T$ from (A2)
with $i=2$, we see that $G[\{x_1, x_2, y_1, y_2\}] \cup z_1x_1\cup
z_1Ty_2\cup (z_1Bq\cup Q)\cup (z_1Tb\cup bb_1\cup b_1Bp\cup P\cup S_2)\cup 
(y_1Aa^*\cup S_1\cup rXp_2\cup p_2Ty_2) $ is a $TK_5$ in $G'$ with branch vertices $x_1,x_2,y_1,y_2,z_1$. 
So assume $e(z_1, B_1)\ge 2$. Let $P_1,P_2$ be independent paths from (A7) with $q^*=p$. Then 
$G[\{x_1, x_2, y_1, y_2\}] \cup z_1x_1\cup (z_1Xp_1\cup Q_1)\cup (P_1\cup Q)\cup (P_2\cup P\cup S_2)\cup 
(y_1Aa^*\cup S_1\cup rXp_2\cup Q_2) $ is a $TK_5$ in $G'$ with branch vertices $x_1,x_2,y_1,y_2,z_1$. 

Thus, we may assume that such $S_1, S_2$ do not exist in $G_2-x_1$. Then by Lemma~\ref{2path}, 
$(G_2-x_1, r, x_2, a^*, c)$ is planar. If $|V(G_2)|\geq 7$ then the
assertion of Theorem~\ref{x1a} follows from Lemma~\ref{apex}.  So assume 
$|V(G_2)|\le 6$.  
If $r=z_2$ and there exists $z\in V(G_2)-\{a^*,c,x_1,x_2,z_2\}$ then $za^*,zc,zx_1, zx_2,zz_2\in E(G)$ (as $G$ is 5-connected); so 
$G[\{c,x_2,z,z_2\}]$ contains $K_4^-$ and $(ii)$ of Theorem~\ref{x1a}
holds. Hence, we may assume that $r\ne z_2$ or $V(G_2)=\{a^*,c,x_1,x_2,z_2\}$.  
Then, $z_2x_1,z_2c\in E(G)$ and $L(C)=\emptyset$ (by (7)).

Recall that  $y_1z_2\notin E(G)$; so  $G[\{x_1, x_2, y_1,z_2\}]\cong
K_4^-$. We complete the proof of Theorem~\ref{x1a} by proving $(iv)$
for this new $K_4^-$.  Let $z_0', z_1'\in N(x_1)-\{x_2,y_1, z_2\}$ be
distinct and let $G'':=G-\{x_1v:v\notin \{x_2,y_1,z_0',z_1', z_2\}\}$. 

Suppose $z_1'\in V(J(A, C))-V(A\cup C)$ or $z_1'\in V(Y_2)$ or
$z_1'\in V(X)$. Then $(J(A,C)\cup Y_2\cup X\cup x_2y_2\cup bb_1)-(A\cup C)$
contains a path from $z_1'$ to $x_2$. Hence, $G-x_1$ contains an
induced path $X'$ from $z_1'$ to $x_2$ such that $A\cup C$ is a cycle
in $(G-x_1)-X'$ and $\{y_1,z_2\}\subseteq V(A\cup C)$. 
So by Lemma~\ref{4A},  we may assume that $X'$ is chosen so that
$y_1,y_2\notin V(X')$ and $(G-x_1)-X'$ is $2$-connected. 
Then by Lemma~\ref{2-connected}, $G''$ contains $TK_5$ (which uses $G[\{x_1,x_2,z_2,y_1\}]$ and $x_1z_1'$). 

So assume $z_1'\in V(L(A, C)-J(A, C))\cup V(A\cup C)$ (as
$L(A)=L(C)=\emptyset$). In fact, $z_1'\in V(C)-\{z_2, y_1\}$. For
otherwise, $(W\cup L(A, C)\cup A)- C$ contains an induced path $X'$
from $z_1'$ to $x_2$, where $W$ comes from (4) and the remark preceding
(5).  Then $(G-x_1)-X'$ contains $C\cup Q\cup qBb_1\cup (X-\{x_1,
x_2\})\cup Y_2$, which has  a cycle containing $\{y_1, z_2\}$. By
Lemma~\ref{4A}, we may assume that $X'$ is chosen so that
$y_1,y_2\notin V(X')$ and $(G-x_1)-X'$ is $2$-connected.  Now the
assertion of Theorem~\ref{x1a} follows from Lemma~\ref{2-connected}. 

If $z_1'\in V(J(A, C))$, then there is a path $P'$ in $J(A, C)$ from
$z_1'$ to some $p'\in V(B)$ and internally disjoint from $A\cup B \cup
C$. So $G[\{x_1, x_2, y_1, z_2\}]\cup z_1'x_1\cup z_1'Cz_2\cup
z_1'Cy_1 \cup (P'\cup p'Bb_1\cup b_1b\cup  Q_3\cup y_2x_2)\cup A$ is a $TK_5$ in $G''$ with 
branch vertices $x_1,x_2,y_1,z_2,z_1'$. 

Thus, we may assume that $z_1'\notin V(J(A, C))$.  So there is a path
$A'$ in $L(A,C)$ from $z_1'$ to some $a'\in V(A)$ 
and internally disjoint from $J(A, C)\cup A\cup C$. Recall the path $W$  from (4) and the remark preceding
(5).  Now $G[\{x_1, x_2, y_1, z_2\}] \cup z_1'x_1\cup z_1'Cz_2\cup z_1'Cy_1\cup (A'\cup a'Aw\cup W)
\cup (Q\cup qBb_1\cup b_1b\cup Q_3\cup Q_2\cup p_2Xz_2)$ is a $TK_5$ in $G''$ 
with branch vertices $x_1,x_2,y_1,z_2,z_1'$.   \qed


\newpage

\begin{thebibliography}{99}
\addtolength{\baselineskip}{-1ex}

%\bibitem {Ca79} P.~Catlin, Haj\'{o}s' graph-coloring conjecture: variations
  %and counterexamples, {\it J.~Combin.~Theory, Ser. B} {\bf 26} (1979) 268--274.

\bibitem {CR79} K. Chakravarti and N. Robertson, Covering three edges with a bond in
a nonseparable graph. {\it Annals of Discrete Math.} (Deza and Rosenberg
eds) (1979) 247.

%\bibitem{CGY} G. Chen, R. Gould, and X. Yu, Graph connectivity
%after path removal, {\it Combinatorica} {\bf 23} (2003) 185--203.

%\bibitem{CY03} S. Curran and X. Yu, Non-separating cycles in 4-connected
%graphs, {\it SIAM J. Discrete Math.} {\bf 16} (2003) 616--629.

%\bibitem {Di06} R. Diestel, Graph Theory (3rd edition), Graduate Text in Mathematics 173, Springer, 2006. 

%\bibitem {EY09} M. N. Ellingham. M. D. Plummer  and G. Yu, Linkage for the diamond and the path with  four vertices. Manuscript. 2010.

%\bibitem {Di52} G. A.~Dirac, A property of 4-chromatic graphs and some
%  remarks on critical graphs, {\it J.~London~Math.~Soc., Ser. B} {\bf 27}
%  (1952) 85-92.

%\bibitem {EF81} P.~Erd\"{o}s and S.~Fajtlowicz, On the conjecture of
%  Haj\'{o}s, {\it Combinatorica} {\bf 1} (1981) 141--143.

%\bibitem {Ha43} H. Hadwiger, \"{U}ber eine Klassifikation der
%  Streckencomplexe, {\it Vierteljahrsschr. Naturforsch. Ges. Z\"{u}rich}
%  {\bf 88} (1943) 133--142. 


\bibitem {HWY15I} D. He, Y. Wang and X. Yu, The Kelmans-Seymour conjecture I: special separations, {\it Submitted}. 

\bibitem {HWY15II} D. He, Y. Wang and X. Yu, The Kelmans-Seymour conjecture II: 2-vertices in $K_4^-$, {\it Submitted}.

%\bibitem {HWY15III} D. He, Y. Wang and X. Yu, The Kelmans-Seymour conjecture III: a strucure,  {\it Submitted}.

%\bibitem {Ka01} K. Kawarabayashi, Note on $k$-contractible edges in  $k$-connected graphs, {\it Australas. J. Combin.} {\bf 24} (2001) 165--168. 

\bibitem {Ka02} K. Kawarabayashi, Contractible edges and triangles in
  $k$-connected graphs, {\it J. Combin. Theory Ser. B} {\bf 85} (2002)
  207--221. 

%\bibitem {KLY05} K. Kawarabayashi, O. Lee, and X. Yu, Non-separating paths
%  in 4-connected graphs, {\it Annals of Combinatorics} {\bf 9} (2005)
%  47--56.

\bibitem {KMY15} K. Kawarabayashi, J. Ma and X. Yu, $K_5$-Subdivisions in graphs containing $K_{2,3}$, 
{\it J. Combin. Theory, Ser. B} {\bf 113} (2015) 18--67.

\bibitem {Ke79} A. K. Kelmans, Every minimal counterexample to the Dirac conjecture
is 5-connected, {\it Lectures to the Moscow Seminar on Discrete Mathematics}
(1979). 


%\bibitem {KM91} A. E.~K\'{e}zdy and P. J.~McGuiness,  Do $3n-5$ edges suffice for a subdivision of $K_5$?
%{\it J. Graph Theory} {\bf 15} (1991) 389-406.
  
% \bibitem{Kr98}  M. Kriesell, Induced paths in $5$-connected graphs,
%{\it J. Graph Theory} {\bf 36} (2001) 52--58.

%\bibitem {KO02} D. K\"{u}hn and D. Osthus, Topological minors in graphs of large
%  girth, {\it J. Combin. Theory Ser. B}  {\bf 86}  (2002) 364--380. 

%\bibitem{Lo75} L. Lov\'{a}sz, Problems, in {\it Recent Advances in
%Graph Theory}, M. Fiedler, ed., Academia, Prague, 1975.

\bibitem{MY10} J. Ma and X. Yu, Independent paths and $K_5$-subdivisions,
{\it J. Combin. Theory Ser. B} {\bf 100} (2010) 600--616.

\bibitem{MY13} J. Ma and X. Yu, $K_5$-subdivisions in graphs containing $K_4^-$, {\it J. Combin. Theory Ser. B} {\bf 103} (2013) 713--732.

%\bibitem {Ma98} W.~Mader, $3n-5$ Edges do force a subdivision of $K_5$,
%{\it Combinatorica} {\bf 18} (1998) 569-595.

\bibitem {Pe68} H. Perfect, Applications of Menger's graph theorem, {\it
    J. Math. Analysis and Applications} {\bf 22}  (1968) 96--111. 

\bibitem {RS90}  Robertson and P. S. Seymour, Graph Minors. IX. Disjoint crossed paths,  {\it J. Combin. Theory Ser. B} {\bf 49} (1990) 40--77.
%\bibitem {RST93} N. Robertson, P. D. Seymour, and R. Thomas, Hadwiger's
%conjecture for $K_6$-free graphs, {\it Combinatorica} {\bf 13}
%(1993) 279--361.

\bibitem{Se77} P. D.~Seymour, Private Communication with X.~Yu.

\bibitem {Se80} P. D. Seymour, Disjoint paths in graphs. {\it Discrete Math.} {\bf
29} (1980) 293--309.

\bibitem{Sh80} Y.~Shiloach,  A polynomial solution to the
undirected two paths problem, {\it J. Assoc. Comp. Mach.} {\bf 27} (1980)
445--456.

\bibitem {Th80}  C. Thomassen, 2-Linked graphs. {\it Europ. J. Combinatorics} {\bf 1}
(1980) 371--378.
 
%\bibitem{Tu63} W. T. Tutte, How to draw a graph, {\it Proc. London Math.
%Soc.} (3) {\bf 13} (1963) 743--767.

%\bibitem{Wa37} K. Wagner, \"{U}ber eine Eigenschaft der ebenen Komplexe,
% {\it Math. Ann.} {\bf 114} (1937) 570--590.

\bibitem {WM67} M. E.~Watkins and D. M.~Mesner, Cycles and connectivity in  graphs, {\it Canadian J.\ Math.} {\bf 19} (1967) 1319-1328.

%\bibitem {Yu98} X. Yu, Subdivisions in planar graphs. {\it
%J. Combinat. Theory Ser. B} {\bf 72} (1998) 10--52. 

%\bibitem {Yu03I} X. Yu, Disjoint paths in graphs I, 3-planar graphs and basic obstructions, {\it 
%Annals of Combinatorics} {\bf 7} (2003) 89--103.

%\bibitem {Yu04II} X. Yu, Disjoint paths in graphs II, a special case, {\it 
%Annals of Combinatorics} {\bf 7} (2003) 105--126.

%\bibitem {Yu04III} X. Yu, Disjoint paths in graphs III, characterization, {\it 
%Annals of Combinatorics} {\bf 7} (2003) 229--246.

%\bibitem {YZ06} X. Yu and F. Zickfeld, 
%Reducing Hajos' coloring conjecture to 4-connected graphs, {\it J. Combin. Theory Ser. B} 
%{\bf 96} (2006) 482--492. 

\end{thebibliography}
\end{document}